\RequirePackage{ifpdf}
\ifpdf 
\documentclass[pdftex]{sigma}
\else
\documentclass{sigma}
\fi

\usepackage{array}
\usepackage{bm}
\begin{document}

 \numberwithin{equation}{section}

 \numberwithin{theorem}{section}
\numberwithin{conjecture}{section}
\numberwithin{proposition}{section}
\numberwithin{lemma}{section}
\numberwithin{corollary}{section}
\numberwithin{definition}{section}
\numberwithin{example}{section}
\numberwithin{remark}{section}
\numberwithin{note}{section}

\allowdisplaybreaks

\renewcommand{\thefootnote}{$\star$}

\renewcommand{\PaperNumber}{058}

\FirstPageHeading

 \def\R{\mathbb{R}}
 \def\Z{\mathbb{Z}}
 \def\gothic#1{\mathfrak{#1}}
 \def\g{\mathfrak{g}}

 \newcommand\parder[1]{\frac{\partial}{\partial #1}}
 \newcommand\Parder[2]{\frac{\partial #1}{\partial #2}}
 \newcommand\Pder[2]{\Parder{#1}{\omega^{#2}}}
 \def\Lieder#1{{\cal L}_{#1}}

 \newcommand\ra{\rightarrow}
 \newcommand\lra{\longrightarrow}
 \newcommand\mat[2]{\left( \begin{array}{#1} #2 \end{array} \right)}
 \def\qRa{\quad\Rightarrow\quad}
 \newcommand\qbox[1]{\quad\mbox{#1}\quad}
 \newcommand\gam[2]{\gamma^{#1}{}_{#2}}
 \newcommand\hgam[2]{\hat\gamma^{#1}{}_{#2}}
 \newcommand\Char[1]{{\cal A}(#1)}
 \newcommand\contact[1]{{\cal C}^{(#1)}}
 \newcommand\intprod{\,\lrcorner\,}
 \newcommand\contactmaps[1]{{\rm Contact}^{#1}(\Sigma_7,\bar\Sigma_7)}
 \newcommand\selfcontactmaps[1]{{\rm Contact}^{#1}(\Sigma_7,\Sigma_7)}
 \newcommand{\hooklongrightarrow}{\lhook\joinrel\longrightarrow}
 \newcommand\class{{\rm class}}
 \def\MA/{Monge--Amp\`ere}

\ShortArticleName{Contact Geometry of Hyperbolic Equations of Generic Type}

\ArticleName{Contact Geometry of Hyperbolic Equations\\ of Generic Type\footnote{This paper is a
contribution to the Special Issue ``\'Elie Cartan and Dif\/ferential Geometry''. The
full collection is available at
\href{http://www.emis.de/journals/SIGMA/Cartan.html}{http://www.emis.de/journals/SIGMA/Cartan.html}}}

 \Author{Dennis THE}

 \AuthorNameForHeading{D. The}

 \Address{McGill University, 805 Sherbrooke Street West, Montreal, QC, H3A 2K6, Canada}
 \Email{\href{mailto:dthe@math.mcgill.ca}{dthe@math.mcgill.ca}}

\ArticleDates{Received April 10, 2008, in f\/inal form August 11,
2008; Published online August 19, 2008}

 \Abstract{We study the contact geometry of scalar second order hyperbolic equations in the plane of generic type.  Following a derivation of parametrized contact-invariants to distinguish Monge--Amp\`ere (class 6-6), Goursat (class 6-7) and generic (class 7-7) hyperbolic equations, we use Cartan's equivalence method to study the generic case.  An intriguing feature of this class of equations is that every generic hyperbolic equation admits at most a~nine-dimensional contact symmetry algebra.  The nine-dimensional bound is sharp: normal forms for the contact-equivalence classes of these maximally symmetric generic hyperbolic equations are derived and explicit symmetry algebras are presented.  Moreover, these maximally symmetric equations are Darboux integrable.  An enumeration of several submaximally symmetric (eight and seven-dimensional) generic hyperbolic structures is also given.}


 \Classification{35A30; 35L70; 58J70}

 \section{Introduction}

 The purpose of this paper is to give a reasonably self-contained account of some key geometric features of a class of (nonlinear) scalar second order hyperbolic partial dif\/ferential equations (PDE) in the plane (i.e.\ in two independent variables) that has received surprisingly very little attention in the literature, namely hyperbolic PDE of {\em generic} (also known as {\em class 7-7}) type.  Even this terminology is not well-known, and so it deserves some clarif\/ication.

 In the geometric study of dif\/ferential equations, there is a natural notion of equivalence associated with the pseudo-group of local point transformations, i.e.\ local dif\/feomorphisms which mix the independent and dependent variables.  Another natural (but coarser) notion is to def\/ine equivalence up to the larger pseudo-group of local contact transformations and one of the principal goals of the geometric theory is to f\/ind invariants to distinguish dif\/ferent contact-equivalence classes.  Restricting now (and for the remainder of this paper) to scalar second order PDE in the plane, we have that given certain nondegeneracy conditions (i.e.\ one can locally solve the equation for one of the highest derivatives), there is a contact-invariant trichotomy into equations of elliptic, parabolic and hyperbolic type.  In the quasi-linear case, invariance of these classes under point transformations appears in \cite{CourantHilbert}.  (Inequivalent normal forms are derived in each case.)  An elegant geometric proof of invariance under contact transformations in the general case was given by R.B. Gardner \cite{Gardner1969-char}.  In the hyperbolic case, there exist two characteristic subsystems which give rise to a f\/iner contact-invariant trichotomy into equations of \MA/ (class 6-6), Goursat (class 6-7), and generic (class 7-7) type.  While this was known to Vranceanu and almost certainly to E.~Cartan and Lie, a modern exposition of these ideas f\/irst appeared in \cite{GK1993}.  To keep our exposition as self-contained as possible, we include the details of these classif\/ications in this paper.  For hyperbolic equations given in the form $z_{xx} = f(x,y,z,z_x,z_y,z_{xy},z_{yy})$,
 the (relative) invariants characterizing the three types of hyperbolic equations were calculated parametrically by Vranceanu (c.f.\ the $B$, $C$ invariants in \cite{Vranceanu1940}).  For a general equation $F(x,y,z,z_x,z_y,z_{xx},z_{xy},z_{yy}) = 0$, these invariants appeared in Chapter~9 of Jur\'a\v{s}' thesis in his characterization of the \MA/ class (c.f.\ $M_\sigma$, $M_\tau$ in \cite{Juras1997}).  Our derivation of these invariants (labelled $I_1$, $I_2$ in this article) is quite dif\/ferent and the novelty in our exposition (see Theorem~\ref{thm:hyp-contact-inv}) is obtaining simpler expressions expressed in terms of certain determinants.  Moreover, we use these invariants to give general examples of hyperbolic equations of Goursat and generic type (see Table~\ref{table:hyp-examples}).

 Hyperbolic equations of \MA/ type have been well-studied in the literature from a geometric point of view (e.g.\ see \cite{Lepage1929, Lychagin1979, Morimoto1979, BGH1-1995, BGH2-1995, Kruglikov1999, Biesecker2003, KLR2007, MVY2007} and references therein).  This class of equations includes the \MA/, wave, Liouville, Klein--Gordon and general $f$-Gordon equations.  At the present time and to the best of our knowledge, there exists only one paper in the literature that has been devoted to the study of hyperbolic equations of generic type.  This paper, {\em ``La g\'eom\'etrisation des \'equations aux d\'eriv\'ees partielles du second ordre''} \cite{Vranceanu1937}, was published by Vranceanu in 1937.  Despite its appearance over 70 years ago, and much attention having been focused on applications of Cartan's equivalence method in the geometric theory of PDE, very few references to \cite{Vranceanu1937} exist.  Currently, the paper does not appear on MathSciNet; the only reference to it by internet search engines appears on Zentralblatt Math.

In \cite{Vranceanu1937}, Vranceanu uses the exterior calculus and Cartan's method of equivalence to study generic hyperbolic equations.  One of the most intriguing results of the paper is that all equations of generic type admit at most a {\em nine}-dimensional local Lie group of (contact) symmetries.  This is in stark contrast to the \MA/ class, where the wave equation is well-known to admit an inf\/inite-dimensional symmetry group.  Vranceanu is able to isolate the correspon\-ding maximally symmetric structure equations as well as some submaximally symmetric structures.  Furthermore, he is able to integrate these abstract structure equations to obtain an explicit parametrization of the corresponding coframe, leading to normal forms for the contact-equivalence classes.  As any practitioner of the Cartan equivalence method can probably attest, this is an impressive computational feat.  Nevertheless, as in the style of Cartan's writings, Vranceanu's arguments are at times dif\/f\/icult to decipher, hypotheses are not clearly stated or are dif\/f\/icult to discern amidst the quite lengthy calculations, and some of his results are not quite correct.  In this paper, we reexamine, clarify, and sharpen some of Vranceanu's results with the perspective of our modern understanding of the geometric theory of dif\/ferential equations through exterior dif\/ferential systems and Cartan's equivalence method.  The hope is that this exposition will provide a clearer understanding of the geometry of this class of equations for a~contemporary audience.

 In Section \ref{background} we recall the contact-invariant classif\/ication of second order scalar PDE into elliptic, parabolic, and hyperbolic classes based on invariants of a (conformal class of a) symmetric bilinear form, and def\/ine the $M_1$ and $M_2$ characteristics in the hyperbolic case.  This leads to a preliminary set of structure equations for hyperbolic equations.  In Section \ref{sec:hyperbolic-subclassify}, the structure equations are further tightened, and using them we show how the class of $M_1$ and $M_2$ leads to the f\/iner classif\/ication into equations of \MA/, Goursat, and generic types.  In Theorem \ref{thm:hyp-contact-inv}, these subclasses of hyperbolic equations are characterized by means of the relative invariants $I_1$, $I_2$. We then restrict to the generic case and derive the generic hyperbolic structure equations.  We note that in Vranceanu's derivation of the generic hyperbolic structure equations (c.f.\ \eqref{StrEqns123} in this paper), the $\epsilon={\rm sgn}(I_1I_2) = \pm 1$ contact invariant was overlooked.  This carries through to the normal forms for the contact-equivalence classes.  Section \ref{str-grp-eq-prb} formulates the equivalence problem for generic hyperbolic equations and recalls some facts from Cartan's theory of $G$-structures applied to our situation.  The structure group that we consider here is strictly larger than Vranceanu's, dif\/fering by certain discrete components.  These naturally arise when considering automorphisms which interchange the $M_1$ and $M_2$ characteristics.  Both Vranceanu and Gardner--Kamran consider only automorphisms which preserve each of $M_1$ and $M_2$.  The nine-dimensional bound on the symmetry group of any generic hyperbolic equation is established in Section \ref{9d-sym}.

 In Section \ref{complete-str-eqs}, we give a clear enumeration of several generic hyperbolic structures which result from Vranceanu's reduction of the structure equations.  These include the aforementioned maximally symmetric (nine-dimensional) structure equations as well as some new submaximally symmetric (eight and seven-dimensional) structures including some with nonconstant torsion.  (Vranceanu derived two eight-dimensional structures with constant torsion in addition to the maximally symmetric structures.)  Finally, Section \ref{maxsym-case} gives a detailed account of the maximally symmetric case.  Integration of the abstract maximally symmetric structure equations leads to the contact-equivalence classes of maximally symmetric generic hyperbolic PDE being parametrized by $(\epsilon, a) \in \{ \pm 1 \} \times (0,1]$, with normal forms given by
 \begin{gather*}
  (\epsilon, a) = (1,1) : \quad 3z_{xx}(z_{yy})^3 + 1 = 0,\\
  (\epsilon, a) \neq (1,1) : \quad (\epsilon + a)^2 \left(2 z_{xy} - (z_{yy})^2 \right)^3 + \epsilon a \left( 3z_{xx} - 6z_{xy}z_{yy} + 2(z_{yy})^3 \right)^2 = 0.
 \end{gather*}
The isomorphism type of the symmetry algebra for the second equation is {\em independent of $(\epsilon,a)$} and is non-isomorphic to the symmetry algebra of the f\/irst equation.  Thus, there are precisely {\em two} non-isomorphic (abstract) symmetry algebras that arise for maximally symmetric generic hyperbolic equations.  These equations are further distinguished in a contact-invariant way using a contact invariant $\Delta_1$ and a relative contact invariant $\Delta_2$.  Both equations satisfy $\Delta_1=0$, but the former satisf\/ies $\Delta_2=0$ while the latter satisf\/ies $\Delta_2 \neq 0$.

 Let us point out two additional points of discrepancy with Vranceanu's calculations: (1) the restriction of the range of the parameter $a$ to $(0,1]$, and (2) a missing factor of 2 for the $z_{xy}z_{yy}$ term in the second equation above.  The f\/irst point is a consequence of the aforementioned larger structure group used in our formulation of the equivalence problem.  The additional discrete factors lead to identif\/ications of dif\/ferent parameter values.  The second point was minor and the error was only introduced by Vranceanu in the last step of his derivation.  To give added justif\/ication to the calculation of the normal forms above, we give explicitly the nine-dimensional symmetry algebras for the normal forms listed above.  Both equations admit the symmetries
  \begin{gather*}
   X_1 = \parder{x}, \qquad X_2 =\parder{y}, \qquad X_3 =\parder{z}, \qquad X_4 =x\parder{z}, \qquad X_5 =y\parder{z},\\
   X_6=x\parder{x} + y\parder{y} + 2z\parder{z}.
 \end{gather*}
 The former admits the additional symmetries
 \begin{gather*}
X_7 = xy\parder{z}, \qquad
   X_8 = 2y\parder{y} + 3z\parder{z}, \qquad
   X_9 = x^2\parder{x} + xz\parder{z},
 \end{gather*}
 while the latter admits the additional symmetries
 \begin{gather*}
  X_7 = y\parder{y} + 3z\parder{z},
\qquad X_8 = x \parder{y} - \frac{1}{2} y^2 \parder{z},
\qquad X_9= x^2 \parder{x} + xy \parder{y} + \left(xz-\frac{1}{6} y^3\right) \parder{z}.
 \end{gather*}
 The calculation of these symmetries (especially in the latter case) is in general a nontrivial task considering the complexity of the equation.

 Numerous appendices provide the details of the proofs of the main statements in the body of this article.

 All considerations in this paper are local, and we will work in the smooth category.  We use the Einstein summation convention throughout.  We will make the convention of using braces enclosing 1-forms to denote the corresponding submodule generated by those 1-forms.   In general, we will abuse notation and not distinguish between a submodule of 1-forms and its corresponding algebraic ideal (i.e.\ with respect to the wedge product) in the exterior algebra.  This is useful when stating structure equations, e.g.\ $d\omega^1 \equiv 0\; \mod I_F$.

 \section{Contact-equivalence of PDE}
 \label{background}

 Consider a scalar second order PDE
 \begin{equation}
 F\left(x,y,z, \frac{\partial z}{\partial x}, \frac{\partial z}{\partial y}, \frac{\partial^2 z}{\partial x^2},  \frac{\partial^2 z}{\partial x \partial y},  \frac{\partial^2 z}{\partial y^2}\right) = 0
 \label{pde}
 \end{equation}
 in two independent variables $x$, $y$ and one dependent variable $z$.   A natural geometrical setting for~\eqref{pde} is the space of 2-jets
$J^2(\R^2,\R)$ with standard local coordinates $(x,y,z,p,q,r,s,t)$ ({\em Monge coordinates}), and the equation above yields a locus
 \begin{gather*}
 L_F = \left\{ (x,y,z,p,q,r,s,t) \in J^2(\R^2,\R): F(x,y,z,p,q,r,s,t) = 0 \right\}.
 \end{gather*}

 We assume that $L_F$ is the image of an open subset $\Sigma_7 \subset \R^7$ under a smooth map $i_F : \Sigma_7 \ra J^2(\R^2,\R)$.

 \begin{definition} We will say that $i_F$ is a {\em nondegenerate parametrization} of the equation $F=0$ if~$i_F$ has maximal rank and $L_F$ is everywhere transverse to the f\/ibers of the natural projection
 \begin{gather*}
 \pi^2_1 : J^2(\R^2,\R) \ra J^1(\R^2,\R),
 \end{gather*}
 i.e.\ ${\rm im}(i_{F*}) + \ker(\pi^2_{1\,*}) = TJ^2(\R^2,\R)$.
 \end{definition}

 We will always work with nondegenerate parametrizations in this paper.  By the transversality assumption $(F_r,F_s,F_t) \neq 0$, and so by the implicit function theorem, one can locally solve \eqref{pde} for one of the highest-order derivatives.
  Since ${\rm im}((\pi^2_1 \circ i_F)_*) = TJ^1(\R^2,\R)$, then $(\pi^2_1 \circ i_F)^* (dx \wedge dy \wedge dz \wedge dp \wedge dq)\neq 0$ and so the standard coordinates $(x,y,z,p,q)$ on $J^1(\R^2,\R)$ along with two additional coordinates $u$, $v$ may be taken as coordinates on $\Sigma_7$.  Thus, without loss of generality, we may assume the parametrization $i_F$ has the form $i_F(x,y,z,p,q,u,v) = (x,y,z,p,q,r,s,t)$, expressing $r$, $s$, $t$ as functions of $(x,y,z,p,q,u,v)$.

 The contact system $\contact{2}$ on $J^2(\R^2,\R)$ is generated by the standard contact forms
 \begin{gather*}
 \theta^1 = dz - pdx - qdy, \qquad \theta^2 = dp - rdx - sdy, \qquad \theta^3 = dq - sdx - tdy
 \end{gather*}
 and pulling back by $i_F$, we obtain a Pfaf\/f\/ian system (i.e.\ generated by 1-forms) $I_F$ on $\Sigma_7$,
 \begin{gather*}
 I_F = i_F^*(\contact{2}) = \{ \omega^1, \omega^2, \omega^3 \},
 \end{gather*}
 where $\omega^\alpha = i_F^* \theta^\alpha$.  There is a correspondence between local solutions of \eqref{pde} and local integral manifolds of $I_F$.

 \begin{definition}
 The equations \eqref{pde} and
 \begin{gather}
  \bar{F}\left(\bar{x},\bar{y},\bar{z}, \frac{\partial \bar{z}}{\partial \bar{x}}, \frac{\partial \bar{z}}{\partial \bar{y}}, \frac{\partial^2 \bar{z}}{\partial \bar{x}^2},  \frac{\partial^2 \bar{z}}{\partial \bar{x} \partial \bar{y}},  \frac{\partial^2 \bar{z}}{\partial \bar{y}^2}\right) = 0,
  \qquad \mbox{(with} \ \  i_{\bar{F}} : \bar\Sigma_7 \ra J^2(\R^2,\R) )\label{pde-Fbar}
 \end{gather}
 are {\em contact-equivalent} if there exists a local dif\/feomorphism $\phi : \Sigma_7 \ra \bar\Sigma_7$ such that $\phi^*I_{\bar{F}} = I_F$.  The collection of all such maps will be denoted $\contactmaps{}$.  A {\em contact symmetry} is a~self-equivalence.
 \end{definition}

 \begin{remark} More precisely, the above def\/inition refers to {\em internal} contact-equivalence.  There is another natural notion of equivalence: namely, \eqref{pde} and \eqref{pde-Fbar} are {\em externally} contact-equivalent if there exists a local dif\/feomorphism $\rho : J^2(\R^2,\R) \ra J^2(\R^2,\R)$ that restricts to a local dif\/feomorphism $\tilde\rho : i_F(\Sigma_7) \ra i_{\bar{F}}(\bar\Sigma_7)$ and preserves the contact system, i.e.\ $\rho^*(\contact{2}) = \contact{2}$.  It is clear that any external equivalence induces a corresponding internal equivalence, but in general the converse need not hold.  The dif\/ference between these two natural notions of equivalence is in general quite subtle and has been investigated in detail in \cite{AKO1993}.  A corollary of their results (c.f.\ Theorem 18 therein) is that for \eqref{pde}, under the maximal rank and transversality conditions, any internal equivalence extends to an external equivalence, and thus the correspondence is one-to-one.
 \end{remark}

 As shown by Gardner \cite{Gardner1969-char}, the (pointwise) classif\/ication of \eqref{pde} into mutually exclusive elliptic, parabolic and hyperbolic classes is in fact a contact-invariant classif\/ication which arises through invariants (i.e.\ rank and index) of a (conformal class of a) symmetric $C^\infty(\Sigma_7)$-bilinear form $\langle \cdot , \cdot \rangle_7$ on $I_F$, namely
 \begin{gather}
 \langle \varphi, \psi \rangle_7 {\rm Vol}_{\Sigma_7} := d\varphi \wedge d\psi \wedge \omega^1 \wedge \omega^2 \wedge \omega^3,
 \label{bilinear-form-defn}
 \end{gather}
 where ${\rm Vol}_{\Sigma_7}$ denotes any volume form on $\Sigma_7$.
 Since $i_F^*$ is surjective, there exists a 7-form $\nu$ on $J^2(\R^2,\R)$ such that $i_F^* \nu = {\rm Vol}_{\Sigma_7}$, and so
 \begin{gather*}
 \langle \varphi, \psi \rangle_7 i_F^* \nu = i_F^*(d\tilde\varphi \wedge d\tilde\psi \wedge \theta^1 \wedge \theta^2 \wedge \theta^3 ),
 \end{gather*}
 where $\tilde\varphi$ and $\tilde\psi$ are any forms on $J^2(\R^2,\R)$ such that $\varphi = i_F^* \tilde\varphi$ and $\psi = i_F^* \tilde\psi$.
 Since $i_{F*} : T\Sigma_7 \ra TJ^2(\R^2,\R)$ is rank 7 (as is $i_F^* : T^*J^2(\R^2,\R) \ra T^*\Sigma_7$) and $i_F^*dF = 0$, then $\ker(i_F^*) = \{ dF \}$, and
 \begin{gather}
 i_F^* \eta = 0 \qquad \mbox{if\/f} \qquad \eta \wedge dF = 0, \qquad \forall \; \eta \in \Omega^*(J^2(\R^2,\R)).
 \label{dF-lemma}
 \end{gather}
 Consequently, letting ${\rm Vol}_{J^2(\R^2,\R)} = \nu \wedge dF$, we see that \eqref{bilinear-form-defn} is equivalent to
 \begin{gather*}
 (\langle \varphi, \psi \rangle_7)_p ({\rm Vol}_{J^2(\R^2,\R)} )_{i_F(p)} = (d\tilde\varphi \wedge d\tilde\psi \wedge \theta^1 \wedge \theta^2 \wedge \theta^3 \wedge dF)_{i_F(p)},
 \end{gather*}
 where $p \in \Sigma_7$.
 This def\/inition is well-def\/ined: it is independent of the choice of $\tilde\varphi$ and $\tilde\psi$ so long as $\varphi = i_F^* \tilde\varphi$ and $\psi = i_F^* \tilde\psi$.

 A computation in the basis $\omega^1$, $\omega^2$, $\omega^3$ reveals that a volume form may be chosen so that
 \begin{gather}
 (\langle \omega^\alpha, \omega^\beta \rangle_7)_p = \left( \begin{array}{ccc} 0 & 0 & 0\\ 0 & F_t & -\frac{1}{2} F_s\\ 0 & -\frac{1}{2} F_s & F_r \end{array} \right)_{i_F(p)}. \label{bilinear-form-matrix}
 \end{gather}
 Our assumption that $i_F$ have maximal rank implies that $\langle \cdot , \cdot \rangle_7$ cannot have rank zero.
 Def\/ining
 \begin{gather*}
 \Delta = i_F^*\left(F_r F_t - \frac{1}{4} F_s{}^2\right),
 \end{gather*}
 we have the following mutually exclusive cases at each point $p \in \Sigma_7$:

 \begin{table}[h] \centering
 \caption{Contact-invariant classif\/ication of scalar second order PDE in the plane.}
 \vspace{1mm}

 $\begin{array}{|c|c|c|} \hline
 \mbox{elliptic} & \mbox{parabolic} & \mbox{hyperbolic}\\ \hline \hline
 \Delta(p) > 0 & \Delta(p) = 0 & \Delta(p) < 0  \\ \hline
 \end{array}$
 \end{table}

 By the commutativity of pullbacks with $d$, it is clear that this classif\/ication is a priori contact-invariant.   We remark that in the classical literature on the geometry of hyperbolic equations, the terminology {\em Monge cha\-rac\-teristics} appears.  These are determined by the roots of the {\em cha\-rac\-teristic equation}
 \begin{gather}
 \lambda^2 - F_s \lambda + F_t F_r = 0. \label{char-eqn}
 \end{gather}
 The discriminant of this equation (with the coef\/f\/icients evaluated on $F=0$) is precisely $-\frac{1}{4} \Delta$, and so the elliptic, parabolic, and hyperbolic cases correspond to the existence of no roots, a~double root, and two distinct roots respectively.

 In the analysis to follow, all constructions for a PDE $F=0$ will implicitly be repeated for a~second PDE $\bar{F}=0$ (if present).  We will concern ourselves exclusively with the hyperbolic case, that is, an open subset of $\Sigma_7$ on which $F=0$ is hyperbolic.
 By the hyperbolicity condition, the two nonzero eigenvalues of $\langle \cdot, \cdot \rangle_7$ have opposite sign, and hence there exists a pair of rank two maximally isotropic subsystems $M_1$ and $M_2$ of $I_F$ at every point of consideration.

 \begin{definition} Given hyperbolic PDE $F=0$ and $\bar{F}=0$, def\/ine
 \begin{gather*}
   \contactmaps{+}  = \{ \phi \in \contactmaps{} : ~\phi^* \bar{M}_1 = M_1, ~\phi^*\bar{M}_2 = M_2 \},\\
   \contactmaps{-}  = \{ \phi \in \contactmaps{} : ~\phi^* \bar{M}_1 = M_2, ~\phi^*\bar{M}_2 = M_1 \}.
 \end{gather*}
 If $\bar{F} = F$, we take $\bar\Sigma_7=\Sigma_7$ and use the notation ${\rm Aut}(I_F) := \selfcontactmaps{}$, etc.
 \end{definition}

 \begin{remark}
 Implicitly, given the Pfaf\/f\/ian system $I_F$ corresponding to a hyperbolic PDE $F=0$, we assume that a choice of labelling for the $M_1$ and $M_2$ characteristics has been made.  This is of course not intrinsic.  All of our f\/inal results will not depend on this choice.
 \end{remark}

 Both Vranceanu \cite{Vranceanu1937} and Gardner--Kamran \cite{GK1993} consider only local dif\/feomorphisms which preserve each of $M_1$ and $M_2$.

 \begin{example}
 For the wave equation written as $z_{xy}=0$, we have the pullbacks of the contact forms on $J^2(\R^2,\R)$ to the parameter space $\Sigma_7 : (x,y,z,p,q,r,t)$,
 \begin{gather*}
 \omega^1 = dz - pdx - qdy, \qquad
 \omega^2 = dp - r dx, \qquad
 \omega^3 = dq - t dy
 \end{gather*}
 and
 \begin{gather*}
 (\langle \omega^\alpha, \omega^\beta \rangle_7)_p = \left( \begin{array}{ccc} 0 & 0 & 0\\ 0 & 0 & -\frac{1}{2} \\ 0 & -\frac{1}{2} & 0 \end{array} \right).
 \end{gather*}
 Thus, $M_1 = \{ \omega^1, \omega^2 \}$ and $M_2 = \{ \omega^1, \omega^3 \}$.  Interchanging the independent variables induces $\phi_0 : \Sigma_7 \ra \Sigma_7$, $(x,y,z,p,q,r,t) \mapsto (y,x,z,q,p,t,r)$, which satisf\/ies
 \begin{gather*}
 \phi_0^*\omega^1 = \omega^1, \qquad
 \phi_0^*\omega^2 = \omega^3, \qquad
 \phi_0^*\omega^3 = \omega^2,
 \end{gather*}
 and hence $\phi_0 \in {\rm Aut}^-(I_F)$.
 \end{example}

 The hyperbolicity condition implies that there exists a local basis of $I_F$ which by abuse of notation we also denote $\omega^1$, $\omega^2$, $\omega^3$ such that
 \begin{gather*}
 M_1 = \{ \omega^1, \omega^2 \}, \qquad
 M_2 = \{ \omega^1, \omega^3 \}
 \end{gather*}
 and the matrix representing $\langle \cdot , \cdot \rangle_7$ is in Witt normal form
 \begin{gather*}
 (\langle \omega^\alpha, \omega^\beta \rangle_7)_p = \left( \begin{array}{ccc} 0 & 0 & 0\\ 0 & 0 & 1\\ 0 & 1 & 0 \end{array} \right).
 \end{gather*}

 \begin{lemma}[{\bf Preliminary hyperbolic structure equations}]
 There exists a (local) coframe $\bm\omega = \{ \omega^i \}_{i=1}^7$ on $\Sigma_7$ such that $I_F = \{ \omega^1, \omega^2, \omega^3 \}$ and
 \begin{gather}
     d\omega^1 \equiv 0, \nonumber\\
     d\omega^2  \equiv  \omega^4 \wedge \omega^5, \qquad\mod I_F,  \label{hyp-str-eqns}\\
     d\omega^3  \equiv  \omega^6 \wedge \omega^7,
\nonumber
 \end{gather}
 with
 \begin{gather*}
 \omega^1 \wedge \omega^2 \wedge \omega^3 \wedge \omega^4 \wedge \omega^5 \wedge \omega^6 \wedge \omega^7 \neq 0.
 \end{gather*}
 \end{lemma}

 \begin{proof}
 In Theorem 1.7 in \cite{BCGGG1991}, an algebraic normal form for a 2-form $\Omega$ is given.  In particular, if $\Omega \wedge \Omega =0$, then $\Omega = \sigma^1 \wedge \sigma^2$ is decomposable.  This statement is also true in a relative sense: if $\Omega \wedge \Omega \equiv 0 \mod I$, then $\Omega \equiv \sigma^1 \wedge \sigma^2 \mod I$, where $I$ is an ideal in the exterior algebra.

 Using this fact, let us deduce consequences of the Witt normal form. By def\/inition of $\langle \cdot, \cdot \rangle_7$, we have (taking congruences below modulo $I_F$)
 \begin{align*}
 &\langle \omega^2, \omega^2 \rangle_7 = 0 \quad\Leftrightarrow\quad d\omega^2 \wedge d\omega^2 \equiv 0 \quad\Leftrightarrow\quad d\omega^2 \equiv \omega^4 \wedge \omega^5, \\
 &\langle \omega^3, \omega^3 \rangle_7 = 0 \quad\Leftrightarrow\quad d\omega^3 \wedge d\omega^3 \equiv 0 \quad\Leftrightarrow\quad d\omega^3 \equiv \omega^6 \wedge \omega^7, \\
 & \langle \omega^2, \omega^3 \rangle_7 = 1 \quad\Leftrightarrow\quad
 d\omega^2 \wedge d\omega^3 \wedge \omega^1 \wedge \omega^2 \wedge \omega^3 = \omega^4 \wedge \omega^5 \wedge \omega^6 \wedge \omega^7 \wedge \omega^1 \wedge \omega^2 \wedge \omega^3 \neq 0.
 \end{align*}
 Using $\langle \omega^1, \omega^2 \rangle_7 = \langle \omega^1, \omega^3 \rangle_7 = 0$, we have
 \begin{align*}
 0 &= d\omega^1 \wedge d\omega^2 \wedge \omega^1 \wedge \omega^2 \wedge \omega^3 = d\omega^1 \wedge \omega^4 \wedge \omega^5 \wedge \omega^1 \wedge \omega^2 \wedge \omega^3, \\
 0 &= d\omega^1 \wedge d\omega^3 \wedge \omega^1 \wedge \omega^2 \wedge \omega^3 = d\omega^1 \wedge \omega^6 \wedge \omega^7 \wedge \omega^1 \wedge \omega^2 \wedge \omega^3,
 \end{align*}
 and thus $d\omega^1 \equiv 0$.
 \end{proof}

 Consequently, $\{ \omega^i \}_{i=1}^7$ is a (local) coframe on $\Sigma_7$, and the structure equations can be written
 \begin{gather}
 d\omega^i = \frac{1}{2} \gamma^i{}_{jk} \omega^j \wedge \omega^k.
 \label{gamma-defn}
 \end{gather}

 \section{\MA/, Goursat and generic hyperbolic equations}
 \label{sec:hyperbolic-subclassify}
 The congruences appearing in the preliminary hyperbolic structure equations can be tightened with a more careful study of integrability conditions and further adaptations of the coframe.  The details are provided in Appendix \ref{app-A}.

 \begin{theorem}[{\bf Hyperbolic structure equations}]  \label{general-hyp-str-eqns}
    Given any hyperbolic equation $F=0$ with nondegenerate parametrization $i_F : \Sigma_7 \ra J^2(\R^2,\R)$, there is an associated coframe $\bm\omega = \{ \omega^i \}_{i=1}^7$ on $\Sigma_7$ such that
    \begin{enumerate}\itemsep=0pt
    \item $I_F = \{ \omega^1, \omega^2, \omega^3 \}$, \ $
            M_1 = \{ \omega^1, \omega^2 \}$, \ $
            M_2 = \{ \omega^1, \omega^3 \}$.
    \item We have the structure equations
     \begin{alignat}{3}
     & d\omega^1  \equiv  \omega^3 \wedge \omega^6 + \omega^2 \wedge \omega^4 \quad && \mod \{ \omega^1 \},&\nonumber \\
    & d\omega^2 \equiv \omega^4 \wedge \omega^5 + U_1 \omega^3 \wedge \omega^7 \quad && \mod \{ \omega^1,\omega^2  \}, & \label{U1U2-streq}\\
   &  d\omega^3 \equiv \omega^6 \wedge \omega^7 + U_2 \omega^2 \wedge \omega^5 \quad && \mod \{ \omega^1,\omega^3 \}&\nonumber
     \end{alignat}
   for some functions $U_1$, $U_2$ on $\Sigma_7$.
  \end{enumerate}
 \end{theorem}

 A f\/iner contact-invariant classif\/ication of hyperbolic equations arises from the study of the {\em class} of $M_1$ and $M_2$.  Let us recall some basic def\/initions.

  \begin{definition}
 Let $I$ be a Pfaf\/f\/ian system on a manifold $\Sigma$.  Def\/ine the
 \begin{enumerate}\itemsep=0pt
 \item {\em Cauchy characteristic space} $\Char{I} =  \{ X \in\mathfrak{X}(\Sigma) : X \in I^\perp,~ X \intprod dI \subset I \}$.
 \item {\em Cartan system} ${\cal C}(I) = \Char{I}^\perp$.  The {\em class} of $I$ is the rank of ${\cal C}(I)$ (as a $C^\infty(\Sigma)$-module).
 \end{enumerate}
 Here, $\perp$ refers to the annihilator submodule.
 \end{definition}

 The hyperbolic structure equations indicate that there are only two possibilities for the class of $M_1$ and $M_2$.

 \begin{corollary} \label{class-cor} For $i=1,2$, $\class(M_i) = 6 \mbox{ or } 7$.  Moreover, $\class(M_i)=6$ iff $U_i =0$.
 \end{corollary}

 \begin{proof} Let $\{ \Pder{}{i} \}_{i=1}^7$ denote the dual basis to $\{ \omega^i \}_{i=1}^7$.  From \eqref{U1U2-streq}, we have
 \begin{gather*}
 \Char{M_1} \subset \left\{ \Pder{}{7} \right\}, \qquad \Char{M_2} \subset \left\{ \Pder{}{5} \right\}.
 \end{gather*}
 Moreover, $\class(M_1) = 6$ if\/f $\Pder{}{7} \in \Char{M_1}$ if\/f $U_1=0$.  Similarly for $M_2$.
 \end{proof}

Consequently, we obtain the subclassif\/ication of hyperbolic equations given in Table~\ref{table:hyp-eqns}.
 \begin{table}[h]\centering
 \caption{Contact-invariant classif\/ication of hyperbolic PDE.}
 \label{table:hyp-eqns}
 \vspace{1mm}

 \begin{tabular}{|c|c|} \hline
 Type & Contact-invariant classif\/ication\\ \hline\hline
 \MA/ (6-6)& $\class(M_1)=\class(M_2)=6$\\
 Goursat (6-7) & $\{ \class(M_1), \class(M_2) \} = \{ 6, 7\}$ \\
 generic (7-7) & $\class(M_1)=\class(M_2)=7$ \\ \hline
 \end{tabular}
 \end{table}

 \begin{example}  \label{hyp-examples} We give some known examples of each type of hyperbolic equation:
 \begin{itemize}\itemsep=0pt
 \item \MA/: wave equation $z_{xy}=0$, Liouville equation $z_{xy} = e^z$, Klein--Gordon equation $z_{xy} = z$, or more generally the $f$-Gordon equation $z_{xy} = f(x,y,z,z_x,z_y)$, and \MA/ equation $z_{xx} z_{yy} - (z_{xy})^2 = f(x,y)$.
 \item Goursat: $z_{xx} = f(z_{xy})$ where $f'' \neq 0$.
 \item generic: $z_{xy} = \frac{1}{2} \sin(z_{xx})\cos(z_{yy})$, or $3z_{xx} (z_{yy})^3 + 1=0$.
 \end{itemize}
 \end{example}

 The terminology for class 6-6 equations is justif\/ied by the following result, known to Vran\-cea\-nu~\cite{Vranceanu1940}.  We refer the reader to Gardner--Kamran \cite{GK1993} for a modern proof.

 \begin{theorem} \label{thm:MA}
 A second-order hyperbolic equation has $\class(M_i)=6$, $i=1,2$ if and only if its locus can be given by an equation of the form
 \begin{gather*}
 a ( z_{xx} z_{yy} - (z_{xy})^2 ) + b z_{xx} + 2c z_{xy} + d z_{yy} + e = 0,
 \end{gather*}
 where $a$, $b$, $c$, $d$, $e$ are functions of $x$, $y$, $z$, $z_x$, $z_y$.
 \end{theorem}

 The examples given above were obtained by constructing explicit coframes which realize the abstract structure equations given in Theorem~\ref{general-hyp-str-eqns}, which in general is a very tedious task and is equation-specif\/ic.  We state here two relative invariants $I_1$, $I_2$ (which are related to the two relative invariants $U_1$, $U_2$) whose vanishing/nonvanishing determine the type of any hyperbolic equation.  Given any hyperbolic equation $F=0$, def\/ine
 \begin{gather*}
 \lambda_\pm = \frac{F_s}{2} \pm \sqrt{|\Delta|},
 \end{gather*}
 which are the roots of the characteristic equation~\eqref{char-eqn}.
 Without loss of generality, we may assume that $F_s \geq 0$.  (If not, take $\hat{F} = -F$ instead.)  By the hyperbolicity assumption $\lambda_+ > 0$.  The proof of the following theorem is given in Appendix~\ref{app:hyp-contact-inv}.

 \begin{theorem}[{\bf Relative contact invariants for hyperbolic equations}] \label{thm:hyp-contact-inv}
 Suppose that $F=0$ is a hyperbolic equation with $F_s \geq 0$.  Let
 \begin{gather*}
  \tilde{I}_1 = \det\left( \begin{array}{ccc} F_r & F_s & F_t\\  \lambda_+ & F_t & 0\\
  \left( \frac{F_t}{\lambda_+} \right)_r &
  \left( \frac{F_t}{\lambda_+} \right)_s &
  \left( \frac{F_t}{\lambda_+} \right)_t
  \end{array} \right), \qquad
 \tilde{I}_2 = \det\left( \begin{array}{ccc} 0 & F_r & \lambda_+ \\ F_r & F_s & F_t\\
  \left( \frac{F_r}{\lambda_+} \right)_r &
  \left( \frac{F_r}{\lambda_+} \right)_s &
  \left( \frac{F_r}{\lambda_+} \right)_t
  \end{array} \right),
 \end{gather*}
 and $I_i = i_F^* \tilde{I}_i$.  Then we have the following classification of $F=0$:
 \begin{center}
 \begin{tabular}{|c|c|} \hline
 Type & Contact-invariant classification\\ \hline\hline
 \MA/ & $I_1=I_2=0$\\
 Goursat & exactly one of $I_1$ or $I_2$ is zero\\
 generic & $I_1I_2 \neq 0$ \\ \hline
 \end{tabular}
 \end{center}
 Moreover, we have the scaling property: If $\phi$ is a function on $J^2(\R^2,\R)$ such that $i_F^* \phi > 0$, then
 \begin{gather*}
 \hat{F} = \phi F \qRa \hat{I}_i = (i_F^*\phi)^2 I_i, \quad i=1,2.
 \end{gather*}
 \end{theorem}

 We note that the scaling property is a fundamental property of these relative invariants: their vanishing/nonvanishing depends only on the equation locus.

 \begin{remark} For a general hyperbolic equation $F(x,y,z,p,q,r,s,t) = 0$, Jur\'a\v{s}' \cite{Juras1997} calculated two (relative) invariants $M_\sigma$, $M_\tau$ whose vanishing characterizes the \MA/ class.  His invariants were given explicitly in terms of two non-proportional real roots $(\mu,\lambda) = (m_x,m_y)$ and $(\mu,\lambda) = (n_x,n_y)$ of the characteristic equation
 \begin{gather*}
 F_r \lambda^2 - F_s \lambda\mu + F_t \mu^2 = 0,
 \end{gather*}
 which he associates to the given PDE.  We note here that the characteristic equation \eqref{char-eqn} that we have used dif\/fers from that of Jur\'a\v{s} (but has the same discriminant).  Our invariants $I_1$, $I_2$ appear to be simpler written in this determinantal form.
 \end{remark}

 Using the relative contact invariants $I_1$, $I_2$ we can identify some more general examples of hyperbolic equations of Goursat and generic type.

 \begin{corollary} \label{cor:hyp-examples} The classification of hyperbolic equations of the form $F(x,y,z,p,q,r,t)=0$, $G(x,y,z,p,q,r,s)=0$, and $rt = f(x,y,z,p,q,s)$ is given in Table~{\rm \ref{table:hyp-examples}} below.
 \end{corollary}

 \begin{proof} The hyperbolicity condition in each case is clear.  Def\/ine the function
 \begin{gather*}
 \Delta_{r,t}^F = F_r{}^2 F_{tt} - 2F_r F_t F_{rt} + F_t{}^2 F_{rr},
 \end{gather*}
 and similarly for $\Delta_{r,s}^G$.
 Without loss of generality $G_s$, $f_s \geq 0$.  The calculation of $\tilde{I}_1$, $\tilde{I}_2$ leads to
 \begin{alignat*}{4}
 &F(x,y,z,p,q,r,t)=0: \qquad && \tilde{I}_1 = \frac{-F_t{}^2\Delta_{r,t}^F}{2(-F_t F_r)^{3/2}},\qquad && \tilde{I}_2 = \frac{-F_r{}^2\Delta_{r,t}^F}{2(-F_t F_r)^{3/2}}, &\\
 &G(x,y,z,p,q,r,s)=0: \qquad && \tilde{I}_1 = 0,\qquad && \tilde{I}_2 = \frac{-\Delta_{r,s}^G}{G_s}, &\\
 &rt = f(x,y,z,p,q,s): \qquad && \tilde{I}_1 = \frac{(f_{ss} - 2) r^2}{\sqrt{f_s{}^2 - 4rt}},\qquad  && \tilde{I}_2 = \frac{(f_{ss} - 2) t^2}{\sqrt{f_s{}^2 - 4rt}}, &\\
 &rt = -f(x,y,z,p,q,s):\qquad  && \tilde{I}_1 = \frac{-(f_{ss} + 2) r^2}{\sqrt{f_s{}^2 - 4rt}},\qquad && \tilde{I}_2 = \frac{-(f_{ss} + 2) t^2}{\sqrt{f_s{}^2 - 4rt}}.&
 \end{alignat*}
 For $F(x,y,z,p,q,r,t)=0$: Since $i_F^*(F_t F_r) < 0$, then either $I_1$, $I_2$ both vanish or both do not vanish, i.e.\ either class 6-6 or class 7-7.  The vanishing of $I_1$, $I_2$ is completely characterized by the vanishing of $i_F^*(\Delta_{r,t}^F)$.  By Theorem \ref{thm:MA}, we know what all class 6-6 equations of the form $F(x,y,z,p,q,r,t)=0$ look like.  Hence,
 \begin{gather*}
 i_F^*(\Delta_{r,t}^F)=0 \qquad \mbox{if\/f its locus can be given by} \quad F(x,y,z,p,q,r,t) = a r + bt + c = 0,
 \end{gather*}
 where $a$, $b$, $c$ are functions of $x$, $y$, $z$, $p$, $q$ only.  The proof for $G(x,y,z,p,q,r,s)=0$ is similar and the result for the last equation is immediate.
 \end{proof}

 \begin{table}[h]\centering
 \caption{General examples of hyperbolic equations and their types.}
 \label{table:hyp-examples}
 \vspace{1mm}
 $\begin{array}{|c|c|c|c|c|} \hline
 \mbox{Equation} & \begin{tabular}{c} Hyperbolicity\\ condition \end{tabular} & \mbox{Type}\\ \hline \hline
 F(x,y,z,p,q,r,t)=0 & i_F^*(F_r F_t) < 0 &\begin{array}{c} \mbox{6-6 if\/f } F \mbox{ is an af\/f\/ine function of } r,t \mbox{ (*)} \\ \mbox{7-7 otherwise} \end{array} \\ \hline
 G(x,y,z,p,q,r,s)=0 & i_G^*(G_s) \neq 0 & \begin{array}{c} \mbox{6-6 if\/f } G \mbox{ is an af\/f\/ine function of } r,s \mbox{ (*)} \\ \mbox{6-7 otherwise } \end{array} \\ \hline
 rt = f(x,y,z,p,q,s) & 4f < f_s{}^2 & \begin{array}{c} \mbox{Assuming $rt \neq 0$:}\\ \mbox{6-6 if\/f } f_{ss}=2 \\ \mbox{7-7 if\/f } f_{ss} \neq 2 \end{array} \\ \hline
 \end{array}$

 \vspace{1mm}

 (*) More precisely, it is the zero-locus of such a function.
 \end{table}

 \begin{remark}
 Hyperbolic equations of Goursat and generic type are necessarily {\em non}-variational.  This is because a variational formulation for a second order PDE requires a f\/irst order Lagrangian (density) $L(x,y,z,p,q)$ and the corresponding Euler--Lagrange equation is
 \begin{gather*}
 \Parder{L}{z} - D_x\left( \Parder{L}{p} \right)- D_y\left( \Parder{L}{q} \right) = 0,
 \end{gather*}
 where $D_x$ and $D_y$ are total derivative operators
 \begin{gather*}
 D_x = \parder{x} + p \parder{z} + r\parder{p} + s\parder{q}, \qquad
 D_y = \parder{y} + q \parder{z} + s\parder{p} + t\parder{q}.
 \end{gather*}
 Thus, the Euler--Lagrange equation is quasi-linear and, if hyperbolic, is of \MA/ type.
 \end{remark}

 For the remainder of the paper we will deal exclusively with the generic case.  In this case~$U_1$,~$U_2$ in \eqref{U1U2-streq} are nonzero and can be further normalized through a coframe adaptation.  Before carrying out this normalization, we recall some more basic def\/initions.

 \begin{definition}
 Given a Pfaf\/f\/ian system $I$ on a manifold $\Sigma$, recall that the {\em first derived system} $I^{(1)} \subset I$ is the Pfaf\/f\/ian system def\/ined by the short exact sequence
 \begin{equation*}
 0 \lra I^{(1)} \hooklongrightarrow I \stackrel{\pi \circ d}{\lra} dI \mbox{ mod } I \lra 0,
 \end{equation*}
 where $\pi : \Omega^*(\Sigma) \ra \Omega^*(\Sigma) / I$ be the canonical surjection.  (Here we abuse notation and identify~$I$ with the {\em algebraic} ideal in $\Omega^*(\Sigma)$ that it generates.)  Iteratively, we def\/ine the {\em derived flag}
 $\cdots \subset I^{(k)} \subset \cdots \subset I^{(1)} \subset I$.
 \end{definition}

 \begin{remark} $I$ is completely integrable (in the Frobenius sense) if\/f $I^{(1)} = I$.
 \end{remark}

 Since $d$ commutes with pullbacks, each derived system $I^{(k)}$ is invariant under any automorphism of $I$, i.e.\ if $\phi \in {\rm Aut}(I)$, then $\phi^* I^{(k)} = I^{(k)}$.

\begin{definition}
 For hyperbolic equations, def\/ine
 \begin{gather*}
 {\rm Char}(I_F,dM_i)  = \{ X \in\mathfrak{X}(\Sigma_7) : X \in I_F^\perp,~ X \intprod dM_i \subset I_F \}, \\
 C(I_F,dM_i)  = {\rm Char}(I_F,dM_i)^\perp.
 \end{gather*}
 \end{definition}

 We now normalize the coef\/f\/icients $U_1$, $U_2$ in the generic case.  Moreover, explicit generators for the f\/irst few systems in the derived f\/lag of $C(I_F,dM_1)$ and $C(I_F,dM_2)$ are obtained.  The proofs of the following theorem and subsequent corollaries are provided in Appendix~\ref{app:gen-hyp}.

 \begin{theorem}[{\bf Generic hyperbolic structure equations}] \label{generic-hyp-str-eqns}
     Given any generic hyperbolic equation $F=0$ with nondegenerate parametrization $i_F : \Sigma_7 \ra J^2(\R^2,\R)$, there is an associated coframe $\bm\omega = \{ \omega^i \}_{i=1}^7$ on $\Sigma_7$ such that
    \begin{enumerate}\itemsep=0pt
    \item[1)] $I_F = \{ \omega^1, \omega^2, \omega^3 \}$, \ $I_F^{(1)} = \{ \omega^1 \}$, \ $
            M_1 = \{ \omega^1, \omega^2 \}$, \  $            M_2 = \{ \omega^1, \omega^3 \}$,
    \item[2)] we have the structure equations
     \begin{alignat}{3}
    & d\omega^1 \equiv \omega^3 \wedge \omega^6 + \omega^2 \wedge \omega^4 \quad && \mod I_F^{(1)}, & \nonumber \\
     &d\omega^2 \equiv \omega^4 \wedge \omega^5 + \omega^3 \wedge \omega^7 \quad && \mod M_1, &\label{StrEqns123}\\
    & d\omega^3 \equiv \omega^6 \wedge \omega^7 + \epsilon \omega^2 \wedge \omega^5 \quad&& \mod M_2, & \nonumber
        \end{alignat}
   where $\epsilon = \pm 1$,
   \item[3)] for {\em some} choice of coframe satisfying the above structure equations, we have
 \begin{alignat}{3}
  & C(I_F,dM_1) = \{ \omega^1, \omega^2, \omega^3, \omega^4, \omega^5 \}, \qquad & &
   C(I_F,dM_2) = \{ \omega^1, \omega^2, \omega^3, \omega^6, \omega^7 \}, &\nonumber \\
 &  C(I_F,dM_1)^{(1)} = \{ \omega^1, \omega^2, \omega^4, \omega^5 \},\qquad &&
   C(I_F,dM_2)^{(1)} = \{ \omega^1, \omega^3, \omega^6, \omega^7 \}, &  \label{charsys}\\
 &  C(I_F,dM_1)^{(2)} = \{ \omega^4, \omega^5 \},\qquad &&
   C(I_F,dM_2)^{(2)}  =  \{ \omega^6, \omega^7 \}. &\nonumber
 \end{alignat}
  \end{enumerate}
 \end{theorem}

 \begin{corollary} \label{gamma-cor}
 For the choice of coframe as in Theorem {\rm \ref{generic-hyp-str-eqns}}, we have the additional structure equations
 \begin{alignat}{3}
 &  d\omega^4 \equiv \epsilon \omega^5 \wedge \omega^6 \quad &&\mod \{ \omega^1, \omega^2, \omega^4 \}, & \nonumber\\
 &  d\omega^5 \equiv 0\quad  &&\mod \{\omega^1, \omega^2, \omega^4, \omega^5\}, & \nonumber\\
 &  d\omega^6 \equiv - \omega^4 \wedge \omega^7 \quad &&\mod \{\omega^1, \omega^3, \omega^6\}, & \label{StrEqns4567}\\
 &  d\omega^7 \equiv 0 \quad &&\mod \{\omega^1, \omega^3, \omega^6, \omega^7\}. &\nonumber
 \end{alignat}
 \end{corollary}

 We will refer to \eqref{StrEqns123} and \eqref{StrEqns4567} collectively as the generic hyperbolic structure equations.

 \begin{corollary} \label{epsilon-cor}
 $\epsilon$ is a contact invariant, and moreover $\epsilon = {\rm sgn}(I_1 I_2)$.
 \end{corollary}

 \begin{example} From Table~\ref{table:hyp-examples} and the proof of Corollary \ref{cor:hyp-examples}, we see $\epsilon=1$ for:{\samepage
 \begin{itemize}\itemsep=0pt
 \item $F(x,y,z,p,q,r,t)=0$ whenever $F$ is {\em not} an af\/f\/ine function of $r$, $t$.
 \item $rt = f(x,y,z,p,q,s)$ whenever $f_{ss} \neq 2$.
 \end{itemize}}
 \end{example}

 \begin{remark} We have the following dictionary of notations for the adapted coframe labelling:
$$
   \begin{array}{c|c|c|c}
   & \mbox{Gardner--Kamran \cite{GK1993}} & \mbox{Vranceanu \cite{Vranceanu1937}} & \mbox{The}\\ \hline
   I_F & \omega^1, \ \pi^2, \ \pi^3  & ds^1, \ ds^2, \ ds^3  & \omega^1, \ \omega^2, \ \omega^3 \\
   M_1 & \omega^1, \ \pi^2  & ds^1, \ ds^2  & \omega^1, \ \omega^2 \\
   M_2 & \omega^1, \ \pi^3  & ds^1, \ ds^3  & \omega^1, \ \omega^3 \\
   C(I_F,dM_1) & \omega^1, \ \pi^2, \ \pi^3, \ \omega^4, \ \omega^5  & ds^1, \ ds^2, \ ds^3, \ ds^5, \ ds^6  & \omega^1, \ \omega^2, \ \omega^3, \ \omega^4, \ \omega^5 \\
   C(I_F,dM_2) & \omega^1, \ \pi^2, \ \pi^3, \ \omega^6, \ \omega^7  & ds^1, \ ds^2, \ ds^3, \  ds^4, \ ds^7  & \omega^1, \ \omega^2, \ \omega^3, \ \omega^6, \ \omega^7
 \end{array}
$$
 \end{remark}

 \section{The structure group and the Cartan equivalence problem}
 \label{str-grp-eq-prb}

 In this section, we reformulate the problem of contact-equivalence of PDE as a Cartan equivalence problem.  The reader will notice the similarities in the calculation of the structure group in this section and in the calculations in the proof of Corollary~\ref{epsilon-cor} provided in Appendix~\ref{app:gen-hyp}.

 For any $\phi \in \contactmaps{+}$,
 \begin{gather*}
 \phi^* I_{\bar{F}}^{(1)} = I_F^{(1)}, \qquad \phi^*(C(I_{\bar{F}},d\bar{M}_i)^{(k)}) = C(I_F,dM_i)^{(k)}, \quad i=1,2, \quad \forall \; k \geq 0.
 \end{gather*}
 Consequently, with respect to the adapted coframe $\bm\omega$ on $\Sigma_7$ (as specif\/ied in Theorem \ref{generic-hyp-str-eqns}) and corresponding coframe $\bm{\bar\omega}$ on $\bar\Sigma_7$, we have
 \begin{gather*}
   \phi^*\mat{c}{\bar\omega^1 \\ \bar\omega^2\\ \bar\omega^3\\ \bar\omega^4\\ \bar\omega^5\\ \bar\omega^6\\ \bar\omega^7}
   = \mat{ccccccc}{
    \lambda_1 & 0 & 0 & 0 & 0 & 0 & 0\\
   \mu_1 & \lambda_2 & 0 & 0 & 0 & 0 & 0\\
   \mu_2 & 0 & \lambda_3 & 0 & 0 & 0 & 0\\
   0 & 0 & 0 & \lambda_4 & \nu_1 & 0 & 0\\
   0 & 0 & 0 & \mu_3 & \lambda_5 & 0 & 0\\
   0 & 0 & 0 & 0 & 0 & \lambda_6 & \nu_2\\
   0 & 0 & 0 & 0 & 0 & \mu_4 & \lambda_7\\
   }
   \mat{c}{\omega^1 \\ \omega^2\\ \omega^3\\ \omega^4\\ \omega^5\\ \omega^6\\ \omega^7}.
 \end{gather*}
 Applying $\phi^*$ to the $d\bar\omega^1$ structure equation in \eqref{StrEqns123} yields
 \begin{gather*}
 \phi^* d\bar\omega^1 = d\phi^*\bar\omega^1 = d(\lambda_1 \omega^1) \equiv
 \lambda_1 (\omega^3 \wedge \omega^6 + \omega^2 \wedge \omega^4) \quad \mod I_F^{(1)},
 \end{gather*}
 and also
 \begin{gather*}
 \phi^* d\bar\omega^1  \equiv \lambda_3 \omega^3 \wedge (\lambda_6 \omega^6 + \nu_2 \omega^7) + \lambda_2 \omega^2 \wedge (\lambda_4 \omega^4 + \nu_1 \omega^5)\quad \mod I_F^{(1)},
 \end{gather*}
 which implies $\nu_1=\nu_2=0$, $\lambda_1 = \lambda_3 \lambda_6 = \lambda_2 \lambda_4$.  Similarly, using the $d\bar\omega^2$, $d\bar\omega^3$ equations yields
 \begin{alignat*}{3}
   &\lambda_1 = \lambda_3 \lambda_6 = \lambda_2 \lambda_4, \qquad && \nu_1=\nu_2=0, &\\
   &\lambda_2 = \lambda_4\lambda_5 = \lambda_3 \lambda_7, \qquad && \mu_1 = \lambda_3 \mu_4, &\\
   &\lambda_3 = \lambda_6\lambda_7 = \lambda_2 \lambda_5, \qquad && \mu_2 = \epsilon \lambda_2 \mu_3.&
 \end{alignat*}
 Then
  \begin{gather*}
   \beta :=\frac{\lambda_6}{\lambda_4} = \frac{\lambda_2}{\lambda_3} = \frac{1}{\lambda_5}  = \frac{\lambda_4}{\lambda_2}, \qquad  \beta= \frac{\lambda_2}{\lambda_3} = \lambda_7 = \frac{\lambda_3}{\lambda_6} \quad\Rightarrow\quad \beta^4 = 1 \quad\Rightarrow\quad \beta = \pm 1,
  \end{gather*}
  and so
  \begin{gather}
  (\lambda_1, \lambda_2, \lambda_3, \lambda_4, \lambda_5, \lambda_6,\lambda_7) = (\beta a_1{}^2, a_1,\beta a_1,\beta a_1,\beta,a_1,\beta), \nonumber\\
  (\mu_1,\mu_2,\mu_3,\mu_4) = (\beta a_1a_2,\epsilon a_1a_3,a_3,a_2),
\qquad \forall \; a_1 \neq 0, a_2,a_3 \in \R.
  \label{str-grp-calc}
  \end{gather}
 This leads us to def\/ine
 \begin{gather*}
 S = {\rm diag}(-1,1,-1,-1,-1,1,-1),
 \end{gather*}
 and the {\em connected} matrix Lie group
 \begin{gather}
  G^0 = \left\{ M({\bf a}) : {\bf a} \in \R^+ \times \R^2 \right\}, \qquad
 M({\bf a}) = \mat{ccccccc}{
   a_1{}^2 & 0 & 0 & 0 & 0 & 0 & 0\\
   a_1 a_2 & a_1 & 0 & 0 & 0 & 0 & 0\\
   \epsilon a_1 a_3 & 0 & a_1 & 0 & 0 & 0 & 0\\
   0 & 0 & 0 & a_1 & 0 & 0 & 0\\
   0 & 0 & 0 & a_3 & 1 & 0 & 0\\
   0 & 0 & 0 & 0 & 0 & a_1 & 0\\
   0 & 0 & 0 & 0 & 0 & a_2 & 1
   }.
 \label{G0-group}
 \end{gather}
 Let us also def\/ine
 \begin{gather*}
 R = \mat{ccccccc}{
   -\epsilon & 0 & 0 & 0 & 0 & 0 & 0\\
   0 & 0 & \epsilon & 0 & 0 & 0 & 0\\
   0 & -\epsilon & 0 & 0 & 0 & 0 & 0\\
   0 & 0 & 0 & 0 & 0 & -1 & 0\\
   0 & 0 & 0 & 0 & 0 & 0 & -\epsilon\\
   0 & 0 & 0 & 1 & 0 & 0 & 0\\
   0 & 0 & 0 & 0 & -\epsilon & 0 & 0\\
   } \quad\Rightarrow\quad R^2 = {\rm diag}(1,-1,-1,-1,1,-1,1).
 \end{gather*}
 We note that $R^4=S^2=1$ and we have the relation $RS = SR^{-1}$, and consequently $R$, $S$ generate the dihedral group of order 8
 \begin{gather*}
 D_8 = \langle R,S : R^4=S^2=SRSR=1 \rangle.
 \end{gather*}

 The results \eqref{str-grp-calc} of the previous calculations establish that
 \begin{gather*}
 \phi \in \contactmaps{+} \qquad \mbox{if\/f} \quad \phi^* \bm{\bar\omega} = g \bm\omega, \qquad \mbox{where} \quad g : \Sigma_7 \ra G^+,
 \end{gather*}
 where $G^+$ is the group generated by $G^0$, $S$, $R^2$, which we can realize as the semi-direct product
 \begin{gather*}
 G^+ = G^0 \rtimes \langle S, R^2 \rangle
 \end{gather*}
 induced by the adjoint action.
 If $\phi_0^* \bm{\bar\omega} = R \bm\omega$, then $\phi_0 \in \contactmaps{-}$.  Conversely, given any $\phi \in \contactmaps{-}$, we have $\phi = \phi_0 \circ \tilde\phi$, where $\tilde\phi \in {\rm Aut}^+(I_F)$.  Thus,
 \begin{gather*}
 \phi^* \bm{\bar\omega} = \tilde\phi^* \phi_0^* \bm{\bar\omega} = \tilde\phi^* R\bm\omega = Rg\bm\omega = {\rm Ad}_R(g) R\bm\omega, \qquad \forall \;g \in G^+.
 \end{gather*}
 Since $G^0$ is ${\rm Ad}_R$-invariant, and ${\rm Ad}_R(S) = SR^2$, then $G^+$ is ${\rm Ad}_R$-invariant and so
  \begin{gather*}
 \phi \in \contactmaps{-} \qquad \mbox{if\/f} \quad \phi^* \bm{\bar\omega} = g \bm\omega, \qquad \mbox{where} \quad g : \Sigma_7 \ra G^-,
 \end{gather*}
 where $G^- = G^+ \cdot R$.  (Note that $G^-$ is {\em not} a group.)  Consequently, we def\/ine
 \begin{gather*}
 G = G^0 \rtimes D_8,
 \end{gather*}
 and we have established:

 \begin{lemma}
 $\phi \in \contactmaps{}$ if and only if
 \begin{gather}
 \phi^* \bm{\bar\omega} = g \bm\omega,  \qquad \mbox{for some} \quad g : \Sigma_7 \ra G.
 \label{base-equivalence}
 \end{gather}
 \end{lemma}

  The group $G$ will play the role of our initial structure group in the application of the Cartan equivalence method \cite{Cartan1953, Gardner1989, Olver1995}.  Specif\/ically, the Cartan equivalence problem for generic hyperbolic equations can be stated as follows: Given the coframes $\bm\omega$, $\bm{\bar\omega}$ on $\Sigma_7$, $\bar\Sigma_7$ respectively, f\/ind necessary and suf\/f\/icient conditions for the existence of a local dif\/feomorphism $\phi : \Sigma_7 \ra \bar\Sigma_7$ satisfying \eqref{base-equivalence}.  This is also known as the isomorphism problem for $G$-structures $(\bm{\omega},G)$ and $(\bm{\bar\omega},G)$.

 \begin{remark}
 Vranceanu (c.f.\ page~366 in~\cite{Vranceanu1937}) considered the equivalence problem with respect to a smaller group which, in our notation, is $G^0 \rtimes \langle R^2 \rangle$.  This has index~4 in~$G$.
 \end{remark}

 The solution of the general Cartan equivalence problem leads to either the structure equations of an $\{e\}$-structure or of an inf\/inite Lie pseudo-group.  However, for the equivalence problem for generic hyperbolic equations only the former case occurs.  In particular, we will show in the next section that we are led to $\{e\}$-structures on $\Sigma_7 \times G_\Gamma$, where $G_\Gamma \subset G$ is a subgroup of dimension at most {\em two}.  (Dif\/ferent $\{e\}$-structures arise due to normalizations of {\em nonconstant} type, and will depend on choices of normal forms $\Gamma$ in dif\/ferent orbits.)  For the moment, let us recall the general solution to the coframe ($\{e\}$-structure) equivalence problem.  Our description below is abbreviated from the presentation given in \cite{Olver1995}.

 Let $\bm\Theta$, $\bm{\bar\Theta}$ be local coframes on manifolds $M$, $\bar{M}$ respectively of dimension $m$, and let $\Phi$ satisfy $\Phi^* \bm{\bar\Theta} = \bm\Theta$.  If the structure equations for the $\{e\}$-structures are correspondingly
 \begin{gather*}
 d\Theta^a = \frac{1}{2} T^a{}_{bc} \Theta^b \wedge \Theta^c, \qquad
 d\bar\Theta^a = \frac{1}{2} \bar{T}^a{}_{bc} \bar\Theta^b \wedge \bar\Theta^c, \qquad 1 \leq a,b,c \leq m,
 \end{gather*}
 then by commutativity of $\Phi^*$ and $d$, the structure functions $T^a{}_{bc}$ are invariants, i.e.
 \begin{gather*}
 \bar{T}^a{}_{bc} \circ \Phi = T^a{}_{bc}.
 \end{gather*}
 For any local function $f$ on $M$, def\/ine the {\em coframe derivatives} $\Parder{f}{\Theta^a}$ by
 \begin{gather*}
 df = \Parder{f}{\Theta^k} \Theta^k.
 \end{gather*}
 Let us write the iterated derivatives of the structure functions as
 \begin{gather*}
 T_\sigma = \frac{\partial^s T^a{}_{bc}}{\partial \Theta^{k_s} \cdots \partial \Theta^{k_1}}, \qquad \mbox{where} \quad \sigma = (a,b,c,k_1,\dots,k_s)\qquad \mbox{and} \quad s = {\rm order}(\sigma)
 \end{gather*}
 and $1 \leq a,b,c,k_1,\dots, k_s \leq m$.  We repeat this construction for the barred variables.
 Necessarily, again as a consequence of commutativity of~$\Phi^*$ and~$d$, the derived structure functions~$T_\sigma$ and~$\bar{T}_\sigma$ satisfy the {\em invariance equations}
 \begin{gather}
 \bar{T}_\sigma(\bar{x}) = T_\sigma(x), \qquad \mbox{when} \quad \bar{x} = \Phi(x), \qquad \forall \; {\rm order}(\sigma) \geq 0. \label{invar}
 \end{gather}
 Note that these equations are not independent: there are generalized Jacobi identities (which we will not describe explicitly here) which allow the permutation of the coframe derivatives, so in general only {\em nondecreasing} coframe derivative indices are needed.

 \begin{definition} Let $\bm\Theta$ be a coframe with def\/ined on an open set $U \subset M$.
 \begin{enumerate}\itemsep=0pt
 \item Let $\mathbb{K}^{(s)}$ be the Euclidean space of dimension equal to the number of multi-indices
 \begin{gather*}
 \sigma = (a,b,c,k_1,\dots,k_r), \qquad b<c, \qquad k_1 \leq \cdots \leq k_r, \qquad 0 \leq r \leq s.
 \end{gather*}
 \item The $s^{\rm th}$ order {\em structure map} associated to $\bm\Theta$ is
 \begin{gather*}
 {\bf T}^{(s)} : U \ra \mathbb{K}^{(s)}, \qquad z_\sigma = T_\sigma(x), \qquad {\rm order}(\sigma) \leq s.
 \end{gather*}
 \item The coframe $\bm\Theta$ is {\em fully regular} if ${\bf T}^{(s)}$ is regular for all $s \geq 0$.  In this case,
 let $\rho_s = {\rm rank}({\bf T}^{(s)})$, and def\/ine the {\em rank} of $\bm\Theta$ as the minimal $s$ such that $\rho_s = \rho_{s+1}$.
 \item The $s^{\rm th}$ order {\em classifying set} is
 ${\cal C}^{(s)}(\bm\Theta,U) = \{ {\bf T}^{(s)}(x) : x \in U \} \subset \mathbb{K}^{(s)}.$
 \end{enumerate}
 \end{definition}

 As a consequence of the invariance equations \eqref{invar}, if $\bm\Theta$ and $\bm{\bar\Theta}$ are equivalent coframes via $\Phi : U \ra \bar{U}$, then
 \begin{gather*}
 {\cal C}^{(s)}(\bm\Theta,U) = {\cal C}^{(s)}(\bm{\bar\Theta},\Phi(U)), \qquad \forall \; s \geq 0.
 \end{gather*}
 This is suf\/f\/icient in the fully regular case.  We refer the reader to \cite{Olver1995} for a proof of the following theorem.

 \begin{theorem} \label{equiv-soln} Suppose $\bm\Theta$, $\bm{\bar\Theta}$ are fully regular coframes on $U$, $\bar{U}$ respectively.
  There exists $\Phi$ satisfying $\Phi^*\bm{\bar\Theta} = \bm\Theta$ if and only if for each $s \geq 0$,  ${\cal C}^{(s)}(\bm\Theta,U) \cap {\cal C}^{(s)}(\bm{\bar\Theta},\bar{U})$ is nonempty.  The set of self-equivalences $\Phi$ (i.e.\ satisfying $\Phi^* \bm\Theta = \bm\Theta$) defines a $p$-dimensional local Lie group of transformations, where $p = m - {\rm rank}(\bm\Theta) \geq 0$.
 \end{theorem}

 \section{Nine-dimensional maximal symmetry}
 \label{9d-sym}

 The solution to the Cartan equivalence problem \eqref{base-equivalence} begins by lifting the problem to the {\em left} principal bundles $\Sigma_7 \times G \stackrel{\pi}{\ra} \Sigma_7$ and $\bar\Sigma_7 \times G \stackrel{\bar\pi}{\ra} \bar\Sigma_7$ by def\/ining
  \begin{gather*}
 \bm{\hat{\omega}}|_{(u,g)} = g \pi^* \bm{\omega}|_{u}, \qquad
 \bm{\hat{\bar\omega}}|_{(\bar{u},g)} = g \bar\pi^* \bm{\bar\omega}|_{\bar{u}}, \qquad
  \mbox{where} \quad u \in \Sigma_7, \quad \bar{u} \in \bar\Sigma_7, \quad g \in G,
 \end{gather*}
 and noting the following key lemma \cite{Gardner1989}.

 \begin{lemma} \label{Gardner-lemma}
 There exists an equivalence $\phi$ as in \eqref{base-equivalence} if and only if there exists a local diffeomorphism $\Phi : \Sigma_7 \times G \ra \bar\Sigma_7 \times G$ satisfying $\Phi^* \bm{\hat{\bar\omega}} = \bm{\hat\omega}$.
 \end{lemma}

 Identifying the coframe $\bm\omega$ on $\Sigma_7$ with its pullback by the canonical projection $\Sigma_7 \times G \ra \Sigma_7$, we can write
 \begin{gather*}
   \hat\omega^i = g^i{}_j \omega^j, \qquad g \in G.
 \end{gather*}
 Using \eqref{gamma-defn}, the structure equations for these lifted forms are then
  \begin{gather*}
   d\hat\omega^i  = (dg \cdot g^{-1})^i{}_j \wedge \hat\omega^j + \frac{1}{2} \hat\gamma^j{}_{k\ell} \hat\omega^k \wedge \hat\omega^\ell,
 \end{gather*}
 where the coef\/f\/icients $\gam{i}{jk}$ transform tensorially under the $G$-action
 \begin{gather}
   \hat\gamma^i{}_{jk} := g^i{}_\ell \gamma^\ell{}_{mn} (g^{-1})^m{}_j (g^{-1})^n{}_k,
 \label{gamma-transform}
 \end{gather}
 and $dg \cdot g^{-1}$ refers to the {\em right}-invariant Maurer--Cartan form on $G$.  Since $D_8$ is discrete, then if $(g,k) \in G^0 \times D_8$,
 \begin{gather*}
 d(gk) \cdot (gk)^{-1} = dg \cdot k \cdot k^{-1} g^{-1} = dg \cdot g^{-1},
 \end{gather*}
 and so we can identify the Maurer--Cartan form on $G$ with that on $G^0$.  For $g = M(a_1,a_2,a_3) \in G^0$ as in \eqref{G0-group}, we have $g^{-1} = M\left(\frac{1}{a_1},-\frac{a_2}{a_1},-\frac{a_3}{a_1} \right)$ and
  \begin{gather*}
   dg \cdot g^{-1} = \mat{ccccccc}{
   2\alpha^1 & 0 & 0 & 0 & 0 & 0 & 0\\
   \alpha^2 & \alpha^1 & 0 & 0 & 0 & 0 & 0\\
   \epsilon \alpha^3 & 0 & \alpha^1 & 0 & 0 & 0 & 0\\
   0 & 0 & 0 & \alpha^1 & 0 & 0 & 0\\
   0 & 0 & 0 & \alpha^3 & 0 & 0 & 0\\
   0 & 0 & 0 & 0 & 0 & \alpha^1 & 0\\
   0 & 0 & 0 & 0 & 0 & \alpha^2 & 0
   },
 \end{gather*}
 where
 \begin{gather*}
 \alpha^1 = \frac{da_1}{a_1}, \qquad \alpha^2 = \frac{da_2}{a_1}, \qquad \alpha^3 = \frac{da_3}{a_1}
 \end{gather*}
 are a basis for the right-invariant 1-forms on $G^0$ and hence $G$.  Identifying $\alpha^i$ on $G$ with their pullback by the canonical projection $\Sigma_7 \times G \ra G$, we have the structure equations for the lifted coframe:
 \begin{gather}
   d\hat\omega^1 = 2\alpha^1 \wedge \hat\omega^1 + \hat\omega^3 \wedge \hat\omega^6 + \hat\omega^2 \wedge \hat\omega^4 + \eta_1 \wedge \hat\omega^1,\nonumber\\
   d\hat\omega^2 = \alpha^2 \wedge \hat\omega^1 + \alpha^1 \wedge \hat\omega^2 + \hat\omega^4 \wedge \hat\omega^5 + \hat\omega^3 \wedge \hat\omega^7 + \eta_2 \wedge \hat\omega^1 + \eta_3 \wedge \hat\omega^2,\nonumber\\
   d\hat\omega^3 = \epsilon \alpha^3 \wedge \hat\omega^1 + \alpha^1 \wedge \hat\omega^3 + \hat\omega^6 \wedge \hat\omega^7 + \epsilon\hat\omega^2 \wedge \hat\omega^5 + \eta_4 \wedge \hat\omega^1 + \eta_5 \wedge \hat\omega^3,\nonumber\\
   d\hat\omega^4 = \alpha^1 \wedge \hat\omega^4 + \eta_6 \wedge \hat\omega^1 + \eta_7 \wedge \hat\omega^2 + \eta_8 \wedge \hat\omega^4 + \eta_9 \wedge \hat\omega^5,\nonumber\\
   d\hat\omega^5 = \alpha^3 \wedge \hat\omega^4 + \eta_{10} \wedge \hat\omega^1 + \eta_{11} \wedge \hat\omega^2 + \eta_{12} \wedge \hat\omega^4 + \eta_{13} \wedge \hat\omega^5 ,\label{G-lifted-coframe}\\
   d\hat\omega^6 = \alpha^1 \wedge \hat\omega^6 + \eta_{14} \wedge \hat\omega^1 + \eta_{15} \wedge \hat\omega^3 + \eta_{16} \wedge \hat\omega^6 + \eta_{17} \wedge \hat\omega^7,\nonumber\\
   d\hat\omega^7 = \alpha^2 \wedge \hat\omega^6 + \eta_{18} \wedge \hat\omega^1 + \eta_{19} \wedge \hat\omega^3 + \eta_{20} \wedge \hat\omega^6 + \eta_{21} \wedge \hat\omega^7,\nonumber\\
   d\alpha^1 = 0,\nonumber\\
   d\alpha^2 = -\alpha^1 \wedge \alpha^2,\nonumber\\
   d\alpha^3 = -\alpha^1 \wedge \alpha^3,\nonumber
   \end{gather}
 where $\eta_i$ are semi-basic 1-forms with respect to the projection $\Sigma_7 \times G \ra \Sigma_7$.
 The structure equations for the lifted forms $\hat\omega^i$ can be written
 \begin{gather}
 d\hat\omega^i = a^i{}_{\rho j} \alpha^\rho \wedge \hat\omega^j + \frac{1}{2} \hgam{i}{jk} \hat\omega^j \wedge \hat\omega^k,
 \label{gen-lifted-streqns}
 \end{gather}
 where $a^i{}_{\rho j}$ are constants (c.f.\ Maurer--Cartan form) and $\hgam{i}{jk}$ is def\/ined as in~\eqref{gamma-transform}.

 \begin{definition}
 The degree of indeterminacy $r^{(1)}$ of a lifted coframe is the number of free variables in the set of  transformations $\alpha^\rho \mapsto \alpha^\rho + \lambda^\rho{}_i \hat\omega^i$ which preserve the structure equations for $d\hat\omega^i$.
 \end{definition}

 For later use, we note the following:

 \begin{lemma} For our lifted coframe $\bm\Theta = \{ \bm{\hat\omega}, \bm\alpha \}$ on $\Sigma_7 \times G$ satisfying \eqref{G-lifted-coframe}, we have $r^{(1)} = 0$.
 \label{lem:indeterminacy}
 \end{lemma}

 \begin{proof} From the $d\hat\omega^1$, $d\hat\omega^2$, $d\hat\omega^3$ equation in \eqref{G-lifted-coframe}, we must have
 \begin{gather*}
   \alpha^1 \mapsto \alpha^1 + \lambda \hat\omega^1, \qquad
   \alpha^2 \mapsto \alpha^2 + \lambda \hat\omega^2, \qquad
   \alpha^3 \mapsto \alpha^3 + \epsilon \lambda\hat\omega^3.
 \end{gather*}
 However, to preserve the form of $d\hat\omega^i$, $i=4,5,6,7$, we must have $\lambda=0$.  Since there are no free variables, then $r^{(1)}=0$.
 \end{proof}

 The goal in Cartan's solution algorithm is to reduce to an $\{e\}$-structure so that Theorem~\ref{equiv-soln} can be invoked.  This amounts to essentially adapting the coframes on the base, i.e.\ f\/ixing a~map $g : \Sigma_7 \ra G$.  Using Lemma~\ref{Gardner-lemma}, coef\/f\/icients in the structure equations are candidates for normalization, from which the structure group $G$ can be subsequently reduced.  However, we only use those coef\/f\/icients which are not af\/fected by the choice of any map $g : \Sigma_7 \ra G$.  Note that pulling the Maurer--Cartan forms back to the base by such a map will express each $\alpha^\rho$ in terms of the new coframe $\bm{\hat\omega}$ (pulled back to the base).  This motivates the following def\/inition.

 \begin{definition} Given a lifted coframe, {\em Lie algebra valued compatible absorption} refers to redef\/ining the right-invariant 1-forms $\alpha^\rho$ by
$\hat\alpha^\rho = \alpha^\rho + \lambda^\rho{}_i \hat\omega^i$, where $\lambda^\rho{}_i$ are functions on the bundle.
 The terms involving the coef\/f\/icients $\hgam{i}{jk}$ which cannot be eliminated by means of Lie algebra valued compatible absorption are called {\em torsion terms} and the corresponding coef\/f\/icients are referred to as {\em torsion coefficients}.
 \end{definition}

 From \eqref{G-lifted-coframe}, the $d\hat\omega^5$ and $d\hat\omega^7$ structure equations indicate that $\hat\gamma^5{}_{56}$, $\hat\gamma^5{}_{57}$, $\hat\gamma^7{}_{47}$, $\hat\gamma^7{}_{57}$ are torsion coef\/f\/icients.  Using \eqref{StrEqns123}, \eqref{StrEqns4567}, and the tensor transformation law \eqref{gamma-transform} for the $\gamma$'s, we see that there is a well-def\/ined $G$-action on $\R^4$ (i.e.\ the range of $(\gam{5}{56},\gam{5}{57},\gam{7}{47}, \gam{7}{57})$) given by the formulas
 \begin{equation}
 \begin{array}{|c|c|c|c|}\hline
 & G^0\mbox{-action by } g=M(a_1,a_2,a_3) & R\mbox{-action} & S\mbox{-action}\\ \hline
 \hat\gamma^5{}_{56} & \frac{1}{a_1} (\gamma^5{}_{56} - \gamma^5{}_{57} a_2 + a_3 \epsilon) & -\gamma^7{}_{47} & \gamma^5{}_{56}\\
 \hat\gamma^7{}_{47} & \frac{1}{a_1} (\gamma^7{}_{47} - a_2 - \gamma^7{}_{57} a_3) & \gamma^5{}_{56} & -\gamma^7{}_{47}\\
 \hat\gamma^5{}_{57} & \gamma^5{}_{57} & \epsilon\gamma^7{}_{57} & -\gamma^5{}_{57}\\
 \hat\gamma^7{}_{57} & \gamma^7{}_{57} & \epsilon\gamma^5{}_{57} & -\gamma^7{}_{57}\\ \hline
 \end{array}
 \label{gam-main-transform}
 \end{equation}

 We can always normalize $\hat\gamma^5{}_{56}$ to zero by using the $G^0$-action and setting
 \begin{gather}
    a_3 = \epsilon(-\gam{5}{56} + \gam{5}{57} a_2).
    \label{a3-normalize}
 \end{gather}
 The matrix factorization
 \begin{gather*}
 M(a_1,a_2,\epsilon(-\gam{5}{56} + \gam{5}{57} a_2)) = M(a_1,a_2,\epsilon \gam{5}{57}  a_2)M(1,0,-\epsilon\gam{5}{56})
 \end{gather*}
 indicates that we can normalize $\gamma^5{}_{56}$ to 0 for the base coframe via
 \begin{gather*}
   \bar\omega^3 = -\gam{5}{56} \omega^1 + \omega^3, \qquad
   \bar\omega^5 = -\epsilon\gam{5}{56} \omega^4 + \omega^5.
 \end{gather*}
 This change of coframe is {\em admissible} in the sense that it preserves the form of the structure equations in \eqref{StrEqns123} and \eqref{StrEqns4567}.  (We henceforth drop the bars.)  Thus, we have the normal form $\Gamma = (\gam{5}{56}=0,\gam{5}{57},\gam{7}{47}, \gam{7}{57})$.  In general, however, this is a normalization of {\em nonconstant type} since $\Gamma$ still may depend on $x \in \Sigma_7$.  Pointwise, we def\/ine the reduced structure group $G_\Gamma$ as the stabilizer of $\Gamma$, i.e.\ it is the subgroup of $G$ preserving the structure equations together with the normalization given by $\Gamma$.  Clearly, the 1-parameter subgroup generated $M(1,0,a_3)$ ($a_3 \in \R$) yields a 1-dimensional orbit through $\Gamma$ and so $\dim(G_\Gamma) \leq 2$ since $\dim(G)=3$.

 The algorithm continues by means of further normalizations and reductions of the structure group until one of two possibilities occurs:
 \begin{enumerate}\itemsep=0pt
 \item[1)] the structure group has been reduced to the identity, i.e.\ get an $\{e\}$-structure on $\Sigma_7$, or
 \item[2)] the structure group has {\em not} been reduced to the identity but the structure group acts trivially on the torsion coef\/f\/icients.
 \end{enumerate}
 By Theorem \ref{equiv-soln}, the former possibility yields a symmetry group of dimension at most seven.  In the latter case, the next step in the algorithm is to prolong the problem to the space $\Sigma_7 \times G_\Gamma$.  Here, we have abused notation and written $G_\Gamma$ also for the structure group in the latter possibility above.  Since, by Lemma \ref{lem:indeterminacy}, $r^{(1)}=0$ with respect to the lifted coframe on $\Sigma \times G$, it is clear that we must have $r^{(1)}=0$ for the lifted coframe on $\Sigma_7 \times G_\Gamma$.  Finally, we invoke the following standard theorem (Proposition~12.1 in \cite{Olver1995}) written here in our notation:

 \begin{proposition} \label{prop:prolong}
 Let $\bm{\hat\omega}$, $\bm{\hat{\bar\omega}}$ be lifts of coframes $\bm\omega$, $\bm{\bar\omega}$ having the same structure group $G_\Gamma$, no group dependent torsion coefficients, and $r^{(1)}=0$.  Let $\bm{\hat\alpha}$, $\bm{\hat{\bar\alpha}}$ be modified Maurer--Cartan forms obtained by performing a full Lie algebra-valued compatible absorption.  Denote $\bm\Theta = \{ \bm{\hat\omega}, \bm{\hat\alpha} \}$, $\bm{\bar\Theta} = \{ \bm{\hat{\bar\omega}}, \bm{\hat{\bar\alpha}} \}$.  Then there exists $\phi : \Sigma_7 \ra \bar\Sigma_7$ satisfying $\phi^*\bm{\bar\omega} = g\bm\omega$ for some $g : \Sigma_7 \ra G_\Gamma$ if and only if there exists $\Phi : \Sigma_7 \times G_\Gamma \ra \bar\Sigma_7 \times G_\Gamma$ satisfying $\Phi^*\bm{\bar\Theta} = \bm\Theta$.
 \end{proposition}

 In other words, we have prolonged to an $\{e\}$-structure on $\Sigma_7 \times G_\Gamma$.  Since $\dim(G_\Gamma) \leq 2$ for any choice of $\Gamma$, then the symmetry group of the coframe is at most nine-dimensional.  Thus, we have proven:

 \begin{theorem} The (contact) symmetry group of any generic hyperbolic equation is finite dimensional and has maximal dimension~$9$.
 \end{theorem}

 In fact, this upper bound is sharp.  We will give explicit normal forms for all contact-equivalence classes of generic hyperbolic equations with 9-dimensional symmetry along with their corresponding symmetry generators and corresponding structure equations.

 Def\/ine
 \begin{gather*}
   m := \gamma^5{}_{57} \in C^\infty(\Sigma_7), \qquad
   n := \gamma^7{}_{57} \in C^\infty(\Sigma_7),
 \end{gather*}
 and note that although $m$ and $n$ are $G^0$-invariant, they are {\em not} $G$-invariant.  However, along each $G$-orbit
 the product $mn$ is invariant.

 We def\/ine two functions which will play an important role in the classif\/ications to follow.  Def\/ine
 \begin{gather*}
  \Delta_1 = mn + \epsilon, \qquad \Delta_2 = m^2 - \epsilon n^2
 \end{gather*}
 Note that $\Delta_1$ is a contact invariant, and $\Delta_2$ is a relative contact invariant: it is $G^+$-invariant, but under the $R$-action, $\hat\Delta_2 = -\epsilon \Delta_2$.

 \begin{corollary}\label{mn-epsilon} If a generic hyperbolic equation has 9-dimensional symmetry group, then \mbox{$\Delta_1\!=\!0$}.
  \end{corollary}

 \begin{proof} Under the assumption of maximal symmetry, all torsion coef\/f\/icients must be constant.  Thus, $\hat{m}$, $\hat{n}$ and consequently $m$, $n$ must be constant.  If $\Delta_1 \neq 0$, then there is a {\em unique} solution to the linear system
 \begin{gather}
   \mat{cc}{m & -\epsilon\\  1 & n}\mat{c}{a_2\\a_3} = \mat{c}{\gam{5}{56} \\ \gam{7}{47}},
 \label{mn-linear-sys}
 \end{gather}
 which yields the normalizations $\hgam{5}{56} = \hgam{7}{47} = 0$ and a {\em two} dimensional reduction of the initial structure group $G$.  Consequently, the stabilizer $G_\Gamma$ would be at most 1-dimensional and the symmetry group would be at most 8-dimensional.  Thus, we must have $\Delta_1 = 0$.
 \end{proof}

 \section{Complete structure equations}
 \label{complete-str-eqs}

  In Appendix \ref{Vranceanu-reduction}, we provide details of Vranceanu's reduction of the generic hyperbolic structure equations which allowed him to isolate the maximally symmetric and two sets of submaximally symmetric structures.

 \begin{theorem} \label{mn-str-eqns}
 Let $K^0 = \{ {\rm diag}(a_1^2,a_1,a_1,a_1,1,a_1,1) : a_1 > 0 \} \subset G$.  Consider a coframe $\{ \omega^i \}_{i=1}^7$ on $\Sigma_7$ satisfying the generic hyperbolic structure equations \eqref{StrEqns123} and \eqref{StrEqns4567}, and the corresponding lifted coframe on $\Sigma_7 \times K^0 \ra \Sigma_7$.  If:
 \begin{enumerate}\itemsep=0pt
 \item[1)] all torsion coefficients on which $K^0$ acts nontrivially are constants, and
 \item[2)] $K^0$ cannot be reduced to the identity,
 \end{enumerate}
 then the structure equations can be put in the form
  \begin{gather}
   d\omega^1 = \omega^3 \wedge \omega^6 + \omega^2 \wedge \omega^4,\nonumber \\
   d\omega^2 = \omega^4 \wedge \omega^5 + \omega^3 \wedge \omega^7 + \omega^2 \wedge \left( -\frac{3n}{2} \omega^5 + \frac{m}{2} \omega^7 \right),\nonumber\\
   d\omega^3 = \omega^6 \wedge \omega^7 + \epsilon \omega^2 \wedge \omega^5 + \omega^3 \wedge \left( -\frac{n}{2} \omega^5 + \frac{3m}{2} \omega^7 \right),\nonumber\\
   d\omega^4 = \epsilon \omega^5 \wedge \omega^6 + \omega^2 \wedge \left(B \omega^5 + \gam{4}{27} \omega^7\right) + \omega^4 \wedge \left(\frac{3n}{2} \omega^5 - \frac{m}{2} \omega^7\right), \label{Vranceanu-red-coframe}\\
   d\omega^5 = m \omega^5 \wedge \omega^7,\nonumber  \\
   d\omega^6 = - \omega^4 \wedge \omega^7 + \omega^3 \wedge \left(\gam{6}{35} \omega^5 + \epsilon B \omega^7\right) + \omega^6 \wedge \left(\frac{n}{2} \omega^5 - \frac{3m}{2} \omega^7\right),\nonumber \\
   d\omega^7= n \omega^5 \wedge \omega^7,\nonumber
   \end{gather}
 where $m,n,B \in C^\infty(\Sigma_7)$,
 \begin{gather}
  dm = m_5 \omega^5 + m_7 \omega^7, \qquad
  dn = n_5 \omega^5 + n_7 \omega^7, \qquad m_{57} = \Parder{m_5}{\omega^7} = \Parder{}{\omega^7} \left( \Parder{m}{\omega^5} \right), \quad etc.\nonumber\\
 dB = \epsilon \left( - 4m\Delta_1 - 2n\epsilon B - 6 mm_5 - m n_7 + nm_7 + \frac{3}{2} m_{57} + \frac{1}{2} n_{77}\right)\omega^5 \nonumber\\
\phantom{dB =}{} + \left(4n\Delta_1 + 2mB - 6 n n_7 - n m_5 + mn_5 - \frac{1}{2} m_{55} - \frac{3}{2} n_{75}\right) \omega^7
 \label{gamma425}
 \end{gather}
and
 \begin{gather*}
  \gam{4}{27} = mn + \epsilon - \frac{3}{2} n_7 - \frac{1}{2} m_5, \qquad
  \gam{6}{35} = mn + \epsilon + \frac{1}{2} n_7 + \frac{3}{2} m_5.
 \end{gather*}
 Finally, the integrability conditions for \eqref{Vranceanu-red-coframe} (i.e.\ $d^2\omega^i=0$ for all $i$) reduce to the integrability conditions for $dm$, $dn$, $dB$ as given above.
 \end{theorem}

 \begin{remark}
 All structures admitting a 9-dimensional symmetry group are included in \eqref{Vranceanu-red-coframe} (since $K^0$ cannot be reduced to the identity).
 \end{remark}

 \begin{remark}
 For all valid structures arising from \eqref{Vranceanu-red-coframe}, the function $\Delta_3 = B := \gam{4}{25}$ is a~relative contact invariant: it is $G^+$-invariant, and under the $R$-action, $\hat\Delta_3 = \epsilon \Delta_3$.
 \end{remark}

 \begin{corollary} \label{HK-reduced} For all valid structures arising from \eqref{Vranceanu-red-coframe}, the original $G$-structure on $\Sigma_7$ can be reduced to an $H$-structure, where $H = H^0 \rtimes D_8$ and
 \begin{gather}
 H^0 = \left\{ \mat{ccccccc}{
   a_1{}^2 & 0 & 0 & 0 & 0 & 0 & 0\\
   a_1 a_2 & a_1 & 0 & 0 & 0 & 0 & 0\\
   m a_1 a_2 & 0 & a_1 & 0 & 0 & 0 & 0\\
   0 & 0 & 0 & a_1 & 0 & 0 & 0\\
   0 & 0 & 0 & \epsilon m a_2 & 1 & 0 & 0\\
   0 & 0 & 0 & 0 & 0 & a_1 & 0\\
   0 & 0 & 0 & 0 & 0 & a_2 & 1
   } : (a_1,a_2) \in \R^+ \times \R \right\}.
   \label{H0-group}
 \end{gather}
 Moreover, wherever $\Delta_1 \neq 0$, or $B \neq 0$ there is a further reduction to a $K$-structure, where $K = K^0 \rtimes D_8$.
\end{corollary}

 \begin{proof}
 For all valid structures satisfying \eqref{Vranceanu-red-coframe}, $\gam{5}{56} = \gam{7}{47} = 0$, so from the $G$-action described in \eqref{gam-main-transform}, the stabilizer $G_\Gamma$ of $\Gamma = (\gam{5}{56},\gam{7}{47},\gam{5}{57},\gam{7}{57}) = (0,0,m,n)$ is contained in $H$ (since we can always keep $\hgam{5}{56}=0$ using $a_3 = \epsilon m a_2$).

 If $\Delta_1 \neq 0$, then $a_2=a_3=0$ is the unique solution to \eqref{mn-linear-sys} and $G_\Gamma \subset K$.
 Alternatively, suppose $B \neq 0$.  Note that $\hgam{4}{15}$ and $\hgam{6}{17}$ are torsion coef\/f\/icients, and for the structure equations~\eqref{Vranceanu-red-coframe}, we have $\gam{4}{15}=\gam{6}{17} = 0$, and the transformation laws (under $H^0$):
 \begin{gather*}
 \hgam{4}{15} = \frac{-Ba_2}{a_1}, \qquad \hgam{6}{17} = \frac{-B\epsilon m a_2}{a_1}.
 \end{gather*}
 Consequently, we can normalize $\hgam{4}{15} = \hgam{6}{17} = 0$ and reduce the connected component of the structure group to $K^0$ by setting $a_2=0$.  The discrete part of the structure group will preserve this reduction since
 \begin{gather*}
 R\mbox{-action}: \quad \hgam{4}{15} = -\gam{6}{17}, \quad \hgam{6}{17} = \gam{4}{15},\\
 S\mbox{-action}: \quad \hgam{4}{15} = -\gam{4}{15}, \quad \hgam{6}{17} = \gam{6}{17}.\tag*{\qed}
 \end{gather*}\renewcommand{\qed}{}
 \end{proof}

  Let us now examine in detail the case when $m$, $n$ are {\em constants}.  Then \eqref{gamma425} becomes
 \begin{gather}
 dB = -2\left( 2\epsilon m\Delta_1 + n B \right)\omega^5 + 2\left(2n\Delta_1 + mB \right) \omega^7.
 \label{gamma425-mn-const}
 \end{gather}
 Applying $d$ to \eqref{gamma425-mn-const} and simplifying, we obtain the integrability condition
 \begin{gather*}
 0 = -12\epsilon\Delta_1 \Delta_2 \omega^5 \wedge \omega^7.
 \end{gather*}

 \begin{corollary} Suppose $m$, $n$ are constants.  Then $\Delta_1 \Delta_2 = 0$ if and only if \eqref{Vranceanu-red-coframe} are valid structure equations.  Moreover, in this case:
 \begin{enumerate}\itemsep=0pt
 \item[1)] $\sigma= -n\omega^5 + m\omega^7$ is closed, so $\sigma = dh$ for some function $h \in C^\infty(\Sigma_7)$;
 \item[2)] $m=0$ iff $n=0$ iff $\sigma=0$ iff $h$ is constant;
 \item[3)] if $\Delta_1 = 0$, then $n = -\frac{\epsilon}{m}$, and $dB = 2B \sigma$, so $B = be^{2h}$, where $b$ is an arbitrary constant;
 \item[4)] if $\Delta_2 = 0$, then:
	\begin{itemize}\itemsep=0pt
	\item if $\epsilon = -1$, then $m=n=0$, and $B$ is an arbitrary constant;
	\item if $\epsilon = 1$, then letting $n = \epsilon_1 m$, $\epsilon_1=\pm 1$ we have $dB = 2( 2(m^2 + \epsilon_1) + B)\sigma$.  If $m\neq 0$, then $B = -2(m^2+\epsilon_1) + be^{2h}$, where $b$ is an arbitrary constant.  If $m=0$, then $B$ is an arbitrary constant;
	\end{itemize}
 \item[5)] If $\Delta_1 = \Delta_2 = 0$, then $\epsilon=1$, and $(m,n)=(1,-1)$ or $(-1,1)$.
 \end{enumerate}
 All of the above structures have a symmetry group with dimension at least seven.
 \end{corollary}

 \begin{proof} We prove only the f\/inal assertion as the others are straightforward to prove.  Let $G_\Gamma$~be the reduced structure group for which there is no group dependent torsion.  By construction (c.f.\ Theo\-rem~\ref{mn-str-eqns}), we must have $K^0 \subset G_\Gamma$, and by Proposition \ref{prop:prolong} we prolong to an $\{e\}$-structure on~$\Sigma_7 \times G_\Gamma$.  If $B$ is constant, then by Theorem \ref{equiv-soln} the symmetry group has dimension $\dim(\Sigma_7 \times G_\Gamma)  \geq  8$.  If $B$ is nonconstant, then by Corollary \ref{HK-reduced}, $G_\Gamma \subset K$.  Note that $\hat{B} = B$, so equation \eqref{gamma425-mn-const} implies that on $\Sigma_7 \times G_\Gamma$, we have
 \begin{gather*}
  dB = -2\left( 2\epsilon m\Delta_1 + n B \right) \hat\omega^5 + 2\left(2n\Delta_1 + mB \right) \hat\omega^7.
 \end{gather*}
 Thus, the coframe derivatives of $B$ are functions of $B$.  Thus, if $B$ is nonconstant, then the rank of the lifted coframe $\bm\Theta$ is 1 and by Theorem \ref{equiv-soln} the symmetry group will be at least $\dim(\Sigma_7 \times G_\Gamma) - {\rm rank}(\bm\Theta) \geq 8-1=7$ dimensional.
 \end{proof}

 \begin{remark}
 In the case $\Delta_2=0$, $\epsilon = 1$, we note that $\epsilon_1$ is a contact invariant.
 \end{remark}

 Certain values of $m$, $n$, $B$ lead to equivalent structures owing to the presence of the $D_8$ discrete subgroup of the original structure group $G$.
 Suppose $\Delta_1 = 0$, so $n = -\frac{\epsilon}{m}$.  Then
 \begin{gather*}
 R\mbox{-action}: \quad \hat{m} = -\frac{1}{m}, \quad \hat{B} = \epsilon B,\\
 S\mbox{-action}: \quad \hat{m} = -m, \quad \hat{B} = B.
 \end{gather*}
 In this case, by choosing a representative element $m \in (0,1]$, we can reduce $D_8$ to $\Z_2 = \langle R^2 \rangle$.  If $\epsilon = 1$, no further reduction occurs.  If $\epsilon = -1$ and $B \neq 0$, we choose a representative out of $\{ B, -B \}$ to reduce the discrete subgroup to the identity.  A similar argument is used in the case $\Delta_1 \neq 0$, where $\Delta_2=0$, $n=\epsilon_1 m$, and
 \begin{gather*}
 R\mbox{-action}: \quad \hat{m} = \epsilon \epsilon_1 m, \quad \hat{B} = \epsilon B,\\
 S\mbox{-action}: \quad \hat{m} = -m, \quad \hat{B} = B.
 \end{gather*}

 The results are organized in Table \ref{streqn-classification}
 according to the dimension of the symmetry group of the resulting $\{e\}$-structures on $\Sigma_7 \times G_\Gamma$.

 \begin{table}[h]
 \centering
  \caption{All generic hyperbolic structures for which $m,n$ are constants and $K^0 \subset G_\Gamma$.}
 \label{streqn-classification}
 \vspace{-5mm}
 \begin{align*}
 \begin{array}{|c|c|c|c|c|c|c|c|} \hline
 \mbox{Sym.\ grp.} & \Delta_1 & \Delta_2 & (\epsilon, m) & n & B & \mbox{Str.\ grp. } G_\Gamma\\ \hline\hline
 9 & 0 & \neq 0 & \{\pm 1\} \times (0,1] & -\frac{\epsilon}{m} & 0 & H^0 \rtimes \langle R^2 \rangle\\
  & & & \mbox{except } (1,1) & & & \\
 9 & 0 & 0 & (1,1) & -1 & 0 & H^0 \rtimes \langle R^2 \rangle\\ \hline\hline
 8 & \neq 0 & 0 & (-1,0) & 0 & b > 0 & K^0 \rtimes \langle R^2,S \rangle\\
 8 & \neq 0 & 0 & (-1,0) & 0 & 0 & K^0 \rtimes D_8\\
 8 & \neq 0 & 0 & (1,0) & 0 & b \in \R & K^0 \rtimes D_8 \\ \hline
 8 & \neq 0 & 0 & \{1\} \times (0,\infty) & m & -2(m^2+1) & K^0 \rtimes \langle R \rangle\\
 8 & \neq 0 & 0 & \{1\} \times (0,\infty) & -m & -2(m^2-1) & K^0 \rtimes \langle R^2 \rangle\\ \hline\hline
 7 & 0 & \neq 0 & \{-1\} \times (0,1] & \frac{1}{m} & be^{2h},\ b > 0 & K^0\\
 7 & 0 & \neq 0 & \{1\} \times (0,1) & -\frac{1}{m} & be^{2h},\ b \in \R^\times & K^0 \rtimes \langle R^2 \rangle\\ \hline
 7 & 0 & 0 & (1,1) & -1 & be^{2h},\ b \in \R^\times & K^0 \rtimes \langle R^2 \rangle\\ \hline
 7 & \neq 0 & 0 & \{1\} \times (0,\infty) & m & -2(m^2+1) + be^{2h},\ b \in \R & K^0 \rtimes \langle R \rangle\\
 7 & \neq 0 & 0 & \{1\} \times (0,\infty) & -m & -2(m^2-1) + be^{2h},\ b \in \R & K^0 \rtimes \langle R^2 \rangle\\ \hline
 \end{array}
 \end{align*}
 ($h$ is a nonconstant function such that $dh = -n\omega^5 + m\omega^7$)
 \end{table}

 \begin{remark}
 Vranceanu explicitly derived the following constant torsion cases:
 \begin{itemize}\itemsep=0pt
 \item 9-dim. symmetry: $\epsilon = 1$, $\Delta_1=0$, $B=0$;
 \item 8-dim. symmetry:
 \begin{enumerate}\itemsep=0pt
 \item[1)] $\epsilon = 1$, $\Delta_1\neq 0$, $\Delta_2 =0$, $m=n=0$,
 \item[2)] $\epsilon = 1$, $\Delta_1\neq 0$, $\Delta_2 =0$, $n=\pm m$, $B = -2(m^2\pm 1)$.
 \end{enumerate}
 \end{itemize}
 \end{remark}

 \begin{theorem}
 All contact-inequivalent generic hyperbolic structures for which:{\samepage
 \begin{enumerate}\itemsep=0pt
 \item[1)] $K^0$ is a subgroup of the structure group, and
 \item[2)] $m$, $n$ are constants,
 \end{enumerate}
 are displayed in Table~{\rm \ref{streqn-classification}}.}
 \end{theorem}

 For ease of reference, we state below the structure equations explicitly for each of the cases above.  For the maximally symmetric cases, we state the structure equations for both the base coframe $\{ \omega^1, \dots, \omega^7 \}$ and the lifted coframe on  $\Sigma_7 \times G_\Gamma$.  In the submaximally symmetric cases, we only display structure equations for the lifted coframe.  (One can obtain the structure equations on the base simply by setting $\hat\alpha^1 = 0$ and removing all hats from the remaining variables.)
 In each case, we assume that $G_\Gamma$ and all parameters are as in Table~\ref{streqn-classification}.  Note that $d\hat\omega^i$ are determined by~\eqref{gen-lifted-streqns}.  Following potentially some Lie algebra valued compatible absorption, $\hat\alpha^\rho = \alpha^\rho + \lambda^\rho{}_i \hat\omega^i$, the structure equations $d\hat\alpha^\rho$ are determined by the integrability conditions $d^2\hat\omega^i=0$.  (We only display the f\/inal results.)

 For those coframes whose structure equations depend explicitly on the (nonconstant) function~$h$, we have $m\neq 0$ (c.f.\ Table~\ref{streqn-classification}) and the symmetry algebra is determined by restricting to the level set $h=h_0$, where $h_0$ is a constant.  (Note: We will abuse notation and identify $h \in C^\infty(\Sigma_7)$ with its pullback to the bundle.)  On this level set, we have $0 = dh = -n \hat\omega^5 + m\hat\omega^7$.  We can choose (the pullback of) $\{ \hat\omega^1, \dots, \hat\omega^6, \hat\alpha^1 \}$ as a coframe on each level set, and the corresponding structure equations will have constant coef\/f\/icients.  Thus, these are Maurer--Cartan equations for a local Lie group.  A well-known fact is that the isomorphism type of the symmetry algebra of a coframe determined in this way is independent of the level set chosen.  Consequently, we make the canonical choice and restrict to the level set $h=0$ in these cases.

 The structure constants for the (contact) symmetry algebra for each of the structures can be read of\/f from the structure equations for the coframe (or its pullback to the level set $h=0$ if~$h$ appears explicitly).  Only the symmetry algebras appearing in the 9-dimensional case will be studied in further detail in this article.

 \subsection[Case 1: $\Delta_1 = 0$, $B = 0$]{Case 1: $\boldsymbol{\Delta_1 = 0}$, $\boldsymbol{B = 0}$}

 This branch consists of precisely all maximally symmetric generic hyperbolic equations.

 Parameters: \qquad $(\epsilon,m) \in \{\pm 1\} \times (0,1]$.

 Base coframe:
 \begin{gather}
 d\omega^1 = \omega^2 \wedge \omega^4 + \omega^3 \wedge \omega^6,\nonumber\\
 d\omega^2 = \omega^4 \wedge \omega^5 + \omega^3 \wedge \omega^7 + \omega^2 \wedge \left( \frac{3\epsilon}{2m} \omega^5 + \frac{m}{2} \omega^7 \right),\nonumber\\
 d\omega^3 = \omega^6 \wedge \omega^7 + \epsilon\omega^2 \wedge \omega^5 + \omega^3 \wedge \left( \frac{\epsilon}{2m} \omega^5 + \frac{3m}{2} \omega^7 \right),\nonumber\\
 d\omega^4 = \epsilon \omega^5 \wedge \omega^6 - \omega^4 \wedge \left( \frac{3\epsilon}{2m} \omega^5 + \frac{m}{2} \omega^7 \right),  \label{9dim-streqns}\\
 d\omega^5 = m \omega^5 \wedge \omega^7,\nonumber\\
 d\omega^6 = - \omega^4 \wedge \omega^7 - \omega^6 \wedge \left( \frac{\epsilon}{2m} \omega^5 + \frac{3m}{2} \omega^7 \right),\nonumber\\
 d\omega^7 = -\frac{\epsilon}{m} \omega^5 \wedge \omega^7.\nonumber
 \end{gather}

 Lifted coframe on $\Sigma_7 \times G_\Gamma$:
 \begin{gather*}
 d\hat\omega^1 = 2\hat\alpha^1 \wedge \hat\omega^1 + \hat\omega^2 \wedge \hat\omega^4 + \hat\omega^3 \wedge \hat\omega^6,\\
 d\hat\omega^2 = \hat\alpha^2 \wedge \hat\omega^1 + \hat\alpha^1 \wedge \hat\omega^2 + \hat\omega^4 \wedge \hat\omega^5 + \hat\omega^3 \wedge \hat\omega^7 + \hat\omega^2 \wedge \left( \frac{3\epsilon}{2m} \hat\omega^5 + \frac{m}{2} \hat\omega^7 \right),\\
 d\hat\omega^3 = m \hat\alpha^2 \wedge \hat\omega^1 + \hat\alpha^1 \wedge \hat\omega^3 + \hat\omega^6 \wedge \hat\omega^7 + \epsilon\hat\omega^2 \wedge \hat\omega^5 + \hat\omega^3 \wedge \left( \frac{\epsilon}{2m} \hat\omega^5 + \frac{3m}{2} \hat\omega^7 \right),\\
 d\hat\omega^4 = \hat\alpha^1 \wedge \hat\omega^4 + \epsilon \hat\omega^5 \wedge \hat\omega^6 - \hat\omega^4 \wedge \left( \frac{3\epsilon}{2m} \hat\omega^5 + \frac{m}{2} \hat\omega^7 \right), \\
 d\hat\omega^5 = \epsilon m \hat\alpha^2 \wedge \hat\omega^4 + m \hat\omega^5 \wedge \hat\omega^7,\\
 d\hat\omega^6 = \hat\alpha^1 \wedge \hat\omega^6 - \hat\omega^4 \wedge \hat\omega^7 - \hat\omega^6 \wedge \left( \frac{\epsilon}{2m} \hat\omega^5 + \frac{3m}{2} \hat\omega^7 \right),\\
 d\hat\omega^7 = \hat\alpha^2 \wedge \hat\omega^6 - \frac{\epsilon}{m} \hat\omega^5 \wedge \hat\omega^7,\\
 d\hat\alpha^1 = \frac{1}{2} \hat\alpha^2 \wedge (\hat\omega^4 + m \hat\omega^6), \\
 d\hat\alpha^2 = \hat\alpha^2 \wedge \left(\hat\alpha^1 + \frac{3}{2} \left(\frac{\epsilon}{m} \hat\omega^5 + m \hat\omega^7\right)\right).
 \end{gather*}

 \subsection[Case 2: $\Delta_2 = 0$, $B$ constant]{Case 2: $\boldsymbol{\Delta_2 = 0}$, $\boldsymbol{B}$ constant}
 \label{8d-str-eqs}

 This branch contains two families of equations with 8-dimensional symmetry.  All coef\/f\/icients in both sets of structure equations are constants.

 \subsubsection[Case 2a: $m=n=0$]{Case 2a: $\boldsymbol{m=n=0}$}
\vspace{-5mm}

  \begin{gather*}
   d\hat\omega^1 = 2 \hat\alpha^1 \wedge \hat\omega^1 + \hat\omega^3 \wedge \hat\omega^6 + \hat\omega^2 \wedge \hat\omega^4, \\
   d\hat\omega^2 = \hat\alpha^1 \wedge \hat\omega^2 + \hat\omega^4 \wedge \hat\omega^5 + \hat\omega^3 \wedge \hat\omega^7, \\
   d\hat\omega^3 = \hat\alpha^1 \wedge \hat\omega^3 + \hat\omega^6 \wedge \hat\omega^7 + \epsilon \hat\omega^2 \wedge \hat\omega^5, \\
   d\hat\omega^4 = \hat\alpha^1 \wedge \hat\omega^4 + \epsilon \hat\omega^5 \wedge \hat\omega^6 + b\hat\omega^2 \wedge \hat\omega^5 +\epsilon \hat\omega^2 \wedge \hat\omega^7, \\
   d\hat\omega^5 = 0, \\
   d\hat\omega^6 = \hat\alpha^1 \wedge \hat\omega^6 - \hat\omega^4 \wedge \hat\omega^7 + \epsilon \hat\omega^3 \wedge \hat\omega^5 + \epsilon b \hat\omega^3 \wedge \hat\omega^7, \\
   d\hat\omega^7 = 0, \\
   d\hat\alpha^1 = 0.
   \end{gather*}

 \subsubsection[Case 2b: $n=\epsilon_1 m \neq 0$ (and $\epsilon = 1$)]{Case 2b: $\boldsymbol{n=\epsilon_1 m \neq 0}$ (and $\boldsymbol{\epsilon = 1}$)}
\vspace{-5mm}

 \begin{gather*}
   d\hat\omega^1 =  2 \hat\alpha^1 \wedge \hat\omega^1 + \hat\omega^3 \wedge \hat\omega^6 + \hat\omega^2 \wedge \hat\omega^4, \\
   d\hat\omega^2 = \hat\alpha^1 \wedge \hat\omega^2 + \hat\omega^4 \wedge \hat\omega^5 + \hat\omega^3 \wedge \hat\omega^7 - \hat\omega^2 \wedge \left( \frac{3\epsilon_1 m}{2} \hat\omega^5 - \frac{m}{2} \hat\omega^7 \right),\\
   d\hat\omega^3 = \hat\alpha^1 \wedge \hat\omega^3 + \hat\omega^6 \wedge \hat\omega^7 + \epsilon \hat\omega^2 \wedge \hat\omega^5 - \hat\omega^3 \wedge \left( \frac{\epsilon_1 m}{2} \hat\omega^5 - \frac{3m}{2} \hat\omega^7 \right),\\
   d\hat\omega^4 = \hat\alpha^1 \wedge \hat\omega^4 + \epsilon \hat\omega^5 \wedge \hat\omega^6 + (m^2+\epsilon_1) \hat\omega^2 \wedge \left(-2 \hat\omega^5 +\epsilon_1 \hat\omega^7\right) + \hat\omega^4 \wedge \left(\frac{3\epsilon_1 m}{2} \hat\omega^5 - \frac{m}{2} \hat\omega^7\right), \\
   d\hat\omega^5 = m \hat\omega^5 \wedge \hat\omega^7, \\
   d\hat\omega^6 = \hat\alpha^1 \wedge \hat\omega^6 - \hat\omega^4 \wedge \hat\omega^7 + (m^2+\epsilon_1) \hat\omega^3 \wedge \left( \epsilon_1 \hat\omega^5 - 2 \hat\omega^7\right) + \hat\omega^6 \wedge \left(\frac{\epsilon_1 m}{2} \hat\omega^5 - \frac{3m}{2} \hat\omega^7\right), \\
   d\hat\omega^7 = \epsilon_1 m \hat\omega^5 \wedge \hat\omega^7,\\
   d\hat\alpha^1 = 0.
   \end{gather*}

 \subsection[Case 3: $B$ nonconstant]{Case 3: $\boldsymbol{B}$ nonconstant}
 \label{7d-str-eqs}
 This branch contains two families of equations with 7-dimensional symmetry.  Note the case $\Delta_1=\Delta_2=0$, $\epsilon=m=-n=1$ is contained in both families.

 \subsubsection[Case 3a: $\Delta_1 = 0$, $B$ nonconstant]{Case 3a: $\boldsymbol{\Delta_1 = 0}$, $\boldsymbol{B}$ nonconstant}
\vspace{-5mm}

 \begin{gather*}
 d\hat\omega^1 = 2\hat\alpha^1 \wedge \hat\omega^1 + \hat\omega^2 \wedge \hat\omega^4 + \hat\omega^3 \wedge \hat\omega^6, \\
 d\hat\omega^2 = \hat\alpha^1 \wedge \hat\omega^2 + \hat\omega^4 \wedge \hat\omega^5 + \hat\omega^3 \wedge \hat\omega^7 + \hat\omega^2 \wedge \left( \frac{3\epsilon}{2m} \hat\omega^5 + \frac{m}{2} \hat\omega^7 \right), \\
 d\hat\omega^3 = \hat\alpha^1 \wedge \hat\omega^3 + \hat\omega^6 \wedge \hat\omega^7 + \epsilon\hat\omega^2 \wedge \hat\omega^5 + \hat\omega^3 \wedge \left( \frac{\epsilon}{2m} \hat\omega^5 + \frac{3m}{2} \hat\omega^7 \right),\\
 d\hat\omega^4 = \hat\alpha^1 \wedge \hat\omega^4 + \epsilon \hat\omega^5 \wedge \hat\omega^6 + be^{2h} \hat\omega^2 \wedge \hat\omega^5 - \hat\omega^4 \wedge \left(\frac{3\epsilon}{2m} \hat\omega^5 + \frac{m}{2} \hat\omega^7\right), \\
 d\hat\omega^5 = m \hat\omega^5 \wedge \hat\omega^7,\\
 d\hat\omega^6 = \hat\alpha^1 \wedge \hat\omega^6 - \hat\omega^4 \wedge \hat\omega^7 + \epsilon be^{2h} \hat\omega^3 \wedge \hat\omega^7 - \hat\omega^6 \wedge \left(\frac{\epsilon}{2m} \hat\omega^5 + \frac{3m}{2} \hat\omega^7\right), \\
 d\hat\omega^7 = - \frac{\epsilon}{m} \hat\omega^5 \wedge \hat\omega^7,\\
 d\hat\alpha^1 = 0.
 \end{gather*}

 On the level set $\{ h = 0 \}$: In this case, $\hat\omega^7 = -\frac{\epsilon}{m^2} \hat\omega^5$.
 \begin{gather*}
 d\hat\omega^1 = 2\hat\alpha^1 \wedge \hat\omega^1 + \hat\omega^2 \wedge \hat\omega^4 + \hat\omega^3 \wedge \hat\omega^6, \\
 d\hat\omega^2 = \hat\alpha^1 \wedge \hat\omega^2 + \left( \frac{\epsilon}{m} \hat\omega^2 - \frac{\epsilon}{m^2} \hat\omega^3 + \hat\omega^4 \right) \wedge \hat\omega^5, \\
 d\hat\omega^3 = \hat\alpha^1 \wedge \hat\omega^3 + \left( \epsilon\hat\omega^2 - \frac{\epsilon}{m} \hat\omega^3 - \frac{\epsilon}{m^2} \hat\omega^6 \right) \wedge \hat\omega^5,\\
 d\hat\omega^4 = \hat\alpha^1 \wedge \hat\omega^4 + \left(b \hat\omega^2 - \frac{\epsilon}{m} \hat\omega^4 - \epsilon \hat\omega^6 \right) \wedge \hat\omega^5, \\
 d\hat\omega^5 = 0,\\
 d\hat\omega^6 = \hat\alpha^1 \wedge \hat\omega^6 - \frac{1}{m^2} ( b \hat\omega^3 - \epsilon \hat\omega^4 - \epsilon m \hat\omega^6 ) \wedge \hat\omega^5, \\
 d\hat\alpha^1 = 0.
 \end{gather*}

 \subsubsection[Case 3b: $\Delta_2 = 0$, $B$ nonconstant]{Case 3b: $\boldsymbol{\Delta_2 = 0}$, $\boldsymbol{B}$ nonconstant}
\vspace{-5mm}

 \begin{gather*}
   d\hat\omega^1 =  2 \hat\alpha^1 \wedge \hat\omega^1 + \hat\omega^3 \wedge \hat\omega^6 + \hat\omega^2 \wedge \hat\omega^4, \\
   d\hat\omega^2 = \hat\alpha^1 \wedge \hat\omega^2 + \hat\omega^4 \wedge \hat\omega^5 + \hat\omega^3 \wedge \hat\omega^7 - \hat\omega^2 \wedge \left( \frac{3\epsilon_1 m}{2} \hat\omega^5 - \frac{m}{2} \hat\omega^7 \right),\\
   d\hat\omega^3 = \hat\alpha^1 \wedge \hat\omega^3 + \hat\omega^6 \wedge \hat\omega^7 + \hat\omega^2 \wedge \hat\omega^5 - \hat\omega^3 \wedge \left( \frac{\epsilon_1 m}{2} \hat\omega^5 - \frac{3m}{2} \hat\omega^7 \right),\\
   d\hat\omega^4 = \hat\alpha^1 \wedge \hat\omega^4 + \hat\omega^5 \wedge \hat\omega^6 + (-2(m^2+\epsilon_1)+be^{2h}) \hat\omega^2 \wedge \hat\omega^5 +(\epsilon_1 m^2+1) \hat\omega^2 \wedge \hat\omega^7 \\
   \phantom{d\hat\omega^4 =}{} + \hat\omega^4 \wedge \left(\frac{3\epsilon_1 m}{2} \hat\omega^5 - \frac{m}{2} \hat\omega^7\right), \\
   d\hat\omega^5 = m \hat\omega^5 \wedge \hat\omega^7, \\
   d\hat\omega^6 = \hat\alpha^1 \wedge \hat\omega^6 - \hat\omega^4 \wedge \hat\omega^7 + (\epsilon_1 m^2+1) \hat\omega^3 \wedge\hat\omega^5 +  (-2(m^2+\epsilon_1)+be^{2h}) \hat\omega^3 \wedge \hat\omega^7 \\
   \phantom{d\hat\omega^6 =}{}  + \hat\omega^6 \wedge \left(\frac{\epsilon_1 m}{2} \hat\omega^5 - \frac{3m}{2} \hat\omega^7\right), \\
   d\hat\omega^7 = \epsilon_1 m \hat\omega^5 \wedge \hat\omega^7,\\
   d\hat\alpha^1 = 0.
   \end{gather*}

  On the level set $\{ h = 0 \}$: In this case, $\hat\omega^7 = \epsilon_1 \hat\omega^5$.
 \begin{gather*}
   d\hat\omega^1 =  2 \hat\alpha^1 \wedge \hat\omega^1 + \hat\omega^3 \wedge \hat\omega^6 + \hat\omega^2 \wedge \hat\omega^4, \\
   d\hat\omega^2 = \hat\alpha^1 \wedge \hat\omega^2 + ( \hat\omega^4 + \epsilon_1 \hat\omega^3 -  \epsilon_1 m \hat\omega^2 ) \wedge \hat\omega^5,\\
   d\hat\omega^3 = \hat\alpha^1 \wedge \hat\omega^3 + \epsilon_1 (\hat\omega^6 + \epsilon_1 \hat\omega^2 +  m \hat\omega^3 ) \wedge \hat\omega^5 ,\\
   d\hat\omega^4 = \hat\alpha^1 \wedge \hat\omega^4 + (- \hat\omega^6 + (-(m^2+\epsilon_1)+b) \hat\omega^2 + \epsilon_1 m \hat\omega^4 ) \wedge \hat\omega^5, \\
   d\hat\omega^5 = 0,\\
   d\hat\omega^6 = \hat\alpha^1 \wedge \hat\omega^6 + \epsilon_1 (- \hat\omega^4 + (-(m^2+\epsilon_1)+b) \hat\omega^3  - m \hat\omega^6 ) \wedge \hat\omega^5, \\
   d\hat\alpha^1 = 0.
   \end{gather*}


 \section{The maximally symmetric case}
 \label{maxsym-case}
 \subsection{A coframing in local coordinates}

 For the remainder of the paper we focus on the maximally symmetric generic hyperbolic structures.  In Appendix~\ref{maxsym-param}, we  outline how Vranceanu arrived at an explicit coframe $\{ \omega^i \}_{i=1}^7$ on~$\Sigma_7$ given in local coordinates which satisf\/ies the structure equations \eqref{9dim-streqns}.  In local coordinates $(x,y,z,p,q,u,v)$ on $\Sigma_7$, the coframe is given by
 \begin{gather}
 \omega^1 = dz - pdx - qdy, \nonumber\\
\omega^2 = \left( \frac{\epsilon m^2}{6} - \frac{ m \alpha v^3}{3u^3} + \frac{\alpha v^2}{2u^2} \right)\omega^6 + \left( -\frac{\epsilon m}{3} - \frac{\alpha v^3}{3u^3} \right) \omega^4 - u^{-3/2} (dp + vdq),\nonumber\\
\omega^3 = \left( \frac{\epsilon m^2}{6} - \frac{ m \alpha v^3}{3u^3} + \frac{\alpha v^2}{2u^2} \right) \omega^4 + \left( -\frac{\epsilon m^3}{3} - \frac{ m^2 \alpha v^3}{3u^3} + \frac{ m \alpha v^2}{u^2} - \frac{\alpha v}{u} \right) \omega^6 \nonumber\\
 \phantom{\omega^3 =}{} -  m u^{-3/2} (dp + v dq) + u^{-1/2} dq, \label{9d-explicit-coframe}\\
 \omega^4 = u^{3/2} dx +  m \sqrt{u}(dy - vdx),
 \qquad \omega^5 = \frac{\epsilon  m( du -  m dv)}{u}, \nonumber\\
 \omega^6 = -\sqrt{u} (dy-vdx),
 \qquad \omega^7 =   \frac{dv}{u},\nonumber
 \end{gather}
 which is valid on the open set $u > 0$, and where $\alpha = 1 - \epsilon m^4$.
 The coordinates $(x,y,z,p,q)$ are identif\/ied with the corresponding coordinates on $J^1(\R^2,\R)$.

 Note that $\alpha$ in the case $\Delta_1=0$ is a relative contact invariant since $\alpha = -\epsilon m^2 \left( m^2 - \frac{\epsilon}{m^2} \right) = -\epsilon m^2 \Delta_2$, $m \neq 0$, and $\Delta_2$ is relative contact invariant.  Since the contact-inequivalent structures are parametrized by $(\epsilon, m) \in \{\pm 1\} \times (0,1]$, then $\alpha \in [0,1) \cup (1,2]$.

 \subsection{Normal forms}

 Let us determine how the coordinates $(u,v)$ on $\Sigma_7$ are related to the standard 2-jet coordinates $(x,y,z,p,q,r,s,t) \in J^2(\R^2,\R)$ .  Let $\chi : \R^2 \ra \Sigma_7$ be any integral manifold of $I_F$ with independence condition $\chi^*(dx \wedge dy) \neq 0$.  Without loss of generality, we identify the coordinates $(x,y)$ on $\R^2$ with the $(x,y)$ coordinates on $\Sigma_7$.  The composition $i_F \circ \chi$ is then an integral manifold of the contact system $\contact{2}$ and on $\R^2$ we can write
 \begin{gather}
 dp = r dx + s dy, \qquad dq = sdx + tdy, \label{dp-dq}
 \end{gather}
 where for convenience $p$ is identif\/ied with $(i_F \circ \chi)^* p$, and similarly for the coordinates $q$, $r$, $s$,~$t$.  Substituting \eqref{dp-dq} into the conditions $0 = \chi^* \omega^2 = \chi^* \omega^3$, and extracting the coef\/f\/icients of $dx$ and $dy$, we obtain the relations
 \begin{gather*}
 0 = 6vs+6r+2\epsilon m u^3-3\epsilon m^2u^2v-\alpha v^3, \\
 0 = 2vt+2s+\epsilon m^2u^2+\alpha v^2, \\
 0 = -6su+ m(6sv + 6r - \alpha v^3)+3v\epsilon m^3u^2+3\alpha v^2u -\epsilon m^2u^3, \\
 0 = -2tu+ m (2tv+2s+\alpha v^2)-\epsilon m^3u^2-2\alpha vu,
 \end{gather*}
 or equivalently, using the coordinate $w = u -  m v$ instead of $u$, we have
 \begin{gather}
  r = -\frac{1}{3} (\epsilon m w^3 + v^3 ),
 \qquad s = -\frac{1}{2} ( \epsilon m^2 w^2 - v^2),
 \qquad t = - (\epsilon m^3 w+ v). \label{rst-param}
 \end{gather}
 Thus, our PDE is of the form
 \begin{gather*}
 F(r,s,t) = 0,
 \end{gather*}
 and we have a nondegenerate parametrization $i_F : \Sigma_7 \ra J^2(\R^2,\R)$ (for $u = w + mv > 0$).

 Consider the case $\alpha = 0$, i.e.\ $(\epsilon, m) = (1,1)$.  In this case, it is straightforward to eliminate both parameters $w$, $v$ and obtain the equation
 \begin{gather}
 rt - s^2 - \frac{t^4}{12} = 0.
 \label{special-maxsym-eqn}
 \end{gather}

 Now consider the general case $\alpha \neq 0$.  Let us write $u = -\frac{1}{\tilde{u}}$, $v = \tilde{v}$ and rewrite \eqref{rst-param} as
 \begin{gather*}
 \tilde{u} t = \epsilon  m^3 - \alpha \tilde{u} \tilde{v}, \\
 \tilde{u}^2 s = -\frac{1}{2} \epsilon  m^2 - \epsilon  m^3 \tilde{u}\tilde{v} + \frac{1}{2} \alpha (\tilde{u} \tilde{v})^2
 = -\frac{\epsilon  m^2}{2\alpha} + \frac{\tilde{u}^2 t^2}{2\alpha},\\
 \tilde{u}^3 r = \frac{1}{3} \epsilon  m + \epsilon  m^2 \tilde{u} \tilde{v} + \epsilon  m^3 (\tilde{u} \tilde{v})^2 - \frac{1}{3} \alpha (\tilde{u} \tilde{v})^3
 = \frac{m(\epsilon+m^4)}{3\alpha^2} - \frac{\epsilon m^2 \tilde{u} t}{\alpha^2} + \frac{\tilde{u}^3 t^3}{3\alpha^2},
 \end{gather*}
 and so using $\nu = (\epsilon m^3 - \alpha \tilde{u} \tilde{v})^{-1}$ as a new parameter, we arrive at
 \begin{gather*}
 2\alpha s - t^2 = -\epsilon  m^2 \nu^2 t^2, \qquad
 3\alpha^2 r = m(\epsilon+m^4) \nu^3 t^3 - 3\epsilon m^2 \nu^2 t^3 + t^3.
 \end{gather*}
  Eliminating the parameter $\nu$, we obtain
 \begin{gather}
 (\epsilon+m^4)^2 (2\alpha s - t^2)^3 + \epsilon m^4 (3\alpha^2 r - 6\alpha st + 2t^3)^2 = 0. \label{gen-eqn-alpha}
 \end{gather}
 Finally, use the scaling $\bar{x} = \frac{1}{\alpha} x$, which induces
 \begin{gather*}
(\bar{r},\bar{s},\bar{t}) = \left( \alpha^2 r, \alpha s, t\right)
 \end{gather*}
 to eliminate $\alpha$ from \eqref{gen-eqn-alpha}.  Dropping bars, and letting $a=m^4$ we obtain
 \begin{gather}
 (\epsilon + a)^2 \left(2 s - t^2 \right)^3 + \epsilon a \left( 3r - 6st + 2t^3 \right)^2 = 0. \label{general-maxsym-eqn}
 \end{gather}
 Note that in the case $\epsilon=1$ considered by Vranceanu, the $st$ term has a missing factor of 2.

 \begin{theorem} \label{param-maxsym-eqn} The contact-equivalence classes of maximally symmetric generic hyperbolic PDE are parametrized by $(\epsilon, a) \in \{ \pm 1 \} \times (0,1]$.  Normal forms from each equivalence class are given by \eqref{special-maxsym-eqn} in the case $(\epsilon,a)=(1,1)$ and \eqref{general-maxsym-eqn} otherwise.
 \end{theorem}

 \begin{remark}
 Letting $\epsilon=a=1$ in \eqref{general-maxsym-eqn}, we have $F= 4 \left(2 s - t^2 \right)^3 + \left( 3r - 6st + 2t^3 \right)^2 = 0$, and
 \begin{gather*}
 \Delta = F_r F_t - \frac{1}{4} F_s{}^2 = -36 (2s - t^2) F.
 \end{gather*}
 On the equation manifold (and hence on $\Sigma_7$), we have $\Delta=0$ and consequently, this limiting equation is parabolic.
 \end{remark}


 \subsection{Nine-dimensional symmetry algebras}

 The calculations leading to \eqref{rst-param} are quite long and consequently to conf\/irm the validity of~\eqref{rst-param} (and in turn,  Theorem \ref{param-maxsym-eqn}), it is useful to describe the nine-dimensional (contact) symmetry algebra explicitly for the normal forms given in the previous section.  Calculating the symmetry algebra is a nontrivial task however -- the standard Lie method of calculating symmetries (by working in $J^2(\R^2,\R)$ on the equation locus) is highly impractical owing to the complexity of the equations.  In Appendix~\ref{9d-sym-alg}, we describe how the symmetry algebra was found by an alternative method.
 In order to give a unif\/ied description of the symmetry algebras, we work with the normal forms~\eqref{special-maxsym-eqn} and~\eqref{gen-eqn-alpha} as these arise from the parametrization~\eqref{rst-param}.

 \begin{proposition} \label{9d-syms} Any equation of the form $F(r,s,t)=0$ admits the symmetries
 \begin{gather*}
  X_1 = \parder{x}, \qquad X_2 =\parder{y}, \qquad X_3 =\parder{z}, \qquad X_4 =x\parder{z}, \qquad X_5 =y\parder{z},\\
  X_6=x\parder{x} + y\parder{y} + 2z\parder{z}.
 \end{gather*}
 The equations \eqref{special-maxsym-eqn} and \eqref{gen-eqn-alpha} have the following additional symmetries:
 \begin{gather*}
 X_7 = y\parder{y} + 3z\parder{z},
 \qquad X_8 = x \parder{y} - \frac{\alpha}{2} y^2 \parder{z},
\qquad X_9= x^2 \parder{x} + xy \parder{y} + \left(xz-\frac{\alpha}{6} y^3\right) \parder{z}.
 \end{gather*}
 In particular, all of these symmetries are projectable point symmetries.
 \end{proposition}

 (Recall that a {\em point} symmetry here is a vector f\/ield on $J^0(\R^2,\R)$.  A point symmetry is {\em projectable} if it projects to a vector f\/ield on the base $\R^2$.)

 The normalization of \eqref{gen-eqn-alpha} to \eqref{general-maxsym-eqn} is carried out by letting $\bar{x} = \frac{1}{\alpha} x$ from which we get:

 \begin{corollary} \label{cor:general-maxsym-eqn}
 The generic hyperbolic equation \eqref{general-maxsym-eqn} has symmetry generators $X_1,\dots, X_6$ as in Proposition~{\rm \ref{9d-syms}} as well as
 \begin{gather*}
  X_7 = y\parder{y} + 3z\parder{z},
\qquad X_8 = x \parder{y} - \frac{1}{2} y^2 \parder{z},
\qquad X_9= x^2 \parder{x} + xy \parder{y} + \left(xz-\frac{1}{6} y^3\right) \parder{z}.
 \end{gather*}
 \end{corollary}

 We will denote the corresponding abstract Lie algebras as $\g_\alpha$ and express their commutator relations in a canonical basis.  Let
 \begin{gather*}
 (e_1,e_2,e_3,e_4,e_5,e_6,e_7,e_8,e_9) = (X_2,X_3,X_4,X_5,-X_8,X_7,X_1,-2X_6+X_7,-X_9).
 \end{gather*}
 The commutator relations in this basis are
 \begin{gather*}
 \begin{array}{c|cccccc|ccc}
 & e_1 & e_2 & e_3 & e_4 & e_5 & e_6 & e_7 & e_8 & e_9 \\ \hline
 e_1 & \cdot & \cdot & \cdot & e_2 & \alpha e_4 & e_1 & \cdot & -e_1 & e_5 \\
 e_2 & & \cdot & \cdot & \cdot  & \cdot & 3 e_2 & \cdot & -e_2 & -e_3 \\
 e_3 & & & \cdot & \cdot & \cdot & 3e_3 & -e_2 & e_3 & \cdot \\
 e_4 & & & & \cdot & e_3 & 2e_4 & \cdot & \cdot & \cdot \\
 e_5 & & & & & \cdot & e_5 & e_1 & e_5 & \cdot \\
 e_6 & & & & & & \cdot & \cdot & \cdot &\cdot \\ \hline
 e_7 & & & & & & & \cdot & -2e_7 & e_8\\
 e_8 & & & & & & & & \cdot & -2e_9\\
 e_9 & & & & & & & & & \cdot
 \end{array}
 \end{gather*}
 In the case $\alpha \neq 0$, redef\/ining
 \begin{gather*}
 (\bar{e}_2, \bar{e}_3, \bar{e}_4 ) =  (\alpha e_2, \alpha e_3, \alpha e_4)
 \end{gather*}
 and dropping the bars, we have the same commutator relations as above except $\alpha$ has been normalized to~1.  Thus, in the case $\alpha \neq 0$, {\em all symmetry algebras are isomorphic}.  (This is also obvious from the fact that the symmetry generators in Corollary \ref{cor:general-maxsym-eqn} are independent of~$\alpha$.)

 Let $\g_1$ denote the abstract symmetry algebra in the case $\alpha \neq 0$, although this is a slight abuse of notation since $\alpha \in (0,1) \cup (1,2]$ in this case.  We calculate for $\g = \g_\delta$ ($\delta=0,1$),
 \begin{alignat*}{3}
& \text{Killing form:} && \kappa = {\rm diag}\left(0,0,0,0,0,24,\left( \begin{array}{ccc} 0 & 0 & 6\\ 0 & 12 & 0\\ 6 & 0 & 0 \end{array}\right)\right),&\\
& \text{derived subalgebra:} & & \g^{(1)} = \langle e_1, e_2, e_3, e_4, e_5, e_7, e_8, e_9 \rangle, & \\
& \text{radical:} &&  \gothic{r} = (\g^{(1)})^{\perp_\kappa} = \langle e_1, e_2, e_3, e_4, e_5, e_6 \rangle, &\\
& \text{(semi-simple) Levi factor:}\quad && \g_{ss} = \langle e_7, e_8, e_9 \rangle \cong \gothic{sl}(2,\R), &\\
 & \text{Levi decomposition:} & & \g = \gothic{r} \rtimes \g_{ss}, &\\
& nilradical: &&  \gothic{n} = \langle e_1, e_2, e_3, e_4, e_5 \rangle, & \\
& \text{derived series of} \ \gothic{r}: && \gothic{r}^{(1)} = \gothic{n}, \quad \gothic{r}^{(2)} = \langle e_2, e_3, \delta e_4 \rangle, \quad
 \gothic{r}^{(\infty)} = \gothic{r}^{(3)} = 0, & \\
& \text{lower central series of} \ \gothic{r}: && \gothic{r}^\infty = \gothic{r}^{1} = \gothic{n}. &
 \end{alignat*}

 An isomorphism between two Lie algebras must restrict to an isomorphism of their radicals and the corresponding derived f\/lags of the radicals.  Since $\gothic{r}^{(2)}$ is two-dimensional for $\g_0$ and three-dimensional for $\g_1$, then we must have $\g_0 \not\cong \g_1$.

  \begin{theorem} The contact symmetry algebra of any maximally symmetric generic hyperbolic PDE is:
 \begin{enumerate}\itemsep=0pt
 \item[1)] nine-dimensional,
 \item[2)] contact-equivalent to a (projectable) point symmetry algebra.
 \end{enumerate}
 Moreover, there are exactly two isomorphism classes of Lie algebras (represented by $\g_0$ and $\g_1$) that arise as such symmetry algebras.
 \end{theorem}

 We remark that Mubarakzjanov has classif\/ied all f\/ive-dimensional real solvable Lie algebras (labelled by~$g_{5,*}$)~\cite{Mubar5-1963} and all six-dimensional non-nilpotent real solvable Lie algebras (labelled by~$g_{6,*}$)~\cite{Mubar6-1963}.  The nilradicals of $\g_0$ and $\g_1$ can be identif\/ied in the former classif\/ication as:
 \begin{gather*}
 \gothic{n}_0 \cong g_{5,1}: \quad (\bar{e}_1,\bar{e}_2,\bar{e}_3,\bar{e}_4,\bar{e}_5) = (e_2,-e_3,e_1,e_5,e_4),\\
 \gothic{n}_1 \cong g_{5,3}: \quad (\bar{e}_1,\bar{e}_2,\bar{e}_3,\bar{e}_4,\bar{e}_5) = (e_3,e_4,-e_2,e_1,e_5).
 \end{gather*}
The radicals of $\g_0$ and $\g_1$ can be identif\/ied in the latter classif\/ication as
\begin{gather*}
 \gothic{r}_0 \cong g_{6,54}: \quad (\bar{e}_1,\bar{e}_2,\bar{e}_3,\bar{e}_4,\bar{e}_5,\bar{e}_6)= \left(e_2,-e_3,e_1,e_5,e_4,\frac{1}{3} e_6\right), \quad\mbox{param.: } (\lambda,\gamma) = \left( 1, \frac{2}{3} \right),\\
\gothic{r}_1 \cong g_{6,76}: \quad (\bar{e}_1,\bar{e}_2,\bar{e}_3,\bar{e}_4,\bar{e}_5,\bar{e}_6) = \left(e_3,e_4,-e_2,e_1,e_5,\frac{1}{3} e_6 \right), \quad\mbox{param.: } h = 1.
 \end{gather*}

 Let us be more explicit about the direct verif\/ication of Proposition~\ref{9d-syms} from the point of view of external symmetries, internal symmetries, and symmetries of the lifted coframe on $\Sigma_7 \times H$.

 \subsubsection{External symmetries}

 Given any vector f\/ield $X$ on $J^0(\R^2,\R)$, there is a corresponding prolonged vector f\/ield $X^{(2)}$ on~$J^2(\R^2,\R)$.  This prolongation is uniquely determined by the condition that $\Lieder{X^{(2)}} \contact{2} \subset \contact{2}$, where $\contact{2}$ is the contact system on $J^2(\R^2,\R)$.  See \eqref{prolongation} for the standard prolongation formula.  For the vector f\/ields in Proposition \ref{9d-syms}, we have
 \begin{gather*}
 X_1^{(2)} = X_1, \qquad X_2^{(2)} = X_2, \qquad X_3^{(2)} = X_3, \\
 X_4^{(2)} = X_4 + \parder{p}, \qquad X_5^{(2)} = X_5 + \parder{q}, \qquad
 X_6^{(2)} = X_6 + p\parder{p} + q\parder{q}, \\
 X_7^{(2)} = X_7 + 3p \parder{p} + 2q \parder{q} + 3r \parder{r} + 2s \parder{s} + t \parder{t},\\
 X_8^{(2)} = X_8 - q \parder{p} - \alpha y \parder{q} - 2s \parder{r} - t \parder{s} - \alpha \parder{t},\\
 X_9^{(2)} = X_9 + (z - xp - yq) \parder{p} - \frac{\alpha}{2} y^2 \parder{q} - (3xr+2ys) \parder{r} - (2xs+yt) \parder{s} - (xt+\alpha y) \parder{t}.
 \end{gather*}

 For \eqref{special-maxsym-eqn} or \eqref{gen-eqn-alpha}, we verify the external symmetry condition
 \begin{gather*}
 \Lieder{X_i{}^{(2)}} F =0 \qbox{whenever} F=0.
 \end{gather*}
 Clearly this is satisf\/ied by $X_i^{(2)}$, $i=1,\dots,6$ since they have no components in the $\parder{r}$, $\parder{s}$, $\parder{t}$ direction and since $F=F(r,s,t)$ for \eqref{special-maxsym-eqn} and~\eqref{gen-eqn-alpha}.  For the remaining vector f\/ields we have
 \begin{gather*}
  \eqref{special-maxsym-eqn}: \qquad  \Lieder{X_7^{(2)}} F = 4 F, \qquad
 \Lieder{X_8^{(2)}} F = 0, \qquad
 \Lieder{X_9^{(2)}} F = -4x F,\\
  \eqref{gen-eqn-alpha}: \qquad  \Lieder{X_7^{(2)}} F = 6 F, \qquad
 \Lieder{X_8^{(2)}} F = 0, \qquad
 \Lieder{X_9^{(2)}} F = -6x F,
 \end{gather*}
 and so the external symmetry condition is satisf\/ied.

 \subsubsection{Internal symmetries}

 The symmetry generators $X_i^{(2)}$ are all tangent to the equation manifold $F=0$, so they induce (via the parametrization \eqref{rst-param}) corresponding vector f\/ields $Z_i$ on $\Sigma_7$.  Letting $X_i^{(1)} = (\pi^2_1)_* X_i^{(2)}$ denote the projection onto $J^1(\R^2,\R)$, and identifying the coordinates $(x,y,z,p,q)$ on $J^1(\R^2,\R)$ with corresponding coordinates on $\Sigma_7$, we have
 \begin{gather*}
 Z_i = X_i^{(1)},  \quad i=1,\dots,6, \qquad
 Z_7 = X_7^{(1)} + w\parder{w} + v\parder{v}, \\
 Z_8 = X_8^{(1)} + \parder{v} - m\parder{w}, \qquad
 Z_9 = X_9^{(1)} - (m y + xw) \parder{w} + (y-xv) \parder{v},
 \end{gather*}
 with $u = w + mv$.  One can verify directly that these vector f\/ields satisfy the internal symmetry condition
 \begin{gather*}
 \Lieder{Z_i} I_F \subset I_F,
 \end{gather*}
 where $I_F = \langle \omega^1, \omega^2, \omega^3 \rangle$ is given by the explicit coframing \eqref{9d-explicit-coframe}.

 \subsubsection[Symmetries of the lifted coframe on $\Sigma_7 \times H'$, where $H' = H^0 \rtimes \langle R^2 \rangle$]{Symmetries of the lifted coframe on $\boldsymbol{\Sigma_7 \times H'}$, where $\boldsymbol{H' = H^0 \rtimes \langle R^2 \rangle}$}

 The lifted coframe $\bm{\hat\omega} = \{ \hat\omega^1, \dots, \hat\omega^7, \hat\alpha^1, \hat\alpha^2 \}$ on $\Sigma_7 \times H'$ is parametrized by
 \begin{gather*}
 \hat\omega^1 = a_1{}^2 \omega^1,
 \qquad \hat\omega^2 = a_1 \omega^2 + a_1 a_2 \omega^1,
 \qquad \hat\omega^3 = a_1 \omega^3 + m a_1 a_2 \omega^1,\\
  \hat\omega^4 = a_1 \omega^4,
\qquad \hat\omega^5 = \omega^5 + \epsilon m a_2 \omega^4,
\qquad \hat\omega^6 = a_1 \omega^6, \qquad
 \hat\omega^7 = \omega^7 + a_2 \omega^6,\\
\hat\alpha^1 = \frac{da_1}{a_1} + \frac{a_2}{2a_1} (\hat\omega^4 +  m \hat\omega^6), \qquad
 \hat\alpha^2 = \frac{da_2}{a_1} + \frac{3a_2}{2a_1} \left( \frac{\epsilon}{ m} \hat\omega^5 +  m\hat\omega^7 \right) - \frac{a_2{}^2}{2a_1{}^2} (\hat\omega^4 +  m \hat\omega^6),
\end{gather*}
 and by construction $\bm{\hat\omega}$ is an $\{e\}$-structure on $\Sigma_7 \times H'$, so that a symmetry is by def\/inition a~map $\Phi : \Sigma_7 \times H' \ra \Sigma_7 \times H'$ such that
 \begin{gather*}
 \Phi^* \bm{\hat\omega}^i = \bm{\hat\omega}^i, \qquad  i=1,\dots, 9,
 \end{gather*}
 with inf\/initesimal analogue
 \begin{gather*}
 {\cal L}_{\hat{Z}} \bm{\hat\omega}^i = 0, \qquad i=1,\dots, 9.
 \end{gather*}
 Explicitly, these lifted vector f\/ields are given by
 \begin{gather*}
 \hat{Z}_i = Z_i, \quad i = 1,\dots,5, \qquad
 \hat{Z}_6 = Z_6 - a_1 \parder{a_1} - a_2 \parder{a_2}, \qquad
 \hat{Z}_7 = Z_7 - \frac{3a_1}{2} \parder{a_1} - \frac{3a_2}{2}  \parder{a_2},\! \\
  \hat{Z}_8 = Z_8, \qquad
 \hat{Z}_9 = Z_9 - \frac{x a_1}{2} \parder{a_1} + \left( \frac{1}{(w+m v)^{3/2}} - \frac{x a_2}{2}\right) \parder{a_2}.
 \end{gather*}

\subsection[Amp\`ere contact transformations and $3z_{xx} (z_{yy})^3 + 1=0$]{Amp\`ere contact transformations and $\boldsymbol{3z_{xx} (z_{yy})^3 + 1=0}$}

 After deriving the normal form
 \begin{gather}
  rt - s^2 - \frac{t^4}{12} = 0, \label{special-maxsym-eqn2}
 \end{gather}
 which appeared in \eqref{special-maxsym-eqn}, Vranceanu remarks that if one makes an {\em Amp\`ere contact transformation}, then \eqref{special-maxsym-eqn2} can be reduced to the simpler form
 \begin{gather}
 rt^3 + \frac{1}{12} = 0. \label{special-maxsym-eqn3}
 \end{gather}
 The notion of an Amp\`ere contact transformation is never def\/ined in Vranceanu's paper and does not appear to be common terminology in the literature.  This terminology is, however, referred to brief\/ly in recent work by Stormark (see page~275 in \cite{Stormark2000}).  Namely, Stormark def\/ines it as the genuine (i.e.\ non-point) contact transformation $\Phi$ of $J^1(\R^2,\R)$ given by
 \begin{gather*}
 (\bar{x},\bar{y},\bar{z},\bar{p},\bar{q}) = \left(p,y,z-px, -x, q\right)
 \end{gather*}
 which is clearly contact since
 \begin{gather*}
 d\bar{z} - \bar{p} d\bar{x} - \bar{q} d\bar{y} =  d(z-px) + x dp - q dy = dz - p dx - q dy.
 \end{gather*}
 This is essentially akin to the Legendre transformation from Hamiltonian mechanics, but only acting with respect to the $x$, $z$, $p$ variables.
 For our purposes, we consider the corresponding Legendre-like transformation acting with respect to $y$, $z$, $q$ variables, namely
 \begin{gather*}
 (\bar{x},\bar{y},\bar{z},\bar{p},\bar{q}) = \left(x,q,z-qy, p, -y\right).
 \end{gather*}
 The prolongation of this transformation to $J^2(\R^2,\R)$ satisf\/ies
 \begin{gather*}
 d\bar{p} - \bar{r} d\bar{x} - \bar{s} d\bar{y}
 = dp - \bar{r} dx - \bar{s} dq \equiv rdx + sdy - \bar{r} dx - \bar{s}( sdx + tdy) \qquad \mod \contact{2}\\
\phantom{ d\bar{p} - \bar{r} d\bar{x} - \bar{s} d\bar{y}}{} \equiv (r - s\bar{s} - \bar{r}) dx + (s - t\bar{s}) dy \qquad \mod \contact{2},\\
 d\bar{q} - \bar{s} d\bar{x} - \bar{t} d\bar{y}
 = -dy - \bar{s} dx - \bar{t} dq \equiv -dy - \bar{s} dx - \bar{t} (sdx + tdy) \qquad \mod \contact{2}\\
 \phantom{d\bar{q} - \bar{s} d\bar{x} - \bar{t} d\bar{y}}{} \equiv -(\bar{s} + s\bar{t}) dx - (1 + t\bar{t}) dy \qquad \mod \contact{2},
 \end{gather*}
 and hence
  \begin{gather*}
 (\bar{r},\bar{s},\bar{t}) = \left(\frac{rt-s^2}{t}, \frac{s}{t}, -\frac{1}{t} \right).
 \end{gather*}
 Consequently
 \begin{gather*}
 0 = rt - s^2 - \frac{t^4}{12} = -\frac{\bar{r}}{\bar{t}} - \frac{1}{12 \bar{t}^4} = -\frac{1}{\bar{t}^4} \left(\bar{r}\bar{t}^3 + \frac{1}{12} \right) \qRa \bar{r}\bar{t}^3 + \frac{1}{12} = 0.
 \end{gather*}
 By applying the subsequent scaling $x = \frac{1}{2} \bar{x}$ (and hence $(r,s,t) = (4\bar{r},2\bar{s},\bar{t})$), we are led to the equation
 \begin{gather}
 3rt^3+1=0, \label{rt3}
 \end{gather}
 which was investigated by Goursat \cite{Goursat1898} who recognized its Darboux integrability.

 Since \eqref{rt3} is contact-equivalent to \eqref{special-maxsym-eqn2}, it is clear that \eqref{rt3} is hyperbolic of generic type with $\Delta_1 = \Delta_2 = 0$ and $\epsilon = a = 1$.  The standard Lie algorithm to calculate symmetries can be applied for this equation in a straightforward manner.  Its contact symmetry algebra consists of (projectable) point symmetries $X_1,\dots, X_6$ as in Proposition~\ref{9d-syms} as well as
 \begin{gather}
  X_7 = xy\parder{z}, \qquad
   X_8 = 2y\parder{y} + 3z\parder{z}, \qquad
   X_9 = x^2\parder{x} + xz\parder{z}.
   \label{rt3-sym-alg}
 \end{gather}
 We note that the vector f\/ields $X_7$, $X_8$, $X_9$ have prolongations
 \begin{gather*}
 X_7^{(2)} = X_7 + y\parder{p} + x\parder{q} + \parder{s},\\
  X_8^{(2)} = X_8 + 3p\parder{p} + q\parder{q} + 3r\parder{r} + s\parder{s} - t\parder{t},\\
   X_9^{(2)} = X_9 + (z-xp)\parder{p} + xq\parder{q} - 3xr\parder{r} + (q-xs)\parder{s} + xt\parder{t}.
 \end{gather*}

\subsection{Darboux integrability}

\begin{definition} For a hyperbolic PDE $F=0$, $I_F$ is said to be Darboux-integrable (at level two) if each of $C(I_F,dM_1)$ and $C(I_F,dM_2)$ contains a completely integrable subsystem of rank two that is independent from $I_F$.
\end{definition}

 Recall that for our adapted coframe as in Theorem \ref{generic-hyp-str-eqns}, we have
 \begin{gather*}
 C(I_F,dM_1)^{(2)} = \{ \omega^4, \omega^5 \} \qquad\mbox{and}\qquad C(I_F,dM_2)^{(2)} = \{ \omega^6, \omega^7 \}.
 \end{gather*}

  \begin{theorem} \label{thm:Darboux-int} Given a generic hyperbolic PDE $F=0$ with (maximal) $9$-dimensional symmetry group, the second derived systems $C(I_F,dM_1)^{(2)}$ and $C(I_F,dM_2)^{(2)}$:
 \begin{enumerate}\itemsep=0pt
   \item[1)] are completely integrable, and hence $I_F$ is Darboux integrable, and
   \item[2)] contain rank one completely integrable subsystems.
 \end{enumerate}
 \end{theorem}

 \begin{proof}
 Referring to the maximally symmetric structure equations \eqref{9dim-streqns}, we have that
 \begin{gather*}
 d\omega^4 \equiv d\omega^5 \equiv 0 \mod C(I_F,dM_1)^{(2)}, \qquad
 d\omega^6 \equiv d\omega^7 \equiv 0 \mod C(I_F,dM_2)^{(2)}.
 \end{gather*}
 Hence, the rank two systems $C(I_F,dM_i)^{(2)}$, $i=1,2$ are complete integrable and $I_F = \{ \omega^1, \omega^2$, $ \omega^3 \}$ is Darboux integrable.  Moreover, since
 \begin{gather*}
 d\omega^5 = m \omega^5 \wedge \omega^7, \qquad
 d\omega^7 = -\frac{\epsilon}{m} \omega^5 \wedge \omega^7,
 \end{gather*}
 then the rank one subsystems $\{ \omega^5 \}$ and $\{ \omega^7\}$ are also completely integrable.
 \end{proof}

 Abstractly, Darboux's integration method for these systems proceeds as follows.  Darboux integrability of $I_F$ implies the existence of completely integrable subsystems $J_i \subset C(I_F,dM_i)$.  Applying the Frobenius theorem to each subsystem $J_i$, there exist local functions $f_i$, $g_i$ called {\em Riemann invariants} such that
 \begin{gather*}
  J_1 = \{ df_1, dg_1 \} \subset C(I_F,dM_1), \qquad
 J_2 = \{ df_2, dg_2 \} \subset C(I_F,dM_2).
 \end{gather*}
 If $\varphi_1$, $\varphi_2$ are arbitrary functions, then restricting to any submanifold determined by
 \begin{gather*}
 S: \quad g_1 = \varphi_1 (f_1), \qquad g_2 = \varphi_2 (f_2),
 \end{gather*}
the structure equations \eqref{hyp-str-eqns} become
 \begin{gather*}
 d\tilde{\omega}^i \equiv 0 \quad \mod \tilde{I}_F, \qquad i=1,2,3,
 \end{gather*}
 where $\tilde{I}_F = \{ \tilde{\omega}^1, \tilde{\omega}^2, \tilde{\omega}^3 \}$ is the restriction of $I_F$ to $S$.  Hence, $\tilde{I}_F$ is completely integrable, and so there exist local functions $h_1$, $h_2$, $h_3$ on $S$ such that
 \begin{gather*}
 \tilde{I}_F = \{ dh_1, dh_2, dh_3 \}.
 \end{gather*}
 Hence, these functions $h_1$, $h_2$, $h_3$ are f\/irst integrals of $\tilde{I}_F$, and together with the constraint $S$ determine f\/irst integrals of $I_F$.

 Explicitly, from our parametrization of the coframe $\{ \omega^i \}_{i=1}^7$ on $\Sigma_7$ (c.f.\ \eqref{9d-explicit-coframe}), we have:
 \begin{gather*}
 \omega^4  = u^{3/2} dx + m\sqrt{u} (dy - vdx) = \sqrt{u} (wdx + mdy)
 = \sqrt{u}(d(my+wx) - xdw),\\
 \omega^5  = \frac{\epsilon m}{u} (du - mdv) = \frac{\epsilon m}{u} dw, \\
 \omega^6  = - \sqrt{u} (dy - vdx) = -\sqrt{u} (d(y-vx) + xdv),\\
 \omega^7  = \frac{dv}{u},
 \end{gather*}
 and so
 \begin{gather*}
  C(I_F,dM_1)^{(2)} = \{ \omega^4, \omega^5 \} = \{ dw, d(my+wx) \},\\
 C(I_F,dM_2)^{(2)} = \{ \omega^6, \omega^7 \} = \{ dv, d(y-vx) \},
 \end{gather*}
 where $w = u - mv$.  Thus,
\begin{gather*}
w, \quad my+wx, \qquad\mbox{and}\qquad v, \quad y-vx
\end{gather*}
 are Riemann invariants and, in principle, Darboux's integration method may be applied to f\/ind solutions or f\/irst integrals to the original equation.  In \cite[Corollary~5.9]{GK1993}
  Gardner--Kamran
  asserted that hyperbolic equations of generic type do not have Riemann invariants.  As f\/irst remarked by Eendebak~\cite{Eendebak2006}, this statement is incorrect and the equation $3rt^3+1=0$ is a~counterexample.  Moreover, as described above, all maximally symmetric generic hyperbolic equations have Riemann invariants.

 We refer the reader to page~130 in Goursat \cite{Goursat1898} for the implementation of Darboux's method to the equation $3rt^3+1=0$.
 The implementation of Darboux's method in the case $(\epsilon,a) \neq (1,1)$ appears to be computationally quite dif\/f\/icult.

 Let us comment on Darboux integrability for the submaximally symmetric cases described in Table~\ref{streqn-classification}.  Recall that the structure equations listed in Sections \ref{8d-str-eqs} and \ref{7d-str-eqs} are those for the lifted coframe.  To obtain the structure equations for the corresponding base coframe, we simply set $\hat\alpha^1 = 0$ and remove all hats.  For all these cases we have either that
 \begin{gather*}
 C(I_F,dM_1)^{(3)} = \{ \omega^4, \omega^5 \}, \qquad C(I_F,dM_2)^{(3)} = \{ \omega^6, \omega^7 \}
 \end{gather*}
 and hence $I_F$ is Darboux integrable, or
 \begin{gather*}
 C(I_F,dM_1)^{(3)} = \{ \omega^5 \}, \qquad C(I_F,dM_2)^{(3)} = \{ \omega^7 \}
 \end{gather*}
 and $I_F$ is not Darboux integrable.  We list the possibilities in Table~\ref{table:Darboux-int}.
 Moreover, for these submaximally symmetric cases, all which are Darboux integrable have one-dimensional subsystems of $C(I_F,dM_1)^{(2)}$ and $C(I_F,dM_2)^{(2)}$ which are completely integrable (namely, $\{ \omega^5 \}$ and $\{ \omega^7 \}$ respectively).  Thus, the converse of Theorem \ref{thm:Darboux-int} is clearly {\em false}.

 \begin{table}[h]
 \centering
  \caption{Darboux integrability of submaximally symmetric generic hyperbolic PDE.}
 \label{table:Darboux-int}

 \vspace{1mm}

 \begin{tabular}{|c|c|c|}\hline
 Case & Darboux integrable? \\ \hline\hline
 2a & no\\
 2b & no in general; yes if $(m,\epsilon_1) = (1,-1)$\\\hline
 3a & yes\\
 3b & no in general; yes if $(m,\epsilon_1) = (1,-1)$\\ \hline
 \end{tabular}
 \end{table}

 \section{Concluding remarks}
 \label{genhyp:conclusions}

 Let us summarize some of the main results of this paper:
 \begin{itemize}\itemsep=0pt
 \item We derived relative invariants $I_1$, $I_2$ (see Theorem \ref{thm:hyp-contact-inv}) given parametrically in terms of an arbitrary hyperbolic equation $F(x,y,z,z_x,z_y,z_{xx},z_{xy},z_{yy}) = 0$.  Their vanishing/non\-va\-ni\-shing distinguishes the three types of hyperbolic equations.
 \item In the generic case, the $\epsilon$ contact invariant is given parametrically as $\epsilon = {\rm sgn}(I_1 I_2) = \pm 1$.
 \item In the abstract analysis of the generic hyperbolic structure equations, we identif\/ied relative contact invariants $m$, $n$, $B$ and $\Delta_1 = mn + \epsilon$, $\Delta_2 = m^2 - \epsilon n^2$ which played a key role in the classif\/ication of various generic hyperbolic structures admitting nine, eight, and seven-dimensional symmetry along with the corresponding complete structure equations.
 \item Integration of maximally symmetric structure equations, leading to normal forms for all contact-equivalence classes of maximally symmetric generic hyperbolic equations.
 \item Nine-dimensional symmetry algebras for these normal forms for generic hyperbolic equations are given explicitly.  There are exactly two such nonisomorphic algebras.
 \item For any maximally symmetric generic hyperbolic equation, the second derived systems of $C(I_F,dM_i)$, $i=1,2$ are rank 2 and completely integrable.  Hence, all maximally symmetric generic hyperbolic equations are Darboux integrable.
 \end{itemize}

 We conclude with some possible points for future investigation:
 \begin{enumerate}\itemsep=0pt
 \item Maximally symmetric equations: (1)~Do ``simpler'' normal forms exist?  (2)~Implement Darboux's integration method in the general case. (3)~Investigate the existence of conservation laws. (4)~Study the local solvability of these equations.
 \item Submaximally symmetric equations: Integrate the structure equations given in Sections~\ref{8d-str-eqs} and~\ref{7d-str-eqs} and f\/ind normal forms for the corresponding PDE equivalence classes.  Address similar questions as above.
 \item The submaximally symmetric structures that we have derived here (see Table~\ref{streqn-classification} and Sections~\ref{8d-str-eqs} and~\ref{7d-str-eqs}) share the common property that $m$, $n$ are constants and $K^0$ is a subgroup of the structure group.  Are there any other reductions of the initial 3-dimensional structure group that lead to valid structures?
 \item In this article, we have carried out a detailed analysis of the generic (7-7) case.  Hyperbolic equations of Goursat (6-7) type are equally poorly understood.
 Some preliminary results on structure equations were stated in \cite{GK1993}, but to our knowledge, Vranceanu's student Petrescu \cite{Petrescu1938} has written the only paper which has made a more detailed investigation into the contact geometry of the Goursat class.  Recasting Petrescu's results for a contemporary audience and building upon his work would make for a natural sequel to our paper.
 \end{enumerate}

\appendix

\section{Hyperbolic structure equations}
\label{app-A}

 We give here the details of the proof of Theorem \ref{general-hyp-str-eqns} starting from the preliminary hyperbolic structure equations \eqref{hyp-str-eqns}.  The main details of this proof have appeared in \cite{Vranceanu1937} and \cite{GK1993}.


 \begin{note} In this section, we will def\/ine changes in the coframe basis using a ``bar'', e.g. $\bar\omega^i = g^i{}_j \omega^j$, but we make the convention that the bar is immediately dropped afterwards, i.e.\ so that $\omega^i$ is redef\/ined to be $\bar\omega^i$.
 \end{note}

 \begin{proof}[Proof of Theorem \ref{general-hyp-str-eqns}]
 From the preliminary hyperbolic structure equations \eqref{hyp-str-eqns}, we have
 \begin{gather}
   d\omega^1 \equiv \omega^3 \wedge \left( a_4 \omega^4 + a_5 \omega^5 + a_6 \omega^6 + a_7 \omega^7 \right)
 \quad \mod M_1, \label{domega1}
 \end{gather}
 but since
 \begin{gather*}
   0  = d^2\omega^1 \wedge \omega^1 \wedge \omega^2 \wedge \omega^3 \wedge \omega^5 = a_4 \omega^4 \wedge \omega^6 \wedge \omega^7 \wedge \omega^1 \wedge \omega^2 \wedge \omega^3 \wedge \omega^5, \\
   0  = d^2\omega^1 \wedge \omega^1 \wedge \omega^2 \wedge \omega^3 \wedge \omega^4 = a_5 \omega^5 \wedge \omega^6 \wedge \omega^7 \wedge \omega^1 \wedge \omega^2 \wedge \omega^3 \wedge \omega^4,
 \end{gather*}
 then $a_4=a_5=0$.   Assuming $a_6=a_7=0$, then $d\omega^1 \equiv 0$ $\mod M_1$, which would imply that $d\omega^1 \wedge d\omega^1 \wedge \omega^1 = 0$, i.e.\ the Darboux rank $\rho$ of $I_F^{(1)} = \{ \omega^1 \}$ would be $\rho \leq 1$.  However, since $\omega^1 = i_F^*(dz - pdx - qdy) \in I_F^{(1)}$ and $d\omega^1 = i_F^*(dx \wedge dp + dy \wedge dq)$, then $d\omega^1 \wedge d\omega^1 \wedge \omega^1 = i_F^*(dx \wedge dp \wedge dy \wedge dq \wedge dz) \neq 0$ since $i_F$ is a nondegenerate parametrization.  Hence $\rho \geq 2$ and we get a contradiction.  Thus, $a_6 \omega^6 + a_7 \omega^7 \neq 0$.

 Choose $A$, $B$ such that $Ba_6 - A a_7=1$, and redef\/ine
 \begin{gather*}
 \bar\omega^6 = a_6 \omega^6 + a_7 \omega^7, \qquad \bar\omega^7 = A \omega^6 + B \omega^7.
 \end{gather*}
 This preserves both the $d\omega^2$ and $d\omega^3$ structure equations (since $\bar\omega^6 \wedge \bar\omega^7 = \omega^6 \wedge \omega^7$).  We have
 \begin{gather*}
   d\omega^1 \equiv \omega^3 \wedge \omega^6,
 \qquad d\omega^2 \equiv \omega^4 \wedge \omega^5 + \omega^3 \wedge  \left( \sum_{k=4}^7 b_k \omega^k \right)
 \quad \mod M_1.
 \end{gather*}
 Redef\/ining $(\bar\omega^2,\bar\omega^4,\bar\omega^5) = (\omega^2 - b_6 \omega^1, \omega^4 + b_5 \omega^3, \omega^5 - b_4 \omega^3)$, we may without loss of generality take $b_4=b_5=b_6=0$.  (Note that the $d\omega^3$ equation is not af\/fected.)  Setting $b_7=U_1$, we have
 \begin{gather}
    d\omega^1 \equiv \omega^3 \wedge \omega^6,
  \qquad d\omega^2 \equiv \omega^4 \wedge \omega^5 + U_1\omega^3 \wedge \omega^7
  \quad \mod M_1. \label{domega12}
 \end{gather}
 A similar argument modulo $\{ \omega^1, \omega^3 \}$ and for the $M_2$ characteristic system yields
 \begin{gather}
    d\omega^1 \equiv \omega^2 \wedge \omega^4,
  \qquad d\omega^3 \equiv \omega^6 \wedge \omega^7 + U_2\omega^2 \wedge \omega^5
  \quad \mod M_2, \label{domega13}
 \end{gather}
 and we emphasize that $\omega^i$ in \eqref{domega12} are the same as $\omega^i$ in~\eqref{domega13}.
 This implies
 \begin{gather*}
 d\omega^1 \equiv \omega^2 \wedge \omega^4 + \omega^3 \wedge \omega^6 + W \omega^2 \wedge \omega^3 \quad \mod \{ \omega^1 \}.
 \end{gather*}
 We calculate
 \begin{gather*}
 0  = d^2 \omega^2 \wedge \omega^1 \wedge \omega^2 \wedge \omega^3 \wedge \omega^4 = d\omega^4 \wedge \omega^5 \wedge \omega^1 \wedge \omega^2 \wedge \omega^3 \wedge \omega^4, \\
 0  = d^2 \omega^1 \wedge \omega^1 \wedge \omega^3 \wedge \omega^4 \wedge \omega^5 = -\omega^2 \wedge (d\omega^4 + W \omega^6 \wedge \omega^7 )\wedge \omega^1 \wedge \omega^3 \wedge \omega^4 \wedge \omega^5.
 \end{gather*}
 The f\/irst equation implies that $\gamma^4{}_{67}=0$ (c.f.\ as def\/ined in~\eqref{gamma-defn})
 and hence the second equation implies $W=0$.
 \end{proof}

 \section{Principal contact invariants for hyperbolic equations}
 \label{app:hyp-contact-inv}

 Here we derive the contact invariants stated in Theorem \ref{thm:hyp-contact-inv}.  We give a brief outline of the computations to follow.
 Ultimately, we are looking for a coframe satisfying the hyperbolic structure equations given in \eqref{U1U2-streq}.  Beginning with the pullback of the basis $\theta^1$, $\theta^2$, $\theta^3$ of the contact ideal on $J^2(\R^2,\R)$, we f\/ind a canonical basis $\tilde\theta^1$, $\tilde\theta^2$, $\tilde\theta^3$ whose pullback brings $\langle \cdot, \cdot \rangle_7$ into the Witt normal form.  This yields the basis $\tilde\omega^1, \dots, \tilde\omega^7$ in Lemma~\ref{lem:tilde-coframe}.  A subsequent normalization in \eqref{semi-normalized-basis} def\/ines the basis $\omega^1, \dots, \omega^7$ which satisf\/ies \eqref{U1U2-streq} and from which (multiples of)~$U_1$ and~$U_2$ can be extracted in parametrized form.  The f\/inal calculation of these parametrized functions is performed in the ambient space $J^2(\R^2,\R)$, essentially using the simple fact given in~\eqref{dF-lemma}.  Pulling these back by $i_F$ yields the desired invariants. We begin with the following lemma.

 \begin{lemma}
  Given a symmetric bilinear form represented by
 \begin{gather}
   Q = \mat{cc}{ a & -b \\ -b & c} \qquad \mbox{with} \quad \Delta = \det(Q) = ac - b^2 < 0,
 \end{gather}
 the change of basis matrix
 \begin{gather*}
      P = \mat{cc}{\lambda_+ & c\\ a & \lambda_+}, \qquad \lambda_+ = b + {\rm sgn}(b) \sqrt{|\Delta|}, \qquad {\rm sgn}(x) = \left\{ \begin{array}{rc} 1, & x \geq 0,\\ -1, & x < 0, \end{array} \right.
 \end{gather*}
 brings $Q$, up to a nonzero scaling, into the Witt normal form, i.e.
 \begin{gather}
   P^t Q P = \rho \mat{cc}{0 & 1 \\ 1 & 0}, \qquad \rho \neq 0.
 \end{gather}
 \end{lemma}

 \begin{proof} Let $v_1$, $v_2$ be the columns of $P$.  Note that $\lambda_\pm = b \pm {\rm sgn}(b) \sqrt{|\Delta|}$ are roots of the polynomial $\lambda^2 - 2b\lambda + ac = 0$ and that $\lambda_+ \neq 0$.  Hence, $v_1$, $v_2$ are nonzero and it is a straightforward verif\/ication that $v_1$ and $v_2$ are null vectors (i.e.\ $v_1^t Q v_1 = v_2^t Q v_2 = 0$).  They are linearly independent since $\det(P) = 2\sqrt{|\Delta|}( |b| + \sqrt{|\Delta|} ) \neq 0$.  Finally, they have nonzero scalar product $\rho = v_1^t Q v_2 = 2\Delta \lambda_+ \neq 0$.
 \end{proof}

 We apply this elementary result to the $2 \times 2$ submatrix appearing in the matrix in \eqref{bilinear-form-matrix} representing the bilinear form $\langle \cdot, \cdot \rangle_7$,
 \begin{gather*}
 \mat{cc}{F_t & -\frac{1}{2} F_s\\ -\frac{1}{2} F_s & F_r}
 \end{gather*}
 with respect to the basis $\theta^2$, $\theta^3$, and arrive at the canonical basis
 \begin{gather*}
     \tilde\theta^2 = \lambda_+ \theta^2 + F_t \theta^3, \qquad
   \tilde\theta^3 = F_r \theta^2 + \lambda_+ \theta^3.
 \end{gather*}

 Without loss of generality, we assume that $F_s \geq 0$.  (If not, consider the equation $-F=0$ which def\/ines the same locus as $F=0$.)
 In the notation of the previous lemma, we have
 \begin{gather*}
 \Delta = F_t F_r - \frac{1}{4} F_s^2 < 0, \qquad \lambda_\pm = \frac{F_s}{2} \pm \sqrt{|\Delta|}.
 \end{gather*}
 Note that $\lambda_+ > 0$, but $\lambda_-$ may vanish.  Both $\lambda_\pm$ are roots of the polynomial
 \begin{gather*}
 \lambda^2 - F_s \lambda + F_t F_r = 0.
 \end{gather*}
 Furthermore, let us note the identities
 \begin{gather*}
  \lambda_+ \lambda_- = F_t F_r, \qquad \lambda_+ + \lambda_- = F_s, \qquad
 \lambda_+^2 - F_t F_r = 2\lambda_+\sqrt{|\Delta|}.
 \end{gather*}
 We also def\/ine the total derivative operators
 \begin{gather*}
 D_x = \parder{x} + p\parder{z} + r\parder{p} + s\parder{q}, \qquad
 D_y = \parder{y} + q\parder{z} + s\parder{p} + t\parder{q}.
 \end{gather*}

 \begin{lemma} \label{lem:tilde-coframe} Suppose that $F_s \geq 0$.  Let $\tilde\omega^i = i_F^* \tilde\theta^i$, where
 \begin{gather*}
    \tilde\theta^1 = dz - p dx - q dy, \qquad
   \tilde\theta^2 = \lambda_+ \theta^2 + F_t \theta^3, \qquad
   \tilde\theta^3 = F_r \theta^2 + \lambda_+ \theta^3,\\
    \tilde\theta^4 = dx - \frac{F_r}{\lambda_+} dy, \qquad
   \tilde\theta^5 = \lambda_+ dr + F_t ds  + D_x F dy , \\
    \tilde\theta^6 = dy - \frac{F_t}{\lambda_+} dx, \qquad
   \tilde\theta^7 = \lambda_+ dt + F_r ds  + D_y F dx.
 \end{gather*}
 Then $\bm{\tilde\omega} = \{ \tilde\omega^i \}_{i=1}^7$ is a (local) coframe on~$\Sigma_7$ satisfying the structure equations
 \begin{alignat*}{3}
& d\tilde\omega^1 \equiv \rho (\tilde\omega^2 \wedge \tilde\omega^4 + \tilde\omega^3 \wedge \tilde\omega^6) + \tilde{W}\tilde\omega^2 \wedge \tilde\omega^3 && \mod I_F^{(1)}, & \\
& d\tilde\omega^2 \equiv \tilde\omega^4 \wedge \tilde\omega^5 + \tilde\omega^3 \wedge \rho\, i_F^* \Xi_1 && \mod M_1, &\\
& d\tilde\omega^3 \equiv \tilde\omega^6 \wedge \tilde\omega^7 + \tilde\omega^2 \wedge \rho\, i_F^* \Xi_2 && \mod M_2,
 \end{alignat*}
 where
 \begin{gather*}
 I_F = \{ \tilde\omega^1, \tilde\omega^2, \tilde\omega^3 \}, \qquad I_F^{(1)} = \{ \tilde\omega^1 \}, \qquad \rho = i_F^*\left(\frac{-1}{2\sqrt{|\Delta|}}\right), \qquad \mbox{and}\\
 \Xi_1 = dF_t -\frac{F_t}{\lambda_+} d\lambda_+  - \left(F_q - \frac{F_t F_p}{\lambda_+} \right) dy,\qquad
 \Xi_2 = dF_r - \frac{F_r}{\lambda_+} d\lambda_+ - \left(F_p - \frac{F_r F_q}{\lambda_+} \right) dx.
 \end{gather*}
 \end{lemma}

 \begin{proof} We need to show that $\bm{\tilde\omega}$ is linearly independent, or equivalently $\tilde\omega^1 \wedge \cdots \wedge \tilde\omega^7 \neq 0$.  By~\eqref{dF-lemma}, it suf\/f\/ices to show that $\tilde\theta^1 \wedge \cdots \wedge \tilde\theta^7 \wedge dF \neq 0$.  We calculate
 \begin{gather}
 \tilde\theta^1 \wedge \tilde\theta^2 \wedge \tilde\theta^3 \wedge \tilde\theta^4 \wedge \tilde\theta^6
  = 4|\Delta| dz \wedge  dp \wedge dq \wedge dx \wedge dy \label{wedge12346}
 \end{gather}
 and so
 \begin{gather*}
   \tilde\theta^1 \wedge \cdots \wedge \tilde\theta^7 \wedge dF \\
    =4|\Delta| dz \wedge  dp \wedge dq \wedge dx \wedge dy \wedge (\lambda_+ dt + F_r ds) \wedge ( \lambda_+ dr + F_t ds ) \wedge (F_r dr + F_s ds + F_t dt)\\
   = 4|\Delta| \det\left( \begin{array}{ccc} 0 & F_r & \lambda_+ \\ \lambda_+ & F_t & 0 \\ F_r & F_s & F_t \end{array} \right) \overbrace{dz \wedge  dp \wedge dq \wedge dx \wedge dy \wedge dr \wedge ds \wedge dt}^{{\rm Vol}_{J^2(\R^2,\R)}}\\
   = 4|\Delta| \lambda_+ (-2F_t F_r + \lambda_+ F_s) {\rm Vol}_{J^2(\R^2,\R)}
   = 4|\Delta| \lambda_+^2 (-2\lambda_- + F_s) {\rm Vol}_{J^2(\R^2,\R)}\\
    = 8|\Delta|^{3/2} \lambda_+^2 {\rm Vol}_{J^2(\R^2,\R)} \neq 0.
 \end{gather*}
 Thus, $\bm{\tilde\omega}$ is a coframe on~$\Sigma_7$.
 Let us verify the structure equations.  Note that
 \begin{gather*}
    i_F^*\theta^2 = \frac{1}{\lambda_+^2 - F_t F_r} \left( \lambda_+ \tilde\omega^2 - F_t \tilde\omega^3 \right), \qquad
   i_F^*\theta^3 = \frac{1}{\lambda_+^2 - F_t F_r} \left( -F_r \tilde\omega^2 + \lambda_+ \tilde\omega^3 \right),\\
    i_F^*dx = \frac{\lambda_+^2}{\lambda_+^2 - F_t F_r} \left(\tilde\omega^4 + \frac{F_r}{\lambda_+} \tilde\omega^6\right), \qquad
   i_F^*dy = \frac{\lambda_+^2}{\lambda_+^2 - F_t F_r} \left(\frac{F_t}{\lambda_+} \tilde\omega^4 + \tilde\omega^6\right),
 \end{gather*}
 where all coef\/f\/icients on the right side are pulled back by $i_F^*$ to $\Sigma_7$ (i.e.\ evaluated on $F=0$)
 \begin{gather*}
    d\tilde\omega^1  = i_F^*(dx \wedge dp + dy \wedge dq) = i_F^*(dx \wedge (\theta^2 + s dy) + dy \wedge (\theta^3 + s dx)) \\
  \phantom{d\tilde\omega^1}{} = i_F^*(dx \wedge \theta^2 + dy \wedge \theta^3).
 \end{gather*}

 Mod $M_1$, we have
 \begin{gather*}
   d\tilde\omega^1  \equiv \frac{\lambda_+^2}{(\lambda_+^2 - F_t F_r)^2} \left( \left(\tilde\omega^4 + \frac{F_r}{\lambda_+} \tilde\omega^6\right) \wedge \left( - F_t \tilde\omega^3 \right) + \left(\frac{F_t}{\lambda_+} \tilde\omega^4 + \tilde\omega^6\right) \wedge \left( \lambda_+ \tilde\omega^3 \right) \right)\\
   \phantom{d\tilde\omega^1}{} \equiv \frac{\lambda_+^2}{(\lambda_+^2 - F_t F_r)^2} \left( \lambda_+ - \frac{F_t F_r}{\lambda_+} \right) \tilde\omega^6 \wedge \tilde\omega^3
   \equiv \frac{\lambda_+}{\lambda_+^2 - F_t F_r} \tilde\omega^6 \wedge \tilde\omega^3
   \equiv \frac{-1}{2\sqrt{|\Delta|}} \tilde\omega^3 \wedge \tilde\omega^6.
 \end{gather*}

 Mod $M_2$, we have
 \begin{gather*}
   d\tilde\omega^1  \equiv \frac{\lambda_+^2}{(\lambda_+^2 - F_t F_r)^2} \left( \left(\tilde\omega^4 + \frac{F_r}{\lambda_+} \tilde\omega^6\right) \wedge \left( \lambda_+ \tilde\omega^2 \right) + \left(\frac{F_t}{\lambda_+} \tilde\omega^4 + \tilde\omega^6\right) \wedge \left( -F_r \tilde\omega^2 \right) \right)\\
   \phantom{d\tilde\omega^1}{} \equiv \frac{\lambda_+^2}{(\lambda_+^2 - F_t F_r)^2} \left( \lambda_+ - \frac{F_t F_r}{\lambda_+} \right) \tilde\omega^4 \wedge \tilde\omega^2
   \equiv \frac{\lambda_+}{\lambda_+^2 - F_t F_r} \tilde\omega^4 \wedge \tilde\omega^2
   \equiv \frac{-1}{2\sqrt{|\Delta|}} \tilde\omega^2 \wedge \tilde\omega^4.
 \end{gather*}

 Now examine $d\tilde\omega^2$ and $d\tilde\omega^3$.  Note the relation
  \begin{gather}
    0  = i_F^*(dF) = i_F^*(F_x dx + F_y dy + F_z dz + F_p dp + F_q dq + F_r dr + F_s ds + F_t dt) \nonumber\\
\phantom{0}{}= i_F^*( F_x dx + F_y dy + F_z (\theta^1 + p dx + q dy) + F_p \left(\theta^2 + r dx + s dy \right) + F_q (\theta^3 + s dx + t dy) \nonumber\\
\phantom{0=}{} + F_r dr + F_s ds + F_t dt) \nonumber\\
 \phantom{0}{}= i_F^*( D_x F dx + D_y F dy + F_r dr + F_s ds + F_t dt + F_z  \theta^1 + F_p \theta^2 + F_q \theta^3) \label{dF-rel}
  \end{gather}
 and so,
 \begin{gather*}
   0  \equiv i_F^*\left(D_x F dx + D_y F dy + F_r dr + F_s ds + F_t dt + \left(F_q - \frac{F_t F_p}{\lambda_+} \right) \theta^3\right)\quad \mod M_1,\\
   0  \equiv i_F^*\left(D_x F dx + D_y F dy + F_r dr + F_s ds + F_t dt + \left(F_p - \frac{F_rF_q}{\lambda_+} \right) \theta^2 \right)\quad \mod M_2.
 \end{gather*}

 Mod $M_1$, we have
 \begin{gather*}
    d\tilde\omega^2 \equiv i_F^*( d\lambda_+ \wedge \theta^2 + \lambda_+ d\theta^2 + F_t d\theta^3 + dF_t \wedge \theta^3)\\
\phantom{d\tilde\omega^2}{} \equiv i_F^*\left( -\frac{F_t}{\lambda_+} d\lambda_+ \wedge \theta^3 + \lambda_+ ( dx \wedge dr + dy \wedge ds ) + F_t ( dx \wedge ds + dy \wedge dt ) + dF_t \wedge \theta^3 \right)\\
\phantom{d\tilde\omega^2}{} \equiv i_F^*\left( \lambda_+ dx \wedge dr - \lambda_- dy \wedge ds + F_t dx \wedge ds + D_x F dx \wedge dy - F_r dy \wedge  dr + \Xi_1 \wedge \theta^3 \right)\\
\phantom{d\tilde\omega^2}{} \equiv i_F^* \left( \left( dx - \frac{F_r}{\lambda_+} dy \right) \wedge \left( \lambda_+ dr + F_t ds  + D_x F dy\right) + \Xi_1 \wedge \theta^3 \right)\\
\phantom{d\tilde\omega^2}{}\equiv \tilde\omega^4 \wedge \tilde\omega^5 + \tilde\omega^3 \wedge \rho\, i_F^* \Xi_1.
 \end{gather*}

 Mod $M_2$, we have
 \begin{gather*}
  d\tilde\omega^3  \equiv i_F^*(dF_r \wedge \theta^2 + F_r d\theta^2 + \lambda_+ d\theta^3 + d\lambda_+ \wedge \theta^3) \\
\phantom{d\tilde\omega^3}{} \equiv i_F^*\left(dF_r \wedge \theta^2 + F_r (dx \wedge dr + dy \wedge ds) + \lambda_+ (dx \wedge ds + dy \wedge dt) - \frac{F_r}{\lambda_+} d\lambda_+ \wedge \theta^2\right)\\
\phantom{d\tilde\omega^3}{}\equiv i_F^*\left( F_r dy \wedge ds + \lambda_+ (dx \wedge ds + dy \wedge dt) - dx \wedge \left(D_y F dy + F_s ds + F_t dt \right) + \Xi_2 \wedge \theta^2 \right)\\
\phantom{d\tilde\omega^3}{} \equiv i_F^*\left( \left(dy - \frac{F_t}{\lambda_+} dx\right) \wedge (\lambda_+ dt + F_r ds  + D_y F dx) + \Xi_2 \wedge \theta^2 \right)\\
\phantom{d\tilde\omega^3}{}\equiv \tilde\omega^6 \wedge \tilde\omega^7 + \tilde\omega^2 \wedge \rho\, i_F^* \Xi_2.
\tag*{\qed}
 \end{gather*}\renewcommand{\qed}{}
  \end{proof}

 Let us write
 \begin{gather*}
 \rho\, i_F^* \Xi_1 = b_4 \tilde\omega^4 + b_5 \tilde\omega^5 + b_6 \tilde\omega^6 + b_7 \tilde\omega^7, \qquad
 \rho\, i_F^* \Xi_2 = c_4 \tilde\omega^4 + c_5 \tilde\omega^5 + c_6 \tilde\omega^6 + c_7 \tilde\omega^7.
 \end{gather*}
 We make the change of basis
 \begin{gather}
 \displaystyle \omega^1 = \tilde\omega^1, \qquad \omega^2 = \tilde\omega^2 - \frac{b_6}{\rho} \tilde\omega^1, \qquad \omega^3 = \tilde\omega^3 - \frac{c_4}{\rho} \tilde\omega^1, \qquad \omega^4 = \rho(\tilde\omega^4 + b_5 \tilde\omega^3), \nonumber\\
 \displaystyle \omega^5 = \frac{1}{\rho} (\tilde\omega^5 - b_4 \tilde\omega^3), \qquad
 \omega^6 = \rho(\tilde\omega^6 + c_7 \tilde\omega^2), \qquad \omega^7 = \frac{1}{\rho} (\tilde\omega^7 - c_6 \tilde\omega^2)
 \label{semi-normalized-basis}
 \end{gather}
 to obtain the structure equations
 \begin{alignat*}{3}
& d\omega^1 \equiv \omega^2 \wedge \omega^4 + \omega^3 \wedge \omega^6 + W\omega^2 \wedge \omega^3 && \mod I_F^{(1)},&\\
& d\omega^2 \equiv \omega^4 \wedge \omega^5 + b_7 \omega^3 \wedge \omega^7 && \mod M_1, &\\
& d\omega^3 \equiv \omega^6 \wedge \omega^7 + c_5 \omega^2 \wedge \omega^5 && \mod M_2, &
 \end{alignat*}
 where $W = \tilde{W} + \rho (b_5 - c_7)$.  Thus, in the notation of Theorem~\ref{general-hyp-str-eqns}, we have that $U_1 = b_7$ and $U_2 = c_5$ and by the argument given at the end of Appendix~\ref{app-A}, we must have that $W=0$.

 We are primarily concerned with the vanishing / nonvanishing of $U_1$ and $U_2$, though the sign will play a later role in the determination of the $\epsilon$ contact invariant.  It suf\/f\/ices to calculate
  \begin{gather*}
   U_1: \ \tilde\omega^7 \mbox{ coef\/f\/icient in } i_F^*\Xi_1, \qquad \mbox{and}\qquad
  U_2: \ \tilde\omega^5 \mbox{ coef\/f\/icient in } i_F^*\Xi_2.
  \end{gather*}
  Note that the sign of these coef\/f\/icients is {\em reversed for both} since $\rho < 0$.
 As a further simplif\/ication, since $i_F^*(dx)$, $i_F^*(dy)$ depend only on $\tilde\omega^4$, $\tilde\omega^6$, and since $i_F^*(\lambda_+) > 0$, it suf\/f\/ices to calculate
 \begin{gather*}
   U_1 : \ \tilde\omega^7 \mbox{ coef\/f\/icient in } i_F^*d\left( \frac{F_t}{\lambda_+} \right), \quad \mbox{and}
 \qquad  U_2 : \ \tilde\omega^5 \mbox{ coef\/f\/icient in } i_F^*d\left( \frac{F_r}{\lambda_+} \right).
 \end{gather*}

 Let us illustrate how to calculate the $\tilde\omega^7$ coef\/f\/icient in $i_F^*d\left( \frac{F_t}{\lambda_+} \right)$.
 Given $\tilde\omega = f_i \tilde\omega^i = i_F^*d\left( \frac{F_t}{\lambda_+} \right)$ on $\Sigma_7$, we are interested in $f_7$.  Since $i_F : \Sigma_7 \ra J^2(\R^2,\R)$ is maximal rank, there exist functions~$\tilde{f}_i$ on $J^2(\R^2,\R)$ such that $i_F^* \tilde{f}_i = f_i$, so $\tilde\omega = i_F^*( \tilde{f}_i \tilde\theta^i )$.  Thus, we calculate $\tilde{f}_7$ using
 \begin{gather*}
 \tilde\theta^1 \wedge \cdots \wedge \tilde\theta^6 \wedge d\left( \frac{F_t}{\lambda_+} \right) \wedge dF = \tilde{f}_7 \tilde\theta^1 \wedge \cdots \wedge \tilde\theta^7 \wedge dF
 \end{gather*}
 and then use $i_F^*$ to obtain $f_7$.

 Using \eqref{wedge12346}, we have
 \begin{gather}
 \nu_1  := \tilde\theta^1 \wedge \tilde\theta^2 \wedge \tilde\theta^3 \wedge \tilde\theta^4 \wedge \tilde\theta^5 \wedge \tilde\theta^6 \wedge d\left( \frac{F_t}{\lambda_+} \right) \wedge dF \nonumber\\
\phantom{\nu_1}{} =-4 |\Delta| dz \wedge dp \wedge dq \wedge dx \wedge dy \wedge dF \wedge (\lambda_+ dr + F_t ds) \wedge d\left( \frac{F_t}{\lambda_+} \right) \nonumber\\
 \phantom{\nu_1}{}=-4|\Delta| \tilde{I}_1 dz \wedge dp \wedge dq \wedge dx \wedge dy \wedge dr \wedge ds \wedge dt \nonumber\\
\phantom{\nu_1}{}= \tilde{f}_7 \tilde\theta^1 \wedge \cdots \wedge \tilde\theta^7 \wedge dF, \label{vol1}
 \end{gather}
 where
 \begin{gather*}
 \tilde{I}_1 = \det\left( \begin{array}{ccc} F_r & F_s & F_t\\ \lambda_+ & F_t & 0\\
  \left( \frac{F_t}{\lambda_+} \right)_r &
  \left( \frac{F_t}{\lambda_+} \right)_s &
  \left( \frac{F_t}{\lambda_+} \right)_t
  \end{array} \right), \qquad \tilde{f}_7 = \frac{-\tilde{I}_1}{2|\Delta|^{1/2}\lambda_+^2}.
 \end{gather*}
 Let $I_1 = i_F^*\tilde{I}_1$, and note $I_1$ and $f_7 = i_F^* \tilde{f}_7$ have opposite sign. Because of the aforementioned sign reversal (i.e.\ $U_1$ and $f_7$ have opposite sign), we have
 \begin{gather}
 {\rm sgn}(U_1) = {\rm sgn}(I_1). \label{sgnU1}
 \end{gather}

 We perform a similar computation for the $\omega^5$ coef\/f\/icient in $i_F^*d\left( \frac{F_r}{\lambda_+} \right)$.
 \begin{gather*}
 \nu_2  := \tilde\theta^1 \wedge \tilde\theta^2 \wedge \tilde\theta^3 \wedge \tilde\theta^4 \wedge  d\left( \frac{F_r}{\lambda_+} \right) \wedge \tilde\theta^6 \wedge \tilde\theta^7 \wedge dF \\
\phantom{\nu_2}{} =-4 |\Delta| dz \wedge dp \wedge dq \wedge dx \wedge dy \wedge (\lambda_+ dt + F_r ds) \wedge dF \wedge d\left( \frac{F_r}{\lambda_+} \right)\\
\phantom{\nu_2}{}=-4|\Delta| \tilde{I}_2 dz \wedge dp \wedge dq \wedge dx \wedge dy \wedge dr \wedge ds \wedge dt\\
\phantom{\nu_2}{}= \tilde{f}_5 \tilde\theta^1 \wedge \cdots \wedge \tilde\theta^7 \wedge dF, \nonumber
 \end{gather*}
 where
 \begin{gather*}
 \tilde{I}_2 = \det\left( \begin{array}{ccc} 0 & F_r & \lambda_+ \\ F_r & F_s & F_t\\
  \left( \frac{F_r}{\lambda_+} \right)_r &
  \left( \frac{F_r}{\lambda_+} \right)_s &
  \left( \frac{F_r}{\lambda_+} \right)_t
  \end{array} \right), \qquad \tilde{f}_5 = \frac{-\tilde{I}_2}{2|\Delta|^{1/2}\lambda_+^2}.
 \end{gather*}
 Let $I_2 = i_F^*\tilde{I}_2$, and note $I_2$ and $f_5 = i_F^* \tilde{f}_5$ have opposite sign. Because of the aforementioned sign reversal (i.e.\ $U_2$ and $f_5$ have opposite sign), we have
 \begin{gather}
 {\rm sgn}(U_2) = {\rm sgn}(I_2). \label{sgnU2}
 \end{gather}
 Since class 6-6, 6-7, and 7-7 hyperbolic equations are determined by the vanishing/nonvanishing of $U_1$, $U_2$, the classif\/ication in Theorem~\ref{thm:hyp-contact-inv} follows.

 For the scaling property, suppose that $\hat{F} = \phi F$ with $i_F^*\phi > 0$.  We may without loss of generality suppose that $\hat\Sigma_7 = \Sigma_7$ and $i_{\hat{F}} = i_F$.  Hence,
 \begin{gather*}
 i_{\hat{F}}^*\hat{F}_r = i_F^*(\phi_r F + \phi F_r) = i_F^*(\phi F_r), \qquad  i_{\hat{F}}^*\hat{F}_s = i_F^*(\phi F_s), \qquad  i_{\hat{F}}^*\hat{F}_t = i_F^*(\phi F_t),\\
  i_{\hat{F}}^*\hat{\Delta} = i_F^*(\phi^2 \Delta), \qquad  i_{\hat{F}}^*\hat\lambda_+ = i_F^*(\phi \lambda_+).
 \end{gather*}
 Since $i_F^*$ commutes with $d$,
 \begin{gather*}
 i_{\hat{F}}^*d\left( \frac{\hat{F}_t}{\hat\lambda_+} \right) = i_F^* d\left( \frac{F_t}{\lambda_+} \right), \qquad
 i_{\hat{F}}^*d\left( \frac{\hat{F}_r}{\hat\lambda_+} \right) = i_F^* d\left( \frac{F_r}{\lambda_+} \right).
  \end{gather*}
 From the coframe def\/inition in Lemma \ref{lem:tilde-coframe}, we see that
 \begin{gather*}
 (\hat{\tilde\omega}^1, \dots, \hat{\tilde\omega}^7) = (\tilde\omega^1, (i_F^*\phi) \tilde\omega^2, (i_F^*\phi) \tilde\omega^3, \tilde\omega^4, (i_F^*\phi) \tilde\omega^5, \tilde\omega^6, (i_F^*\phi) \tilde\omega^7 ).
 \end{gather*}
 Hence, from the expression for $\nu_1$ in \eqref{vol1},
 \begin{gather}
 i_{\hat{F}}^* \hat\nu_1 = (i_F^*\phi)^4 i_F^*\nu_1. \label{nu1-scaling}
 \end{gather}
 Since $i_{\hat{F}}^*\hat\Delta = i_F^*(\phi^2 \Delta)$, then by \eqref{vol1} and \eqref{nu1-scaling} we see that necessarily
 \begin{gather*}
 \hat{I}_1 = (i_F^*\phi)^2 I_1.
 \end{gather*}
 Similarly, we f\/ind that
 \begin{gather*}
 \hat{I}_2 = (i_F^*\phi)^2 I_2.
 \end{gather*}

 \section{Generic hyperbolic structure equations}
 \label{app:gen-hyp}

 Starting with the hyperbolic structure equations (c.f.\ Theorem \ref{general-hyp-str-eqns}), we specialize here to the generic case and prove Theorem \ref{generic-hyp-str-eqns} and Corollaries \ref{gamma-cor} and \ref{epsilon-cor}.

 \begin{proof}[Proof of Theorem \ref{generic-hyp-str-eqns}]
 Suppose $\class(M_1)=\class(M_2)=7$ so that $U_1 U_2 \neq 0$.  Let $(\bar\omega^1, \dots,$ $\bar\omega^7) = (\lambda_1 \omega^1, \dots, \lambda_7 \omega^7)$, where
   \begin{gather*}
      \lambda_1 = |U_1|^{-\frac{3}{4}}|U_2|^{-\frac{1}{4}},
     \!\qquad \lambda_2 = |U_1|^{-\frac{1}{2}}|U_2|^{-\frac{1}{2}},
\!\qquad \lambda_3 = |U_1|^{-\frac{3}{4}}|U_2|^{-\frac{1}{4}},
\!\qquad \lambda_4 = |U_1|^{-\frac{1}{4}}|U_2|^{\frac{1}{4}},\\
 \lambda_5 = |U_1|^{-\frac{1}{4}}|U_2|^{-\frac{3}{4}},
\qquad \lambda_6 = 1,
 \qquad \lambda_7 = |U_1|^{-\frac{3}{4}}|U_2|^{-\frac{1}{4}}.
   \end{gather*}
 This normalizes $(U_1,U_2)$ to $(\epsilon_1,\epsilon_2)=({\rm sgn}(U_1),{\rm sgn}(U_2))=(\pm 1,\pm 1)$.  Dropping bars and making the further coframe change
 $\bar\omega^3 = \epsilon_1 \omega^3$, $\bar\omega^6 = \epsilon_1 \omega^6$ normalizes
 $(\epsilon_1,\epsilon_2)$ to $(1,\epsilon)$, where $\epsilon = \epsilon_1 \epsilon_2 = \pm 1$.  Thus, we have obtained the structure equations \eqref{StrEqns123} in Theorem \ref{generic-hyp-str-eqns}.

 Using the structure equations \eqref{StrEqns123}, let us now establish explicit generators for $C(I_F,dM_1)$ and its derived systems.  Let $\{ \Pder{}{1}, \dots, \Pder{}{7} \}$ be the dual basis to $\{ \omega^1, \dots, \omega^7 \}$.  From the structure equations, we have
 \begin{gather*}
{\rm  Char}(I_F,dM_1) = \left\{ \parder{\omega^6}, \parder{\omega^7} \right\} \qRa C(I_F,dM_1) = \{ \omega^1, \omega^2, \omega^3, \omega^4, \omega^5 \}.
 \end{gather*}
 Next, we clearly have $d\omega^1 \equiv 0$, $d\omega^2 \equiv 0$, $d\omega^3 \not\equiv 0$ $\mod C(I_F,dM_1)$
 while
 \begin{gather*}
 0  = d^2\omega^2 \wedge \omega^1 \wedge \omega^2 \wedge \omega^3 \wedge \omega^4 = d\omega^4 \wedge \omega^5 \wedge \omega^1 \wedge \omega^2 \wedge \omega^3 \wedge \omega^4,\\
 0 ,= d^2\omega^2 \wedge \omega^1 \wedge \omega^2 \wedge \omega^3 \wedge \omega^5 = d\omega^5 \wedge \omega^4 \wedge \omega^1 \wedge \omega^2 \wedge \omega^3 \wedge \omega^5
 \end{gather*}
 implies that $d\omega^4 \equiv 0$, $d\omega^5 \equiv 0$ $\mod C(I_F,dM_1)$.
 Thus, $C(I_F,dM_1)^{(1)} = \{ \omega^1, \omega^2, \omega^4, \omega^5 \}$ and
 \begin{gather*}
    d\omega^1 \equiv \omega^3 \wedge \omega^6,
 \qquad d\omega^4 \equiv I \omega^3 \wedge \omega^6 + J \omega^3 \wedge \omega^7,\\
    d\omega^2 \equiv \omega^3 \wedge \omega^7,
  \qquad d\omega^5 \equiv K \omega^3 \wedge \omega^6 + L \omega^3 \wedge \omega^7
\quad \mod\{\omega^1,\omega^2,\omega^4,\omega^5\}.
 \end{gather*}
 We have $d\omega^1 \not\equiv 0$, $d\omega^2 \not\equiv 0$ $\mod C(I_F,dM_1)^{(1)}$.  Redef\/ining
 \begin{gather*}
    \bar\omega^4 = \omega^4 - I \omega^1 - J \omega^2, \qquad
   \bar\omega^5 = \omega^5 - K \omega^1 - L \omega^2
 \end{gather*}
 preserves the structure equations \eqref{StrEqns123}.  Dropping bars, we have
 \begin{gather*}
 d\omega^4 \equiv 0, \qquad  d\omega^5 \equiv 0 \quad \mod C(I_F,dM_1)^{(1)}.
 \end{gather*}
 Hence, $C(I_F,dM_1)^{(2)} = \{ \omega^4, \omega^5 \}$.  Similarly, we describe the derived f\/lag of $C(I_F,dM_2)$.
 \end{proof}

 \begin{proof}[Proof of Corollary \ref{gamma-cor}]  Using the integrability condition $d^2 \omega^1 = 0$, we have
 \begin{gather}
   0  = d^2 \omega^1 \wedge \omega^1 \wedge \omega^2 \wedge \omega^6 = (\omega^7 \wedge \omega^4 - d\omega^6 ) \wedge \omega^3 \wedge \omega^1 \wedge \omega^2 \wedge \omega^6, \label{gamlem1}\\
   0  = d^2 \omega^1 \wedge \omega^1 \wedge \omega^3 \wedge \omega^4 = (\epsilon \omega^5 \wedge \omega^6 - d\omega^4 ) \wedge \omega^2 \wedge \omega^1 \wedge \omega^3 \wedge \omega^4, \label{gamlem2}\\
   0  = d^2 \omega^1 \wedge \omega^1 \wedge \omega^4 \wedge \omega^6 = -(\omega^3 \wedge d\omega^6 + \omega^2 \wedge d\omega^4 ) \wedge \omega^1 \wedge \omega^4 \wedge \omega^6. \label{gamlem3}
 \end{gather}
 Equation \eqref{gamlem1} implies $\gamma^6{}_{47} = -1$ and $\gamma^6{}_{57} = 0$;
 \eqref{gamlem2} implies $ \gamma^4{}_{56} = \epsilon$ and $\gamma^4{}_{57} = 0$;
 \eqref{gamlem3} implies $\gamma^4{}_{35} = \gamma^6{}_{25}$ and $\gamma^4{}_{37} = \gamma^6{}_{27}$.
 Referring to $C(I_F,dM_i)^{(2)}$, in \eqref{charsys}, we know $\gamma^6{}_{25} = \gamma^4{}_{37} = 0$, so we obtain $\gamma^4{}_{35} = \gamma^6{}_{27} = 0$.  Using these values, we have
 \begin{gather*}
 0 = d^2 \omega^4 \wedge \omega^1 \wedge \omega^2 \wedge \omega^4 \wedge \omega^6
 = (\gam{4}{25} - \epsilon \gam{6}{37}) \wedge \omega^5 \wedge \omega^3 \wedge \omega^7 \wedge \omega^1 \wedge \omega^2 \wedge \omega^4 \wedge \omega^6. \tag*{\qed}
 \end{gather*}
\renewcommand{\qed}{}
  \end{proof}

 \begin{proof}[Proof of  Corollary \ref{epsilon-cor}]  We are interested in all coframe changes which preserve the form of the structure equations~\eqref{StrEqns123} and~\eqref{StrEqns4567} with the exception of possibly changing the value of the $\epsilon$ coef\/f\/icient to its negative.  If necessarily the value of the $\epsilon$ coef\/f\/icient is preserved under these coframe changes, then it is a contact invariant.  The change of coframe $\bar\omega^i = R^i{}_j \omega^j$, where
 \begin{gather*}
 R = \mat{ccccccc}{
   -\epsilon & 0 & 0 & 0 & 0 & 0 & 0\\
   0 & 0 & \epsilon & 0 & 0 & 0 & 0\\
   0 & -\epsilon & 0 & 0 & 0 & 0 & 0\\
   0 & 0 & 0 & 0 & 0 & -1 & 0\\
   0 & 0 & 0 & 0 & 0 & 0 & -\epsilon\\
   0 & 0 & 0 & 1 & 0 & 0 & 0\\
   0 & 0 & 0 & 0 & -\epsilon & 0 & 0\\
   }
   \end{gather*}
 interchanges the labelling of $M_1$ and $M_2$ but preserves the structure equations, including the value of $\epsilon$.  Thus, without loss of generality, we may restrict to change of coframe which preserve~$M_1$ and~$M_2$ and hence each $C(I_F,dM_i)^{(k)}$.  Thus, we have
  \begin{gather*}
   \mat{c}{\bar\omega^1 \\ \bar\omega^2\\ \bar\omega^3\\ \bar\omega^4\\ \bar\omega^5\\ \bar\omega^6\\ \bar\omega^7}
   = \mat{ccccccc}{
    \lambda_1 & 0 & 0 & 0 & 0 & 0 & 0\\
   \mu_1 & \lambda_2 & 0 & 0 & 0 & 0 & 0\\
   \mu_2 & 0 & \lambda_3 & 0 & 0 & 0 & 0\\
   0 & 0 & 0 & \lambda_4 & \nu_1 & 0 & 0\\
   0 & 0 & 0 & \mu_3 & \lambda_5 & 0 & 0\\
   0 & 0 & 0 & 0 & 0 & \lambda_6 & \nu_2\\
   0 & 0 & 0 & 0 & 0 & \mu_4 & \lambda_7\\
   }
   \mat{c}{\omega^1 \\ \omega^2\\ \omega^3\\ \omega^4\\ \omega^5\\ \omega^6\\ \omega^7}.
 \end{gather*}
 Using \eqref{StrEqns123} yields
 \begin{gather*}
 d\bar\omega^1  = d(\lambda_1 \omega^1) \equiv
 \lambda_1 (\omega^3 \wedge \omega^6 + \omega^2 \wedge \omega^4) \quad \mod I_F^{(1)}
 \end{gather*}
 and also
 \begin{gather*}
  d\bar\omega^1  \equiv \bar\omega^3 \wedge \bar\omega^6 + \bar\omega^2 \wedge \bar\omega^4 \equiv \lambda_3 \omega^3 \wedge (\lambda_6 \omega^6 + \nu_2 \omega^7) + \lambda_2 \omega^2 \wedge (\lambda_4 \omega^4 + \nu_1 \omega^5)\quad \mod I_F^{(1)} ,
 \end{gather*}
  which implies $\nu_1=\nu_2=0$, $\lambda_1 = \lambda_3 \lambda_6 = \lambda_2 \lambda_4$.   We also have
 \begin{gather*}
 d\bar\omega^2  = d( \mu_1 \omega^1 + \lambda_2 \omega^2 ) \equiv
 \mu_1 \omega^3 \wedge \omega^6 + \lambda_2 (\omega^4 \wedge \omega^5 + \omega^3 \wedge \omega^7)  \quad \mod M_1, \\
 d\bar\omega^3 = d( \mu_2 \omega^1 + \lambda_3 \omega^3 ) \equiv
 \mu_2 \omega^2 \wedge \omega^4 + \lambda_3 (\omega^6 \wedge \omega^7 + \epsilon \omega^2 \wedge \omega^5)  \quad \mod M_2
 \end{gather*}
 and
 \begin{gather*}
  d\bar\omega^2  \equiv \bar\omega^4 \wedge \bar\omega^5 + \bar\omega^3 \wedge \bar\omega^7 \equiv \lambda_4 \omega^4 \wedge\lambda_5 \omega^5 + \lambda_3 \omega^3 \wedge (\mu_4 \omega^6 + \lambda_7 \omega^7) \quad \mod M_1,\\
  d\bar\omega^3  \equiv \bar\omega^6 \wedge \bar\omega^7 + \delta\bar\omega^2 \wedge \bar\omega^5
   \equiv \lambda_6 \omega^6 \wedge\lambda_7 \omega^7 + \delta \lambda_2 \omega^2 \wedge (\mu_3 \omega^4 + \lambda_5 \omega^5)  \quad \mod M_2,
  \end{gather*}
 where $\delta = \pm \epsilon$.  Consequently, we have the system of equations
 \begin{alignat*}{3}
&  \lambda_1 = \lambda_3 \lambda_6 = \lambda_2 \lambda_4, \qquad && \nu_1=\nu_2=0, &\\
 &  \lambda_2 = \lambda_4\lambda_5 = \lambda_3 \lambda_7, \qquad && \mu_1 = \lambda_3 \mu_4, & \\
&  \lambda_3 = \lambda_6\lambda_7 = \delta \epsilon \lambda_2 \lambda_5, \qquad && \mu_2 = \delta \lambda_2 \mu_3.&
 \end{alignat*}
 Then
  \begin{gather*}
   \beta :=\frac{\lambda_6}{\lambda_4} = \frac{\lambda_2}{\lambda_3} = \frac{\delta\epsilon}{\lambda_5}  = \frac{\delta\epsilon\lambda_4}{\lambda_2}, \quad\quad \beta= \frac{\lambda_2}{\lambda_3} = \lambda_7 = \frac{\lambda_3}{\lambda_6},
  \end{gather*}
  and so
  \begin{gather*}
  0 \leq \beta^4 = \left( \frac{\lambda_6}{\lambda_4} \right)\left( \frac{\lambda_2}{\lambda_3} \right)\left( \frac{\delta\epsilon\lambda_4}{\lambda_2} \right)\left( \frac{\lambda_3}{\lambda_6} \right) = \delta\epsilon.
  \end{gather*}
  Consequently, we must have $\delta = \epsilon$ and so $\epsilon$ is a contact invariant.

 Finally, from \eqref{sgnU1}, \eqref{sgnU2}, and the proof of Theorem \ref{generic-hyp-str-eqns} given in Appendix \ref{app:gen-hyp}, we have
 \begin{gather*}
 \epsilon = {\rm sgn}(U_1 U_2) = {\rm sgn}(I_1 I_2).\tag*{\qed}
 \end{gather*}\renewcommand{\qed}{}
 \end{proof}

 \begin{remark}
 Vranceanu \cite{Vranceanu1937} incorrectly asserts that $\epsilon$ can be normalized to~1.  In particular, in the formulas preceding his equation~(4) he writes $\nu = (BC^3)^{\frac{1}{4}}$.  This is invalid if $BC<0$.
 \end{remark}

\section{Isolating maximally and submaximally symmetric structures}
\label{Vranceanu-reduction}

 Starting with the structure equations \eqref{StrEqns123}, \eqref{StrEqns4567} (and not assuming $\gam{5}{56}=0$) we clarify Vran\-ceanu's method of isolating all structures admitting maximal symmetry as well as several structures admitting submaximal symmetry.  The key to understanding Vranceanu's method occurs in a single rather cryptic paragraph on pages~367--368 of \cite{Vranceanu1937}.   In our notation, Vranceanu writes:

 \medskip

 {\em
 In effect, we can remark that if~$a_1$ were equal to $1$ we could reduce~$a_2$, and consequently~$a_3$, to zero, by cancelling in the covariant $d\omega^6$ the term $\gam{6}{46}$ with the help of the coefficient~$a_2$.  This indicates that for~$a_2$,~$a_3$ different than zero, we must have $a_1 \neq 1$, and consequently the systems for which we cannot reduce~$a_2$ and $a_3$ to zero are found amongst those for which $a_1 \neq 1$, $a_2=a_3=0$, and such that we cannot reduce~$a_1$ to~$1$.}

 \medskip

 Let us refer back to the structure equations \eqref{G-lifted-coframe} for the lifted coframe on $\Sigma_7 \times G \ra \Sigma_7$.  Let~$H^0$ denote the subgroup obtained by setting $a_1=1$ in $G^0$.  Since Vranceanu did not consider discrete symmetries, let us
consider the corresponding lifted coframe on $\Sigma_7 \times H^0 \ra \Sigma_7$. We would have structure equations as in \eqref{G-lifted-coframe} except no $\alpha^1$ terms would appear (and $d\alpha^2 = d\alpha^3=0$).  Thus, using Lie algebra valued compatible absorption using $\alpha^2$ and $\alpha^3$ only, $\gam{6}{46}$ would be a torsion coef\/f\/icient.  With respect to the group $H^0$, we have the transformation laws
 \begin{gather}
 \hgam{6}{46} = \gam{6}{46} - \gam{6}{56} a_3 + a_2, \qquad  \hat\gamma^6{}_{56} = \gamma^6{}_{56}. \label{gamma646}
 \end{gather}
 By setting $a_2 = \gam{6}{56} a_3 - \gam{6}{46}$, we can normalize $\hgam{6}{46}=0$.  In the group reduced from $H^0$, we only obtain $a_2 = \gam{6}{56} a_3$, so it is unclear whether $a_2$, $a_3$ can {\em both} be reduced to 0.  Thus, Vranceanu's claim may be true, but if so it is certainly not obvious.  Are there any structures for which $a_1$ can be reduced to 1, but at least one of $a_2$, $a_3$ are nonzero?

 The second sentence in the paragraph appears to be quite cryptic and at f\/irst sight appears even self-contradictory.  One can make sense of this as follows: Vranceanu sets out to f\/ind all structures for which $K^0 = \{ {\rm  diag} (a_1{}^2,a_1,a_1,a_1,1,a_1,1) : a_1 > 0 \} \subset G$ {\em cannot} be reduced to the identity.  Thus, he is considering the restricted equivalence problem with respect to the subgroup~$K^0$, and a lifted coframe on the bundle $\Sigma_7 \times K^0 \ra \Sigma_7$ satisfying the structure equations
 \begin{gather}
   d\hat\omega^1  = 2\alpha^1 \wedge \hat\omega^1 + \hat\omega^3 \wedge \hat\omega^6 + \hat\omega^2 \wedge \hat\omega^4 + \eta_1 \wedge \hat\omega^1,\nonumber\\
   d\hat\omega^2  = \alpha^1 \wedge \hat\omega^2 + \hat\omega^4 \wedge \hat\omega^5 + \hat\omega^3 \wedge \hat\omega^7 + \eta_{21} \wedge \hat\omega^1 + \eta_{22} \wedge \hat\omega^2,\nonumber\\
   d\hat\omega^3  = \alpha^1 \wedge \hat\omega^3 + \hat\omega^6 \wedge \hat\omega^7 + \epsilon\hat\omega^2 \wedge \hat\omega^5 + \eta_{31} \wedge \hat\omega^1 + \eta_{33} \wedge \hat\omega^3,\nonumber\\
   d\hat\omega^4  = \alpha^1 \wedge \hat\omega^4 + \eta_{41} \wedge \hat\omega^1 + \eta_{42} \wedge \hat\omega^2 + \eta_{44} \wedge \hat\omega^4 + \epsilon\hat\omega^5 \wedge \hat\omega^6, \label{K-lifted-coframe}\\
   d\hat\omega^5  = \eta_{51} \wedge \hat\omega^1 + \eta_{52} \wedge \hat\omega^2 + \eta_{54} \wedge \hat\omega^4 + \eta_{55} \wedge \hat\omega^5,\nonumber\\
   d\hat\omega^6  = \alpha^1 \wedge \hat\omega^6 + \eta_{61} \wedge \hat\omega^1 + \eta_{63} \wedge \hat\omega^3 + \eta_{66} \wedge \hat\omega^6 - \hat\omega^4 \wedge \hat\omega^7,\nonumber\\
   d\hat\omega^7  = \eta_{71} \wedge \hat\omega^1 + \eta_{73} \wedge \hat\omega^3 + \eta_{76} \wedge \hat\omega^6 + \eta_{77} \wedge \hat\omega^7,\nonumber\\
   d\alpha^1  = 0.\nonumber
   \end{gather}
 where $\eta_i$ are semibasic with respect to the canonical projection $\Sigma_7 \times K^0 \ra \Sigma_7$.  Note that under the $K^0$-action, all coef\/f\/icients transform by a scaling action:
 \begin{gather*}
 \hgam{i}{jk} =  (a_1)^p \gam{i}{jk}
 \end{gather*}
 for some scaling weight $p \in \R$.  (In other words, by setting $a_2=a_3=0$ in the general transformation formulas under $G^0$, c.f.\ last sentence in Vranceanu's paragraph above.)  For specif\/ic $i$,~$j$,~$k$, if $\gam{i}{jk}\neq 0$, then we can normalize to $\hgam{i}{jk}=\pm 1$ by setting $a_1 = |\gam{i}{jk}|^{-1/p}$, thereby reducing~$K^0$ to the identity.  Since we are assuming that $K^0$ cannot be reduced, Vranceanu's conclusion would be that $\gam{i}{jk}=0$.  A technical assumption that should be made here is that if $p \neq 0$, we are considering $\gam{i}{jk}$ to be a {\em constant} torsion coef\/f\/icient.  For example, if $\gam{i}{jk}\neq 0$ (as functions) but vanishes at a point, then $K$ cannot be normalized since $a_1$ must be nonzero.  In summary, we have:

 \begin{lemma}
 Suppose $\hat\kappa = (a_1)^p \kappa$ is a torsion coefficient with respect to $K^0$.  If (1) $K^0$ cannot be reduced to the identity, (2) $p \neq 0$, and (3) $\kappa$ is a constant, then $\kappa=0$.
 \label{vanishing-lemma}
 \end{lemma}

 \begin{note} We do {\em not} need to assume that {\em all} torsion terms are constant.  Only those torsion terms with a nontrivial scaling action by $K^0$ are assumed to be constant.
 \end{note}

 Vranceanu f\/inishes the paragraph with the consequences of this simplif\/ication:

 \medskip

 {\em However, it is easy to see that for these systems we must have in the formulas \eqref{StrEqns4567}, all null coefficients, except for $\gamma^4{}_{25}$, $\gamma^4{}_{27}$, $\gamma^4{}_{53}$, $\gamma^4{}_{56}=1$, $\gamma^5{}_{57}$, $\gamma^6{}_{35}$, $\gamma^6{}_{37}$, $\gamma^6{}_{72}$, $\gamma^6{}_{74}=1$, $\gamma^7{}_{75}$ and likewise one must have $\gamma^2{}_{3\alpha}=0 $ $(\alpha \neq 2)$, $\gamma^3{}_{2\alpha}=0 $ $(\alpha \neq 3)$.}

 \medskip

 \begin{note} The values $\gam{4}{53}=\gam{6}{72}=0$, $\gam{4}{56} = \epsilon$, $\gam{6}{47} = -1$ were established in Corollary \ref{gamma-cor}.  Also, the f\/inal assertions in Vranceanu's sentence above should read $\gam{2}{1\alpha}$ and $\gam{3}{1\alpha}$ respectively~-- see explanations below.
 \end{note}

 We give an outline of these details.  With respect to the subgroup $K^0$, the following are torsion terms (c.f.\ \eqref{K-lifted-coframe}):
 \begin{gather*}
 \hgam{2}{1k} \quad (k\neq 2), \qquad
 \hgam{3}{1k} \quad (k\neq 3), \qquad
 \hgam{4}{1k}, \  \hgam{4}{2k} \quad (k \neq 4), \qquad
 \hat\gamma^6{}_{1k}, \ \hat\gamma^6{}_{3k} \quad (k\neq 6),\\
 \hgam{5}{1k},  \ \hgam{5}{2k},  \ \hgam{5}{4k}, \
  \hgam{7}{1k}, \  \hgam{7}{3k}, \ \hgam{7}{6k} \quad (k \mbox{ arbitrary}), \qquad
 \hat\gamma^5{}_{35}, \ \hat\gamma^7{}_{27}, \ \hat\gamma^7{}_{47}.
 \end{gather*}
 The scaling weight for all these terms is nonzero with the exception of $\hgam{4}{25}$, $\hgam{4}{27}$, $\hgam{6}{35}$, $\hgam{6}{37}$ for which $p=0$.  By Lemma \ref{vanishing-lemma}, all of the corresponding coef\/f\/icients must vanish with the exception of  $\gam{4}{25}$, $\gam{4}{27}$, $\gam{6}{35}$, $\gam{6}{37}$.  Vranceanu did not make use of the following torsion terms:
 \begin{gather*}
 2\hgam{2}{k2} - \hgam{1}{k1} \quad (k \neq 1,2),
\qquad  2\hgam{3}{k3} - \hgam{1}{k1} \quad (k \neq 1,3),
 \qquad 2\hgam{4}{k4} - \hgam{1}{k1} \quad (k \neq 1,4),\\
 2\hgam{6}{k6} - \hgam{1}{k1} \quad (k \neq 1,6),
\qquad \hgam{3}{13} - \hgam{2}{12}, \qquad \hgam{4}{14} - \hgam{2}{12}, \qquad
 \hgam{6}{16} - \hgam{2}{12}.
 \end{gather*}
 All of these terms have nonzero scaling weight (and hence must vanish) with the exception of
 \begin{gather*}
 2\hgam{i}{ki} - \hgam{1}{k1} \quad\mbox{(no sum)},  \qquad i=2,3,4,6, \quad k=5,7.
 \end{gather*}
 Thus, the structure equations take the form
 \begin{gather}
   d\omega^1  = \omega^3 \wedge \omega^6 + \omega^2 \wedge \omega^4 + \omega^1 \wedge \eta_1,\nonumber\\
   d\omega^2  = \omega^4 \wedge \omega^5 + \omega^3 \wedge \omega^7 + \omega^2 \wedge \eta_2,\nonumber\\
   d\omega^3  = \omega^6 \wedge \omega^7 + \epsilon \omega^2 \wedge \omega^5 + \omega^3 \wedge \eta_3,\nonumber\\
   d\omega^4  = \epsilon \omega^5 \wedge \omega^6 + \omega^2 \wedge (\gam{4}{25} \omega^5 + \gam{4}{27} \omega^7) + \omega^4 \wedge \eta_4,\label{StrEqn-pre-eta1-red}\\
   d\omega^5  = m \omega^5 \wedge \omega^7,\nonumber\\
   d\omega^6  = - \omega^4 \wedge \omega^7 + \omega^3 \wedge (\gam{6}{35} \omega^5 + \gam{6}{37} \omega^7) + \omega^6 \wedge \eta_6,\nonumber\\
   d\omega^7  = n \omega^5 \wedge \omega^7,\nonumber
   \end{gather}
 where $\eta_k = \gam{k}{k\ell} \omega^\ell$ (no sum on $k$), and
 \begin{gather*}
  2\gam{4}{46} = 2\gam{3}{36} = 2\gam{2}{26} = \gam{1}{16},
\qquad 2\gam{6}{63} = 2\gam{4}{43} = 2\gam{2}{23} = \gam{1}{13}, \\
  2\gam{6}{64} = 2\gam{3}{34} = 2\gam{2}{24} = \gam{1}{14},
 \qquad 2\gam{6}{62} = 2\gam{4}{42} = 2\gam{3}{32} = \gam{1}{12},\\
   \gam{6}{16} = \gam{4}{14} = \gam{3}{13} = \gam{2}{12},
 \end{gather*}
 or more succinctly,
 \begin{gather}
  \gam{k}{k\ell} = \frac{1}{2} \gam{1}{1\ell} \quad \mbox{(no sum),}\qquad k,\ell =2,3,4,6, \label{gamma-red1}\\
  \gam{k}{1k} = \gam{2}{12} \quad \mbox{(no sum),} \qquad k =3,4,6. \label{gamma-red2}
 \end{gather}

 \begin{lemma}
 There exists an admissible change of coframe so that $\eta_1=0$ and $\gam{2}{12}=0$.
 \end{lemma}
 \begin{proof} Let us evaluate some integrability conditions.
 \begin{gather}
 0  = d^2 \omega^1
 = ( \eta_1 - \eta_2 - \eta_4 ) \wedge \omega^2 \wedge \omega^4 + (\eta_1 - \eta_3 - \eta_6) \wedge \omega^3 \wedge \omega^6 - \omega^1 \wedge d\eta_1. \label{eta1-int0}
 \end{gather}
 Isolating terms that are in the ideal generated by $\omega^1$ and using \eqref{gamma-red2}, we have
 \begin{gather}
 0  = \omega^1 \wedge ( 2\gam{2}{12} d\omega^1 -  d\eta_1). \label{eta1-int1}
 \end{gather}
 Note that for $k=2,3,4,6$,
 \begin{gather}
0  = d^2 \omega^k \wedge \omega^5 \wedge \omega^7 = - \omega^k \wedge d\eta_k \wedge \omega^5 \wedge \omega^7 \qquad\mbox{(no sum on $k$)} \nonumber\\
 \phantom{0}{} = - \omega^k \wedge d( \gam{k}{k\ell} \omega^\ell ) \wedge \omega^5 \wedge \omega^7
 = - \omega^k \wedge d\left( \gam{2}{21} \omega^1 + \frac{1}{2} \eta_1 \right) \wedge \omega^5 \wedge \omega^7 \nonumber\\
 \phantom{0}{} = - \omega^k \wedge \left( d\gam{2}{21} \wedge \omega^1 + \gam{2}{21} d\omega^1 + \frac{1}{2} d\eta_1 \right) \wedge \omega^5 \wedge \omega^7 .\label{eta1-int2}
 \end{gather}
 Thus, \eqref{eta1-int1} and \eqref{eta1-int2} imply
 \begin{gather*}
 d\eta_1 = d(2\gam{2}{12} \omega^1) + \omega^1 \wedge (f_1 \omega^5 + f_2 \omega^7)
 \end{gather*}
 for some functions $f$, $g$.  Applying $d$ and wedging with $\omega^1$ yields
 \begin{gather*}
 0 = \omega^1 \wedge d\omega^1 \wedge (f_1 \omega^5 + f_2 \omega^7) = \omega^1 \wedge (\omega^2 \wedge \omega^4 + \omega^3 \wedge \omega^6) \wedge (f_1 \omega^5 + f_2 \omega^7)
 \end{gather*}
 and so we must have $f_1 = f_2 = 0$.  By Poincar\'{e}'s lemma, we have
 \begin{gather*}
 \eta_1 = 2\gam{2}{12} \omega^1 + 2dh
 \end{gather*}
 for some function $h$.  Def\/ine the change of coframe
 \begin{gather*}
(\bar\omega^1,\bar\omega^2,\bar\omega^3,\bar\omega^4,\bar\omega^6)  = (a_1{}^2 \omega^1, a_1 \omega^2, a_1 \omega^3, a_1 \omega^4, a_1 \omega^6 ), \qquad \mbox{where} \quad a_1 = e^h.
 \end{gather*}
 Then
 \begin{gather*}
 d\bar\omega^1  = 2e^{2h} dh \wedge \omega^1 + e^{2h} d\omega^1
 = 2e^{2h} dh \wedge \omega^1 + e^{2h} (\omega^2 \wedge \omega^4 + \omega^3 \wedge \omega^6 + \omega^1 \wedge \eta_1) \\
 \phantom{ d\bar\omega^1}{}  = 2 dh \wedge \bar\omega^1 + \bar\omega^2 \wedge \bar\omega^4 + \bar\omega^3 \wedge \bar\omega^6 + \bar\omega^1 \wedge \eta_1 \\
  \phantom{ d\bar\omega^1}{} = \bar\omega^2 \wedge \bar\omega^4 + \bar\omega^3 \wedge \bar\omega^6.
 \end{gather*}
 Dropping bars, we have a new coframe satisfying \eqref{StrEqn-pre-eta1-red} with $\eta_1=0$.  Moreover, using \eqref{eta1-int1}, we have that $\gam{2}{12}=0$.
 \end{proof}

 From \eqref{eta1-int0}, we have that
 \begin{gather*}
 \gam{4}{45} = -\gam{2}{25}, \qquad
 \gam{4}{47} = -\gam{2}{27}, \qquad
 \gam{6}{65} = -\gam{3}{35}, \qquad
 \gam{6}{67} = -\gam{3}{37},
 \end{gather*}
 and using \eqref{gamma-red1}, \eqref{gamma-red2}, the structure equations take the form
 \begin{gather*}
   d\omega^1  = \omega^3 \wedge \omega^6 + \omega^2 \wedge \omega^4, \\
   d\omega^2  = \omega^4 \wedge \omega^5 + \omega^3 \wedge \omega^7 + \omega^2 \wedge ( \gam{2}{25} \omega^5 + \gam{2}{27} \omega^7 ), \\
   d\omega^3  = \omega^6 \wedge \omega^7 + \epsilon \omega^2 \wedge \omega^5 + \omega^3 \wedge ( \gam{3}{35} \omega^5 + \gam{3}{37} \omega^7 ),\\
   d\omega^4  = \epsilon \omega^5 \wedge \omega^6 + \omega^2 \wedge (\gam{4}{25} \omega^5 + \gam{4}{27} \omega^7) - \omega^4 \wedge (\gam{2}{25} \omega^5 + \gam{2}{27} \omega^7), \\
   d\omega^5  = m \omega^5 \wedge \omega^7, \\
   d\omega^6  = - \omega^4 \wedge \omega^7 + \omega^3 \wedge (\gam{6}{35} \omega^5 + \gam{6}{37} \omega^7) - \omega^6 \wedge (\gam{3}{35} \omega^5 + \gam{3}{37} \omega^7), \\
   d\omega^7  = n \omega^5 \wedge \omega^7.
   \end{gather*}

 Further integrability conditions reveal
 \begin{gather*}
  0 = d^2 \omega^5 = dm \wedge \omega^5 \wedge \omega^7 \qRa dm = m_5 \omega^5 + m_7 \omega^7,\\
 0 = d^2 \omega^7 = dn \wedge \omega^5 \wedge \omega^7 \qRa dn = n_5 \omega^5 + n_7 \omega^7.
 \end{gather*}
 Next,
 \begin{gather*}
  0 =d^2 \omega^2 \wedge \omega^2 = (n - \gam{3}{35} + \gam{2}{25}) \omega^2 \wedge \omega^3 \wedge \omega^5 \wedge \omega^7 + (m - 2\gam{2}{27}) \omega^2 \wedge \omega^4 \wedge \omega^5 \wedge \omega^7, \\
 0 =d^2 \omega^3 \wedge \omega^3 = \epsilon(-m + \gam{3}{37} - \gam{2}{27}) \omega^2 \wedge \omega^3 \wedge \omega^5 \wedge \omega^7 - (n + 2\gam{3}{35}) \omega^3 \wedge \omega^5 \wedge \omega^6 \wedge \omega^7, \\
  \hspace{-0.25in}\Rightarrow \ \gam{2}{25} = -\frac{3n}{2}, \qquad \gam{2}{27} = \frac{m}{2},  \qquad \gam{3}{35} = -\frac{n}{2}, \qquad \gam{3}{37} = \frac{3m}{2},\\
 0 = d^2 \omega^2 = ( mn + \epsilon - \gam{4}{27} ) \omega^2 \wedge \omega^5 \wedge \omega^7 + \frac{1}{2} \omega^2 \wedge (3 dn \wedge \omega^5 - dm \wedge \omega^7),\\
 0 = d^2 \omega^3  = ( \gam{6}{35} - mn - \epsilon) \omega^3 \wedge \omega^5 \wedge \omega^7 + \frac{1}{2} \omega^3 \wedge (dn \wedge \omega^5 - 3dm \wedge \omega^7),\\
 0 = d^2 \omega^4 \wedge \omega^1 \wedge \omega^2 \wedge \omega^4 \wedge \omega^6 = (\epsilon\gam{6}{37} - \gam{4}{25}) \omega^3 \wedge \omega^5 \wedge \omega^7 \wedge \omega^1 \wedge \omega^2 \wedge \omega^4 \wedge \omega^6,\\
  \hspace{-0.25in}\Rightarrow   \  \gam{4}{25} = B = \epsilon\gam{6}{37}, \qquad  \gam{4}{27} = mn + \epsilon - \frac{3}{2} n_7 - \frac{1}{2} m_5, \qquad
  \gam{6}{35} = mn + \epsilon + \frac{1}{2} n_7 + \frac{3}{2} m_5.
 \end{gather*}
  The f\/inal two integrability conditions are
 \begin{gather*}
 0  =d^2 \omega^4 = (-4n\Delta_1 - 2mB + 6 n n_7 + n m_5  - mn_5) \omega^2 \wedge \omega^5 \wedge \omega^7 + \omega^2 \wedge \omega^5 \wedge dB\\
\phantom{0  =d^2 \omega^4 =}{}  - \frac{1}{2} \omega^2 \wedge \omega^7 \wedge ( dm_5 + 3 dn_7 ),\\
 0  =d^2 \omega^6 =  ( - 4m\Delta_1 - 2n\epsilon B - 6 mm_5 - m n_7 + nm_7)\omega^3 \wedge \omega^5 \wedge \omega^7 + \epsilon \omega^3 \wedge \omega^7 \wedge dB\\
\phantom{0  =d^2 \omega^6 =}{}   + \frac{1}{2} \omega^3 \wedge \omega^5 \wedge ( 3dm_5 + dn_7 ),
 \end{gather*}
which implies
 \begin{gather*}
 0  = d^2 \omega^4 \wedge \omega^7 = \omega^2 \wedge \omega^5 \wedge dB \wedge \omega^7,\\
 0  = d^2 \omega^6 \wedge \omega^5 = \epsilon \omega^3 \wedge \omega^7 \wedge dB \wedge \omega^5,
 \end{gather*}
 and hence $d^2\omega^4 = d^2 \omega^6 = 0$ is equivalent to
 \begin{gather*}
 dB  = \epsilon \left( - 4m\Delta_1 - 2n\epsilon B - 6 mm_5 - m n_7 + nm_7 + \frac{3}{2} m_{57} + \frac{1}{2} n_{77}\right)\omega^5 \\
\phantom{dB  =}{} + \left(4n\Delta_1 + 2mB - 6 n n_7 - n m_5 + mn_5 - \frac{1}{2} m_{55} - \frac{3}{2} n_{75}\right) \omega^7.
 \end{gather*}
 This reduces the structure equations to those given in Theorem~\ref{mn-str-eqns}.


 \section{Parametrization of maximally symmetric structures}
 \label{maxsym-param}

 While in general it is very dif\/f\/icult to go from structure equations specif\/ied in an abstract coframe to a coframe described parametrically (i.e.\ in local coordinates), Vranceanu succeeded in doing this for the maximally symmetric structure equations \eqref{9dim-streqns}.  Certain steps in his computation were unjustif\/ied or contained errors (e.g.\ Vranceanu missed the contact invariant $\epsilon$ in his structure
equations and hence his normal forms), and so we provide a more complete outline here.

 Let us consider a coframe $\{ \omega^i \}_{i=1}^7$ on $\Sigma_7$ satisfying \eqref{9dim-streqns}, and let us consider $\omega^1$ in its canonical form
 \begin{gather}
 \omega^1 = dz - p dx - qdy. \label{omega1-canonical}
 \end{gather}

 \subsection*{$\boldsymbol{\omega^5}$ and $\boldsymbol{\omega^7}$}

 The subsystem $\{ \omega^5, \omega^7 \}$ is completely integrable, so by the Frobenius theorem there exist local functions $u$, $v$ on $\Sigma_7$ such that $\{ du, dv \}$ generate the same subsystem.  We introduce $u$, $v$ as follows.  By using the $d\omega^5$ and $d\omega^7$ equations,
 \begin{gather*}
 0  = d\left( \frac{\epsilon}{ m} \omega^5 +  m \omega^7 \right) \qRa \frac{\epsilon}{ m} \omega^5 +  m \omega^7 = \frac{du}{u}.
 \end{gather*}
 Consequently,
 \begin{gather*}
 d\omega^7  = -\frac{du}{u} \wedge \omega^7 \qRa d(u \omega^7) = 0,
 \end{gather*}
and hence
 \begin{gather}
 \omega^7 = \frac{dv}{u}, \qquad \omega^5 = \frac{\epsilon  m( du -  m dv)}{u}. \label{omega57}
 \end{gather}

  \subsection*{An ansatz for $\boldsymbol{\omega^2}$, $\boldsymbol{\omega^3}$, $\boldsymbol{\omega^4}$, $\boldsymbol{\omega^6}$}

 We can now determine an ansatz for $\omega^2$, $\omega^3$, $\omega^4$, $\omega^6$.  Plugging \eqref{omega1-canonical}, into the $d\omega^1$ equation, we obtain
 \begin{gather*}
 dx \wedge dp + dy \wedge dq = \omega^2 \wedge \omega^4 + \omega^3 \wedge \omega^6.
 \end{gather*}
 Taking the interior product with $\parder{u}$, we f\/ind that
 \begin{gather*}
 0  = (i_{\partial_u} \omega^2) \omega^4 - (i_{\partial_u} \omega^4) \omega^2 + (i_{\partial_u} \omega^3) \omega^6 - (i_{\partial_u} \omega^6) \omega^3,
 \end{gather*}
 and hence $i_{\partial_u} \omega^k = 0$ and similarly $i_{\partial_v} \omega^k = i_{\partial_z} \omega^k = 0$ for $k=2,3,4,6$.  Thus,
 \begin{gather}
 \omega^k = A^k dp + B^k dq + C^k dx + D^k dy. \label{omega2346-ansatz}
 \end{gather}
 Vranceanu asserts without justif\/ication that the coef\/f\/icients $A^k$, $B^k$, $C^k$, $D^k$ can be taken to only depend on $u$, $v$.  In Lemma \ref{ansatz-lemma}, we show that we can reduce to dependence on $x$, $y$, $u$, $v$, but it is unclear why a further simplif\/ication is a priori necessary.

 \subsection*{$\boldsymbol{\omega^4}$ and $\boldsymbol{\omega^6}$}

 We observe from the $d\omega^4$ and $d\omega^6$ equations, that
 \begin{gather*}
 d\left( \omega^4 +  m \omega^6 \right) = -\frac{3}{2} \left( \omega^4 +  m \omega^6 \right) \wedge \frac{du}{u} \qRa 0 = d\big( u^{-3/2} ( \omega^4 +  m \omega^6) \big),
 \end{gather*}
 and we make the choice of introducing the $x$-coordinate into the coframe via
 \begin{gather*}
 \omega^4 +  m \omega^6 = u^{3/2} dx.
 \end{gather*}
 Now substitute into the $d\omega^6$ equation and simplify to get
 \begin{gather*}
d\left( \frac{\omega^6}{\sqrt{u}} \right)  = d(vdx)
 \end{gather*}
 and we make the choice of introducing the $y$-coordinate into the coframe via
 \begin{gather*}
 \omega^6 = -\sqrt{u} (dy-vdx), \qquad \omega^4 = u^{3/2} dx +  m \sqrt{u}(dy - vdx).
 \end{gather*}

 \subsection*{$\boldsymbol{\omega^2}$ and $\boldsymbol{\omega^3}$}

 Using the ansatz \eqref{omega2346-ansatz}, we write
 \begin{gather*}
 \omega^2  = \sigma_1 \omega^6 + \mu \omega^4 + \alpha_1 dp + \beta_1 dq,\\
 \omega^3  = \sigma_2 \omega^4 + \rho \omega^6 + \alpha_2 dp + \beta_2 dq.
 \end{gather*}
 Plugging these and \eqref{omega1-canonical} into $d\omega^1$,
 we obtain $\sigma_1 = \sigma_2 =: \sigma$ and
 \begin{gather*}
 dx \wedge dp + dy \wedge dq
  = -(\sqrt{u}(u -  m v) \alpha_1 + v\sqrt{u} \alpha_2) dx \wedge dp + \sqrt{u}( m \alpha_1 - \alpha_2) dp \wedge dy, \\
 \phantom{dx \wedge dp + dy \wedge dq=}{}  + (\sqrt{u}(u- m v)\beta_1 + v\sqrt{u} \beta_2) dq \wedge dx - \sqrt{u}( m \beta_1 - \beta_2) dy \wedge dq,
 \end{gather*}
 which can be solved to obtain
 \begin{gather}
 \omega^2  = \sigma \omega^6 + \mu \omega^4 - u^{-3/2} (dp + vdq), \label{omega2-ansatz}\\
 \omega^3   = \sigma \omega^4 + \rho \omega^6 -  m u^{-3/2} (dp + v dq) + u^{-1/2} dq, \label{omega3-ansatz}
 \end{gather}
 for some functions $\sigma$, $\mu$, $\rho$.

\begin{lemma} \label{ansatz-lemma} $\sigma$, $\mu$, $\rho$ only depend on $x$, $y$, $u$, $v$.
\end{lemma}
\begin{proof}
 Since $\omega^4$, $\omega^6$ depend only on $dx$, $dy$ and $\omega^5$, $\omega^7$ depend only on $du$, $dv$, it suf\/f\/ices to show that $d\sigma$, $d\mu$, $d\rho$ depend only on $\omega^4$, $\omega^5$, $\omega^6$, $\omega^7$.
The structure equations \eqref{9dim-streqns} indicate that $d\omega^k \equiv 0$ $\mod \{ \omega^5, \omega^7 \}$ for $k=2,3,4,6$.  Applying $d$ to \eqref{omega2-ansatz}, \eqref{omega3-ansatz}, we obtain
 \begin{gather*}
 d\omega^2 \equiv d\sigma \wedge \omega^6 + d\mu \wedge \omega^4 \quad \mod \{ \omega^5, \omega^7 \},\\
 d\omega^3 \equiv d\sigma \wedge \omega^4 + d\rho \wedge \omega^6 \quad \mod \{ \omega^5, \omega^7 \},
 \end{gather*}
 so necessarily,
 \begin{align*}
 d\sigma &\equiv \sigma_4 \omega^4 + \sigma_6 \omega^6 \quad\mod \{ \omega^5, \omega^7 \},\\
 d\mu &\equiv \mu_4 \omega^4 + \sigma_4 \omega^6 \quad\mod \{ \omega^5, \omega^7 \},\\
 d\rho & \equiv \sigma_6 \omega^4 + \rho_6 \omega^6 \quad\mod \{ \omega^5, \omega^7 \}.\tag*{\qed}
 \end{align*}\renewcommand{\qed}{}
 \end{proof}

 We now make the simplifying search for solutions where $\sigma$, $\mu$, $\rho$ are functions of $u$, $v$ only.  Thus,
 \begin{gather*}
 d\sigma = \Pder{\sigma}{5} \omega^5 + \Pder{\sigma}{7} \omega^7, \qquad
 d\mu = \Pder{\mu}{5} \omega^5 + \Pder{\mu}{7} \omega^7, \qquad
 d\rho = \Pder{\rho}{5} \omega^5 + \Pder{\rho}{7} \omega^7.
 \end{gather*}
 Then substitution of \eqref{omega2-ansatz}, \eqref{omega3-ansatz} into the structure equations for $\omega^2$ and $\omega^3$ yields
 \begin{gather*}
   d\omega^2 = \left( \Parder{\sigma}{\omega^5} + \frac{2\epsilon\sigma}{ m} + \epsilon \mu \right) \omega^5 \wedge \omega^6 + \left( \Parder{\sigma}{\omega^7} + 2 m \sigma + \rho \right) \omega^7 \wedge \omega^6 + \left( \Parder{\mu}{\omega^5} + \frac{3\epsilon\mu}{ m}\right)\omega^5 \wedge \omega^4 \\
  \phantom{d\omega^2 =}{}  + \left( \Parder{\mu}{\omega^7} + 2\sigma +  m \mu \right) \omega^7 \wedge \omega^4 + \omega^3 \wedge \omega^7 - \left( \frac{3\epsilon}{2 m} \omega^5 + \frac{ m}{2} \omega^7 \right) \wedge \omega^2,\\
   d\omega^3  = \left( \Parder{\sigma}{\omega^5} + \frac{2\epsilon\sigma}{ m} + \epsilon \mu \right) \omega^5 \wedge \omega^4 + \left( \Parder{\sigma}{\omega^7} + 2 m\sigma + \rho \right) \omega^7 \wedge \omega^4 \\
   \phantom{d\omega^3  =}{}
   + \left( \Parder{\rho}{\omega^5} + 2 \epsilon \sigma + \frac{\epsilon \rho}{ m} \right) \omega^5 \wedge \omega^6
 + \left( \Parder{\rho}{\omega^7} + 3 m \rho \right) \omega^7 \wedge \omega^6 + \epsilon \omega^2 \wedge \omega^5 \\
\phantom{d\omega^3  =}{}
   - \left( \frac{\epsilon}{2 m} \omega^5 + \frac{3 m}{2} \omega^7 \right) \wedge \omega^3.
 \end{gather*}
 Comparison with the structure equations yields
 \begin{alignat}{4}
   &\Parder{\mu}{\omega^5} = - \frac{3\epsilon\mu}{ m} -1, \qquad
   && \Parder{\sigma}{\omega^5} = -\frac{2\epsilon\sigma}{ m} - \epsilon \mu, \qquad
   &&\Parder{\rho}{\omega^5} = - 2 \epsilon \sigma - \frac{\epsilon}{ m} \rho, & \label{sigma-mu-rho1}\\
   &\Parder{\mu}{\omega^7} = - 2\sigma -  m \mu, \qquad
   &&\Parder{\sigma}{\omega^7} = - 2 m\sigma - \rho, \qquad
   &&\Parder{\rho}{\omega^7} = - 3 m \rho - 1 . & \label{sigma-mu-rho2}
 \end{alignat}
 For $\gam{i}{jk}$ def\/ined as in \eqref{gamma-defn}, we have
 \begin{gather*}
 \left[\Parder{}{\omega^j},\Parder{}{\omega^k}\right] = - \gam{i}{jk} \Parder{}{\omega^i} \qRa
\left[ \Parder{}{\omega^5}, \Parder{}{\omega^7} \right] = - m \Parder{}{\omega^5} + \frac{\epsilon}{ m} \Parder{}{\omega^7}.
 \end{gather*}
 Using this commutator relation, we can verify that all integrability conditions for the system \eqref{sigma-mu-rho1}, \eqref{sigma-mu-rho2} are satisf\/ied.  Using \eqref{omega57}, we have for $F=F(u,v)$,
 \begin{gather*}
 dF = F_u du + F_v dv = \frac{u F_u}{\epsilon m} \omega^5 + u( m F_u +F_v) \omega^7,
 \end{gather*}
 and so
 \begin{gather*}
 \Pder{F}{5}  = \frac{uF_u}{\epsilon m}, \qquad \Pder{F}{7} = ( m F_u + F_v)u = u F_v + \epsilon m^2 \Pder{F}{5}.
 \end{gather*}
 We use these expressions for the coframe derivatives to rewrite the system \eqref{sigma-mu-rho1}, \eqref{sigma-mu-rho2} and obtain  expressions for $\sigma$, $\mu$, $\rho$.  The three equations in \eqref{sigma-mu-rho1} are respectively equivalent to
 \begin{gather*}
  (u^3 \mu)_u = -\epsilon m u^2 \quad \Rightarrow \quad \mu = -\frac{\epsilon m}{3} +\frac{a_1(v)}{u^3}, \\
   (u^2\sigma)_u = - m u \mu \quad \Rightarrow \quad \sigma = \frac{\epsilon m^2}{6} + \frac{ m a_1(v)}{u^3} + \frac{a_2(v)}{u^2}, \\
 (u\rho)_u = - 2 m\sigma \quad \Rightarrow \quad \rho = -\frac{\epsilon m^3}{3} + \frac{ m^2 a_1(v)}{u^3} + \frac{2 m a_2(v)}{u^2} + \frac{a_3(v)}{u}.
 \end{gather*}
 Using these expressions for $\sigma$, $\mu$, $\rho$, the three equations in \eqref{sigma-mu-rho2} are respectively equivalent to
 \begin{gather}
 a_1'(v) = - 2 a_2(v), \qquad  a_2'(v) = -a_3(v), \qquad  a_3'(v) = -\alpha, \label{a1a2a3}
 \end{gather}
 where $\alpha = 1 - \epsilon m^4$.
 The general solution of \eqref{a1a2a3} is
 \begin{gather*}
 a_3(v) = -\alpha v + c_3,
\qquad a_2(v) = \frac{1}{2} \alpha v^2 - c_3 v + c_2,
\qquad a_1(v) = -\frac{1}{3}\alpha v^3 + c_3 v^2 - 2 c_2 v + c_1,
 \end{gather*}
 for arbitrary constants $c_1$, $c_2$, $c_3$.  Taking the simplest choice $c_1=c_2=c_3=0$ yields
 \begin{gather*}
 \mu = -\frac{\epsilon m}{3} - \frac{\alpha v^3}{3u^3},
\qquad \sigma = \frac{\epsilon m^2}{6} - \frac{ m \alpha v^3}{3u^3} + \frac{\alpha v^2}{2u^2},
\qquad \rho = -\frac{\epsilon m^3}{3} - \frac{ m^2 \alpha v^3}{3u^3} + \frac{ m \alpha v^2}{u^2} - \frac{\alpha v}{u},
 \end{gather*}
 and these can be substituted into our ansatz for $\omega^2$, $\omega^3$.
 Finally, our explicit parametrization for the coframe on $\Sigma_7$ satisfying the structure equations~\eqref{9dim-streqns} is given by~\eqref{9d-explicit-coframe}.


 \section{Calculation of the nine-dimensional symmetry algebras}
 \label{9d-sym-alg}

 We provide here the calculation of the (contact) symmetries of the maximally symmetric generic hyperbolic PDE's $F=0$ corresponding to the normal forms \eqref{special-maxsym-eqn} and \eqref{gen-eqn-alpha}.  A direct approach by computing one of:
 \begin{enumerate}\itemsep=0pt
 \item[1)] internal symmetries of $I_F$,
 \item[2)] external symmetries of $F$,
 \item[3)] symmetries of the lifted coframe on $\Sigma_7 \times H$,
 \end{enumerate}
 without any prior simplif\/ications is highly impractical owing to the complexity of the equation.  We f\/irst make several observations which allow us to simplify the calculation dramatically.

 Any equation of the form $F(r,s,t)=0$ admits a six-parameter family of symmetries
 \begin{gather*}
  \tilde{x} = c_1 + c_6 x, \qquad  \tilde{y} = c_2 + c_6 y, \qquad  \tilde{z} = c_3 + c_4 x + c_5 y + c_6{}^2 z
 \end{gather*}
 corresponding to the symmetry generators
 \begin{gather}
  \parder{x}, \qquad \parder{y}, \qquad  \parder{z}, \qquad
 x\parder{z}, \qquad  y\parder{z}, \qquad  x\parder{x} + y\parder{y} + 2z\parder{z}. \label{symgen6}
 \end{gather}
 These are clearly symmetries because the induced action on the second derivative coordina\-tes~$r$, $s$, $t$ is trivial.  (Equivalently, the prolongations of the vector f\/ields \eqref{symgen6} to $J^2(\R^2,\R)$ have no components along $\partial_r$, $\partial_s$, $\partial_t$.)

 We will treat both cases \eqref{special-maxsym-eqn} and \eqref{gen-eqn-alpha} simultaneously by using the parametric form of the equation.  Using \eqref{rst-param}, we have the following vector f\/ields tangent to the equation manifold in $J^2(\R^2,\R)$:
 \begin{gather*}
 \parder{w}  = -\epsilon m w^2 \parder{r} - \epsilon m^2 w \parder{s} - \epsilon m^3 \parder{t},\\
 \parder{v}  = -v^2 \parder{r} + v \parder{s} - \parder{t}.
 \end{gather*}
 The key is noticing the following $C^\infty(\Sigma_7)$-linear combinations of $\parder{w}$, $\parder{v}$ that can be expressed purely in terms of the 2-jet variables $x$, $y$, $z$, $p$, $q$, $r$, $s$, $t$:
 \begin{gather}
 w\parder{w} + v\parder{v}  = 3r\parder{r} + 2s\parder{s} + t\parder{t}, \label{hard-vf1}\\
  m\parder{w} - \parder{v}  = 2s\parder{r} + t\parder{s} + \alpha\parder{t}, \label{hard-vf2}\\
  (my+wx)\parder{w} - (y-vx)\parder{v}  = (2ys+3xr) \parder{r} + (yt+2xs) \parder{s} + (\alpha y+xt) \parder{t}. \label{hard-vf3}
 \end{gather}
 These three vector f\/ields are also tangent to the equation manifold, but they are clearly not contact vector f\/ields since: (1) they do not preserve the contact ideal, and (2) they do not arise as prolongations of vector f\/ields on $J^1(\R^2,\R)$ (c.f.\ B\"acklund's theorem).  However, this leads us to the following problem: Do \eqref{hard-vf1}--\eqref{hard-vf3} describe the $r$, $s$, $t$ components of contact vector f\/ields on $J^2(\R^2,\R)$?  If so, then those contact vector f\/ields would also be tangent to the equation manifold and hence would correspond to contact symmetries of the equation.  The search is greatly simplif\/ied by the observation that the components in  \eqref{hard-vf1}--\eqref{hard-vf3} are {\em linear} in $x$, $y$, $r$, $s$,~$t$ and independent of $z$, $p$, $q$.

 By B\"acklund's theorem, any contact vector f\/ield on $J^2(\R^2,\R)$ is the prolongation of a contact vector f\/ield on $J^1(\R^2,\R)$.  Hence, we look at a generalized vector f\/ield of order one on $J^0(\R^2,\R)$
 \begin{gather*}
 X  = \xi^1(x,y,z,p,q) \parder{x} + \xi^2(x,y,z,p,q) \parder{y} + \eta(x,y,z,p,q) \parder{z}
 \end{gather*}
 and examine its prolongation $X^{(2)} = pr^{(2)}(X)$ on $J^2(\R^2,\R)$.  If we write
 \begin{gather*}
 X^{(2)} = X + \eta^x \parder{p} + \eta^y \parder{q} + \eta^{xx} \parder{r} + \eta^{xy} \parder{s} + \eta^{yy} \parder{t},
 \end{gather*}
 then the standard prolongation formula \cite{Olver1995} is
 \begin{gather}
 \eta^{J,i} = D_i \eta^J - (D_i \xi^j) z_{J,j}, \label{prolongation}
 \end{gather}
 where $D_i$ are total derivative operators and we have used the notation $x^1= x$, $x^2=y$.  $J$~is an unordered (symmetric) multi-index, so that for example $z_1 = p$, $z_{12} = z_{21} = s$, etc.   Ian Anderson's {\tt JetCalculus} package in Maple~v.11 was very useful for computing and manipulating these prolongations.  We give an outline of the calculation here.

 In order that the components of $X^{(2)}$ in the directions $\parder{r}$, $\parder{s}$, $\parder{t}$, have: (1) no dependence on third-order terms, and (2) no {\em quadratic} dependence on $r$, $s$, $t$, we must have (using \eqref{prolongation}),
 \begin{gather}
 0=\eta_p - p (\xi^1)_p - q (\xi^2)_p, \label{der3-1}\\
 0=\eta_q - p (\xi^1)_q- q (\xi^2)_q, \label{der3-2}\\
 0=\eta_{pp} - p (\xi^1)_{pp} - q (\xi^2)_{pp} - 2 (\xi^1)_p, \label{der2-1}\\
 0=\eta_{qq} - p (\xi^1)_{qq} - q (\xi^2)_{qq} - 2 (\xi^2)_q, \label{der2-2}\\
 0=\eta_{pq} - p (\xi^1)_{pq} - q (\xi^2)_{pq} - (\xi^1)_q - (\xi^2)_p. \label{der2-3}
 \end{gather}
 Dif\/ferentiating \eqref{der3-1}, \eqref{der3-2} with respect to $p$ or $q$ and comparing with \eqref{der2-1}--\eqref{der2-3} leads to the conclusion that $\xi^1$, $\xi^2$, $\eta$ are all independent of $p$, $q$.  Consequently, the vector f\/ields we derive are necessarily inf\/initesimal {\em point} transformations, i.e.\ they project to vector f\/ields on $J^0(\R^2,\R)$.

 We f\/irst examine the coef\/f\/icient of $\parder{s}$ in $X^{(2)}$, or equivalently ${\cal L}_{X^{(2)}}(s)$.  Using \eqref{prolongation},
 \begin{gather}
  {\cal L}_{X^{(2)}} (s)  = -p^2 q (\xi^1)_{zz} - pq^2 (\xi^2)_{zz} - p^2 (\xi^1)_{yz} - 2ps (\xi^1)_z - pt(\xi^2)_z - q^2(\xi^2)_{xz} \nonumber\\
 \phantom{{\cal L}_{X^{(2)}} (s)  =}{}  - qr (\xi^1)_z - 2qs (\xi^2)_z + (\eta_{zz} - (\xi^1)_{xz} - (\xi^2)_{yz})pq - r(\xi^1)_y - t(\xi^2)_x \nonumber\\
 \phantom{{\cal L}_{X^{(2)}} (s)  =}{} + (\eta_{yz} - (\xi^1)_{xy})p + (\eta_{xz} - (\xi^2)_{xy})q + (\eta_z - (\xi^1)_x - (\xi^2)_y)s + \eta_{xy}. \label{Lie-s}
 \end{gather}
 Referring back to \eqref{hard-vf1}--\eqref{hard-vf3}, we require ${\cal L}_{X^{(2)}} (s) \in {\rm span}_\R\{ s,t,yt,xs \}$.  From \eqref{Lie-s}, we must have $\xi^1 = \xi^1(x)$, $\xi^2=\xi^2(x,y)$, and $\eta$ is linear in $z$ with $\eta_{yz}=\eta_{xy}=0$.  Thus, $\eta = C_1(x) + C_2(y) + C_3(x) z$.
Recalculating $X^{(2)}$, we have
 \begin{gather*}
 {\cal L}_{X^{(2)}} (r) = \eta_{xx} + (C_3(x) - 2(\xi^1)_x)r - 2s(\xi^2)_x + (2C_3'(x) - (\xi^1)_{xx})p - q(\xi^2)_{xx},  \\
  {\cal L}_{X^{(2)}} (s) =  - t(\xi^2)_x  + (C_3'(x) - (\xi^2)_{xy})q + (C_3(x) - (\xi^1)_x - (\xi^2)_y)s, \\
 {\cal L}_{X^{(2)}} (t) = C_2''(y) + (C_3(x) - 2(\xi^2)_y)t - (\xi^2)_{yy} q.
 \end{gather*}
 From \eqref{hard-vf1}--\eqref{hard-vf3}, we require ${\cal L}_{X^{(2)}} (r) \in {\rm span}_\R\{ r,s,xr,ys \}$ and ${\cal L}_{X^{(2)}} (t) \in {\rm span}_\R\{ 1,t,y,tx \}$.  Thus, $\eta_{xx}= (\xi^2)_{xx}=(\xi^2)_{yy}=0$, $C_2(y)$ is at most cubic in $y$, $(\xi^1)_x$ is linear in $x$, and $C_3'(x) = (\xi^2)_{xy}$.   Referring to \eqref{symgen6}, we have six symmetries with trivial action on $r,s,t$ and so without loss of generality we can require $\xi^1$ to have no constant or linear term, $\xi^2$ to have no constant term, and $\xi^3$ to have no constant term or terms purely linear in~$x$,~$y$.  Thus,
 \begin{gather*}
  \xi^1 = a_2 x^2,
\qquad \xi^2 = b_{01} x + b_{10} y + a_2 xy,
\qquad \eta =c_{020} y^2 + c_{030} y^3 + c_{001}z + a_2 xz,
 \end{gather*}
 and we have
 \begin{gather}
 {\cal L}_{X^{(2)}} (r)  = (-3a_2 x - 2a_1 + c_{001})r -2(b_{01} - a_2 y)s, \label{r-coeff}\\
 {\cal L}_{X^{(2)}} (s)  = (-2a_2 x - a_1 - b_{10} + c_{001})s + (b_{01} - a_2 y)t, \label{s-coeff}\\
 {\cal L}_{X^{(2)}} (t)  = 2 c_{020} + 6c_{030} y +(-x a_2 - 2 b_{10} + c_{001})t. \label{t-coeff}
 \end{gather}
 Comparing \eqref{r-coeff}--\eqref{t-coeff} to each of the desired expressions arising in \eqref{hard-vf1}--\eqref{hard-vf3}, we obtain the remaining symmetry generators
 \begin{gather*}
 X_7 = y\parder{y} + 3z\parder{z}, \qquad
 X_8 = x \parder{y} - \frac{\alpha}{2} y^2 \parder{z}, \qquad
 X_9= x^2 \parder{x} + xy \parder{y} + \left(xz - \frac{\alpha}{6} y^3 \right) \parder{z}.
 \end{gather*}

\subsection*{Acknowledgements}

 It is my pleasure to thank Niky Kamran for his lucid explanations of exterior dif\/ferential systems and the Cartan equivalence method, his guidance while studying \cite{GK1993}, and for bringing to my attention Vranceanu's work \cite{Vranceanu1937}.  Many of the calculations in this paper were either facilitated by or rechecked using the {\tt DifferentialGeometry}, {\tt LieAlgebras}, and {\tt JetCalculus} packages (in Maple v.11) written by Ian Anderson.  I would also like to thank Thomas Ivey and the three anonymous referees for their comments and corrections to help improve the exposition of this article.  This work was supported by funding from NSERC and McGill University.

\pdfbookmark[1]{References}{ref}
\LastPageEnding

\end{document}